\documentclass[preprint,3p,12pt]{elsarticle}
\usepackage{mathrsfs}
\usepackage{amsmath}
\usepackage{stmaryrd}
\usepackage{bbding}
\usepackage{dcolumn}
\usepackage{graphicx}
\usepackage{amsfonts}
\usepackage{amssymb}
\usepackage{psfrag}
\usepackage{wrapfig}
\usepackage{subfigure}
\usepackage{makeidx}
\usepackage{bm}
\usepackage{epsf}
\usepackage{epsfig}
\usepackage{setspace}
\usepackage{graphicx}
\usepackage{epstopdf}
\usepackage{psfrag}
\usepackage{subfigure}
\usepackage{color}

\begin{document}
\title{A compact high-order gas-kinetic scheme on unstructured mesh for acoustic and shock wave computations}

\author[HKUST1]{Fengxiang Zhao}
\ead{fzhaoac@connect.ust.hk}

\author[HKUST2]{Xing Ji}
\ead{xjiad@connect.ust.hk}

\author[HKUST1]{Wei Shyy}
\ead{weishyy@ust.hk}

\author[HKUST1,HKUST2,HKUST3]{Kun Xu\corref{cor}}
\ead{makxu@ust.hk}

\address[HKUST1]{Department of Mechanical and Aerospace Engineering, Hong Kong University of Science and Technology, Clear Water Bay, Kowloon, HongKong}
\address[HKUST2]{Department of Mathematics, Hong Kong University of Science and Technology, Clear Water Bay, Kowloon, HongKong}
\address[HKUST3]{Shenzhen Research Institute, Hong Kong University of Science and Technology, Shenzhen, China}
\cortext[cor]{Corresponding author}

\begin{abstract}
	
Following the development of a third-order compact gas-kinetic scheme (GKS) for the Euler and Navier-Stokes equations (Journal of Computational Physics 410 (2020) 109367),
in this paper an even higher-order compact GKS up to sixth order of accuracy will be constructed for the shock and acoustic wave computation on unstructured mesh.
The compactness is defined by the physical domain of dependence for an unstructured triangular cell, which may involve the closest neighbors of
neighboring cells.
The compactness and high-order accuracy of the scheme are coming from the consistency between the high-order initial reconstruction and the
high-order gas evolution model under GKS framework.
The high-order evolution solution at a cell interface provides not only a time-accurate flux function, but also the time-evolving flow variables. Therefore, the cell-averaged flow variables and their gradients can be explicitly updated at the next time level from the  moments of the same time-dependent gas distribution function.
Based on the cell averages and cell-averaged derivatives, both linear and nonlinear high-order reconstruction can be obtained for macroscopic flow variables in the evaluation of local equilibrium and non-equilibrium states. The current nonlinear reconstruction is a combination of WENO and ENO methodology, which is specifically suitable for compact GKS on unstructured mesh with a high-order ($\geq 4$) accuracy.
The initial piecewise discontinuous reconstruction is used for the determination non-equilibrium state and an evolved smooth reconstruction for the equilibrium state. The evolution model in gas-kinetic scheme is based on a relaxation process from non-equilibrium to equilibrium state.
The time-accurate gas distribution function in GKS provides the Navier-Stokes flux function directly without separating the inviscid and viscous terms, which simplifies the numerical method on unstructured mesh.
Based on the time-accurate flux solver, the two-stage fourth-order time discretisation can be
applied to get a fourth-order time-accurate solution with only two stages, which reduces two reconstructions in comparison with the same time-accurate method with Runge-Kutta time stepping.
The current high-order GKS can uniformly capture acoustic and shock waves without identifying trouble cells and implementing additional limiting procedure. In addition, the fourth- up to sixth-order compact GKS can use almost the same time step as a second-order shock capturing scheme. The fourth-order GKS on unstructured mesh will be used in the computations from low speed incompressible viscous flow to the high Mach number shock interaction.
The accuracy, efficiency, and robustness of the scheme have been validated. The main conclusion of the paper is that beyond the first-order Riemann solver, the use of high-order gas evolution model seems necessary in the development of high-order schemes.
\end{abstract}

\begin{keyword}
Gas-kinetic model, High-order reconstruction, Compact scheme, Triangular mesh
\end{keyword}

\maketitle
\section{Introduction}

The development of high-order methods is important for the compressible flow computations in many engineering applications,
such as compressible turbulent flows \cite{hill_hybrid}, aero-acoustics \cite{tam1995AIAA}, and complex flows with shock and boundary layer interactions \cite{garnier-shock-boundlayer}.
With the advantages of high-order methods in terms of solution accuracy and computational efficiency \cite{WangZJ-high-order}, extensive effort has been spent on search of high-order schemes in the past decades \cite{Case-Woodward,harten2,liu,jiang,shu,reed,cockburn1,CPR,dumbser,pan1}.
The incompatibility among the compactness, high-order accuracy, and shock capturing still remains in the algorithm development.
The discrepancy between the linear scheme for smooth flow and nonlinear scheme for discontinuous one has not been solved uniformly.
For the flow with co-existing shock interaction and acoustic wave propagation, the unification of linear and nonlinear formulation
is a preferred property in a high-order scheme.

The high-order compact scheme is attractive due to its high resolution, high parallel efficiency, low storage, and simple implementation on unstructured mesh \cite{lele,CGKSAIA,shu,CPR,dumbser}. Compact scheme has  physical foundation because the CFL condition determines the domain of dependence, where only neighboring cells are dynamically connected. Theoretically, all high-order non-compact schemes have weakness by including the cells without dynamic dependence into the evolution of local solution.
The popular Lele-type compact schemes are based on the implicit relation between flow variables and their derivatives within a compact stencil \cite{lele}, which influences greatly on the development of linear schemes \cite{WangQ,bai}.
For smooth flows, the successful compact schemes include the finite difference (FD) and finite volume (FV) schemes \cite{lele,mahesh,WangQ}, DG method \cite{reed,cockburn1}, and correction procedure via reconstruction (CPR) method \cite{FR,CPR}.
The DG method has the same stencil as the second-order scheme, and achieves high accuracy by using high-order piecewise polynomials within cells and evolving the multiple degrees of freedom (DOFs) based on their distinguishable governing equations from the weak formulations. Each cell in DG method only interacts with its neighboring cells and the method become very efficient for parallel computation \cite{shu} and the application on unstructured meshes. For discontinuous flows, the success of the above compact schemes is limited.
The research on identifying the  trouble cells and using additional limiting procedure seems unavoidable in DG \cite{cockburn2,qiu1,DG_limiters} and CPR \cite{CPR_limiter-1,CPR_limiter-2} methods.

Nonlinear schemes have been designed for the flows with discontinuities. The successful nonlinear schemes include total variation diminishing (TVD) schemes \cite{harten1}, essentially non-oscillation (ENO) schemes \cite{harten2}, and weighted ENO (WENO) schemes \cite{liu,jiang}. In past twenty years, the WENO-type high-order methods have received the most attention among nonlinear schemes. The central ingredients in WENO schemes are to construct several low order polynomials and to design smoothness indicators to adaptively assemble them to get a high-order one.
Most recent effort is about the optimization of stencil selection and the design of weighting functions  \cite{taylor_WENO-optimize,WENO-M,WENO-Z,zhao2017_wneo6th}. WENO schemes can achieve very high-order accuracy in the smooth region and maintain non-oscillatory property across shock waves \cite{shu1}. But, as mentioned before, the large stencils used in WENO use the information of the cells which have no any dynamic connection under the CFL condition with the local reconstructed cell. As a result, the numerical dissipation in WENO is still high in comparison with compact schemes \cite{taylor_WENO-optimize}. Modified WENO schemes, such as WENO-M and WENO-Z \cite{WENO-M,WENO-Z}, have been proposed to improve the performance at critical points. The hybrid linear and nonlinear schemes have been investigated as well \cite{hill_hybrid,ren}.

Both non-compact and compact high-order GKS have been developed in the past years on structured and unstructured meshes \cite{pan1,pan3,ji1,ji2}.
On unstructured mesh, only third-order accurate compact GKS has been developed \cite{pan4,ji4}.
However, on structured uniform mesh, compact GKS up to eighth-order accuracy have been successfully constructed for compressible flow simulation with shock and aero-acoustics waves \cite{CGKSAIA,Compact-GKSacoustic}.
The high-order compact GKS demonstrates high-order accuracy, spectral-like resolution and excellent robustness for both smooth and discontinuous solutions under a large CFL number ($CFL\geq0.5$).
In this paper, the compact high-order GKS on  unstructured mesh will be further developed with even higher-order accuracy.
Here the compactness is related to the domain of dependence for physical wave propagation under the CFL condition, which involves only nine cells around the central triangle, see Fig.\ref{0-stencil}. Based on such a stencil, up to sixth-order compact scheme will be developed.
For the non-compact WENO schemes with the same accuracy on unstructured mesh, a much large stencil with different order of magnitude has to be used.
The high-order compact GKS can capture smooth and discontinuous solutions with the update of the numerical solution without identifying the trouble cells and including additional limiting procedure.
Due to the time-accurate evolution solution at a cell interface, both cell-averaged conservative flow variables and their gradients can be explicitly updated \cite{xu,ji2}. The cell-averaged derivatives are coming from the strong evolution solution rather than the weak formulation in DG method for the update of the similar DOFs \cite{cockburn1,cockburn2}.

Besides gas evolution model, the high-order linear and nonlinear reconstructions are also involved in the GKS. Based on the cell-averaged flow variables and their derivatives, the methodology of WENO method can be extended here on unstructured mesh for the nonlinear reconstruction to deal with discontinuous flow, even though it has difficulty to get high-order reconstruction in the traditional WENO methods with the cell-averaged flow variables alone \cite{hu1999WENO,zhao2018WENO}. The nonlinear WENO reconstruction will be combined with a linear reconstruction obtained through a constrained least square (CLS) method to get an initial piecewise discontinuous reconstruction.
 The nonlinear reconstruction will be for the determination of non-equilibrium state around a cell interface.
A smooth reconstruction for the equilibrium state there is obtained from the nonlinear reconstructions through a kinetic collision model.
Then, the time accurate gas distribution function at a cell interface is constructed from a relaxation process from the initial non-equilibrium state to an equilibrium one, and it provides Navier-Stokes solution for the flow variable and flux evaluation.
The convergence study shows that the compact GKS can work well from fourth to sixth order of accuracy.
Based on the time-accurate flux solver, another advantage of the current method is to implement the two-stage fourth-order (S2O4) time discretisation, which is developed for the Lax-Wendroff type flow solvers \cite{li,du-hermite}, to improve the temporal accuracy with less middle stages.
The S2O4 method has been validated in the previous non-compact and compact GKS \cite{pan1,pan2,ji2,ji4}.
Due to the use of only two reconstructions in the fourth-order scheme, the GKS becomes efficient in comparison with the schemes based on the Runge-Kutta time stepping method.
Besides, the multi-stage multi-derivative time integrators for hyperbolic conservation laws \cite{seal} can be adopted into GKS for higher-order temporal accuracy \cite{ji1}.

This paper is organized as follows. The GKS will be introduced in Section 2. Section 3 is about the updates of cell averaged flow variables and
cell-averaged derivatives through S2O4 method. Section 4 and Section 5 will present the compact linear and nonlinear reconstructions for
the determination  of a final piecewise discontinuous high-order polynomial inside each control volume.
In Section 6, the  GKS will be validated in a wide range of test cases from the strong shock interaction to acoustic wave propagation.
The last section is the conclusion.

\section{Gas-kinetic schemes}
The GKS provides a time-accurate gas evolution model from a piecewise discontinuous initial polynomials \cite{xu2,xu1}.
Starting from high-order spatial reconstruction, various temporal discretizations have been used in GKS \cite{liQB2010high,pan0,ji1,ji2,pan4}.
The successful applications include compressible multi-component flow \cite{pan2}, DNS at high Mach number \cite{cao2019}, and turbulence simulation \cite{cao2019implicit}.
A brief introduction of GKS and the special features will be presented in this section.

The gas-kinetic evolution model in GKS is based on the BGK equation \cite{BGK-1},
\begin{equation}\label{bgk}
f_t+\textbf{u}\cdot\nabla f=\frac{g-f}{\tau},
\end{equation}
where $\textbf{u}=(u,v)$ is the particle velocity, $f$ is the gas distribution function, $g$ is the corresponding equilibrium state that $f$ approaches, and $\tau$ is particle collision time.
The equilibrium state $g$ is a Maxwellian distribution,
\begin{equation*}
\begin{split}
g=\rho(\frac{\lambda}{\pi})^{\frac{K+2}{2}}e^{-\lambda((u-U)^2+(v-V)^2+\xi^2)},
\end{split}
\end{equation*}
where $\lambda =m/2kT $, and $m, k, T$ are the molecular mass, the Boltzmann constant, and temperature, respectively.
$K$ is the number of internal degrees of freedom, i.e. $K=(4-2\gamma)/(\gamma-1)$ for two-dimensional flow,
and $\gamma$ is the specific heat ratio. $\xi$ is the internal variable with $\xi^2=\xi^2_1+\xi^2_2+...+\xi^2_K$.
Due to the conservation of mass, momentum and energy during particle collisions, $f$ and $g$ satisfy the compatibility condition,
\begin{equation}\label{compatibility}
\int \frac{g-f}{\tau}\pmb{\psi} \mathrm{d}\Xi=0,
\end{equation}
at any point in space and time, where $\pmb{\psi}=(\psi_1,\psi_2,\psi_3,\psi_4)^T=(1,u,v,\displaystyle \frac{1}{2}(u^2+v^2+\xi^2))^T$, $\text{d}\Xi=\text{d}u\text{d}v\text{d}\xi_1...\text{d}\xi_{K}$.

The macroscopic mass $\rho$, momentum ($\rho U, \rho V$), and energy $\rho E$ can be evaluated from the gas distribution function,
\begin{equation}\label{g-to-convar}
{\textbf{W}} =
\left(
\begin{array}{c}
\rho\\
\rho U\\
\rho V\\
\rho E\\
\end{array}
\right)
=\int f \pmb{\psi} \mathrm{d}\Xi.
\end{equation}
The corresponding fluxes for mass, momentum, and energy in $i$-th direction is given by
\begin{equation}\label{g-to-flux}
{\textbf{F}_i} =\int u_i f \pmb{\psi} \mathrm{d}\Xi,
\end{equation}
with $u_1 = u$ and $u_2 = v$ in the 2D case.

Based on the BGK equation, the GKS provides a time-accurate evolution solution $f$ at a cell interface \cite{xu1}.
On the mesh size scale, the conservations of mass, momentum and energy in a control volume become
\begin{equation}\label{semifvs}
\frac{\text{d}\textbf{W}_{j}}{\text{d}t}=-\frac{1}{\big|\Omega_j\big|}\int_{\partial \Omega_j} \textbf{F}\cdot \textbf{n} \mathrm{d} l,
\end{equation}
where $\textbf{W}_{j}$ is the cell-averaged conservative variables, $\textbf{F}=(\textbf{F}_1,\textbf{F}_2)$ is the time dependent fluxes at cell interfaces. The $\textbf{W}_{j}$ is defined as
\begin{align}\label{cell-average}
\textbf{W}_{j}&\equiv \frac{1}{\big| \Omega_j \big|} \iint_{\Omega_j} \textbf{W}(x,y) \text{d}x\text{d}y.
\end{align}
The line integral on the right hand side (RHS) of Eq.(\ref{semifvs}) is discretized by a q-point Gaussian integration formula,
\begin{align}\label{semifvs-rhs}
\int_{\partial \Omega_j} \textbf{F}\cdot \textbf{n} \mathrm{d} l=\sum_{l=1}^{l_0}\big( \big|\Gamma_{l} \big| \sum _{k=1}^q \omega_k \textbf{F}(x_k,y_k)\cdot \textbf{n}_l \big),
\end{align}
where $l_0$ is the side number of cell $\Omega_j$, and $l_0=3$ for a triangular mesh.
Due to the connection among the flow variables $\textbf{W}$, the fluxes $\textbf{F}$ , and the gas distribution function $f$,
the central point of GKS is to construct a time-dependent gas distribution function $f$ at the cell interface.
The integral solution  of BGK equation is \cite{xu2},
\begin{equation}\label{integral1}
\begin{split}
f(\textbf{x}_0,t,\textbf{u},\xi)=&\frac{1}{\tau}\int_0^t g(\textbf{x}^{'},t',\textbf{u},\xi)e^{-(t-t')/\tau}\mathrm{d}t' \\
&+e^{-t/\tau}f_0(\textbf{x}_0-\textbf{u}(t-t_0),u,v,\xi),
\end{split}
\end{equation}
where $\textbf{x}_0$ is the numerical quadrature point at the cell interface for flux evaluation, and $\textbf{x}_0=\textbf{x}^{'}+\textbf{u}(t-t^{'})$ is the particle trajectory.
Here $f_0$ is the initial state of gas distribution function $f$ at $t=0$.
The integral solution basically presents a physical process from the particle free transport in $f_0$ in the kinetic scale to the hydrodynamic flow evolution in the integration of $g$.
The contributions from $f_0$ and $g$ in the determination of $f$ at the cell interface depend on the ratio of time step to the local particle collision time, i.e., $e^{-t/\tau}$.
For the NS solution, based on the Chapman-Enskog expansion $f_0$ can be explicitly determined from initial reconstructions of macroscopic flow variables.
In the current high-order GKS, the high-order nonlinear reconstruction with compact stencils will be used in the determination of $f_0$.
And based on the nonlinear reconstruction on both sides of the cell interface, the smooth reconstruction is obtained dynamically with the same high-order accuracy for the equilibrium $g$.
Therefore, the above integral solution not only incorporates a physical evolution process from initial discontinuous non-equilibrium state
to a continuous equilibrium one, but also unifies the discontinuous and smooth reconstruction in the evolution process.
This fact is crucially important for the scheme to capture both nonlinear shock and linear acoustic wave accurately
in the computation with its dynamic adaptation factor of $\exp{(-t/\tau)}$.

In order to obtain the solution $f$, both $f_0$ and $g$ in Eq.(\ref{integral1}) need to be modeled \cite{xu2,liQB2010high}.
Based on the integral solution, the gas distribution function with a second-order accuracy is \cite{xu2}
\begin{align}\label{2nd-f}
f(\textbf{x}_0,t,\textbf{u},\xi)=&(1-e^{-t/\tau})g_0+((t+\tau)e^{-t/\tau}-\tau)(\overline{a}_1u+\overline{a}_2v)g_0\nonumber\\
+&(t-\tau+\tau e^{-t/\tau}){\bar{A}} g_0\nonumber\\
+&e^{-t/\tau}g_r[1-(\tau+t)(a_{1r}u+a_{2r}v)-\tau A_r)]H(u)\nonumber\\
+&e^{-t/\tau}g_l[1-(\tau+t)(a_{1l}u+a_{2l}v)-\tau A_l)](1-H(u)),
\end{align}
where the terms related to $g_0$ are from the integral of the equilibrium state and the terms related to $g_l$ and $g_r$ are from the initial term $f_0$ in the Eq.(\ref{integral1}). All the coefficients in Eq.(\ref{2nd-f}) can be determined from the initially reconstructed  macroscopic flow variables.

\section{Update of cell-averaged variables and derivatives}
In this section, the updates of cell-averaged conservative variables and their derivatives are presented.
Based on the time-accurate gas distribution function in Eq.(\ref{2nd-f}), both fluxes and conservative flow variables at the cell
interface can be evaluated. Besides the updates of conservative flow variables through fluxes,
the cell-averaged derivatives can be calculated as well through Gauss's theorem with the integration of flow variables around a closed
boundary of the target cell.
The original complicated time relaxation process in gas distribution function can be approximated as a
linear function of time and is used in the MSMD framework in time discretization.

\subsection{Time-accurate gas distribution function and temporal discretisation}
For the second-order evolution model, the time accurate $f$ in Eq.(\ref{2nd-f}) at cell interface can be approximated as a linear function of time,
\begin{equation}\label{f-linearize-1}
\hat{f}(t)=f^n+ t f_t^n, \\
\end{equation}
where the unknowns $f^n$ and $f_t^n$ can be determined by two conditions from Eq.(\ref{2nd-f}). Define the time integration of $f(t)$ as
\begin{align*}
\overline{f}(t)=\int_0^tf(t')dt'.
\end{align*}
By using $\overline{f}(t/2)$ and $\overline{f}(t)$, we can get the unknowns in Eq.(\ref{f-linearize-1})
\begin{align}\label{f-linearize-2}
\begin{split}
&f^n=\frac{1}{t}\big( 4\overline{f}(t/2)-\overline{f}(t) \big),\\
&f^n_t=\frac{4}{t^2}\big(-2\overline{f}(t/2)+\overline{f}(t) \big).
\end{split}
\end{align}
Taking moments $u\pmb\psi$ on $\hat{f}(t)$ and $\hat{f}_t(t)$ at $t=0$, the numerical flux and its time derivative can be obtained,
\begin{align}\label{f-F-numer}
\begin{split}
&\textbf{F}^n=\int u f^n \pmb{\psi} \mathrm{d}\Xi,\\
&\textbf{F}^n_t=\int u f^n_t \pmb{\psi} \mathrm{d}\Xi.
\end{split}
\end{align}
Both $\textbf{F}^n$ and $\textbf{F}^n_t$ will be used to update the cell averages of flow variables through MSMD method.

The S2O4 method \cite{li,pan1,ji2} is adopted in the compact GKS as temporal discretisation to achieve a fourth-order accuracy.
It can be extended to higher-order temporal accuracy with the multi-stage multi-derivative time integrator \cite{hairer,seal} under the Lax-Wendroff type solver-based finite volume framework \cite{li,pan1}.
The details can be found in \cite{ji1}. In this paper the fourth-order time discretization will be adopted.
For conservation laws, the semi-discrete finite volume scheme Eq.(\ref{semifvs}) is rewritten as
\begin{align*}
\frac{\text{d} \textbf{W}_{j}}{\text{d} t}\equiv \mathcal{L}_j(\textbf{W}),
\end{align*}
where $\mathcal{L}_j(\textbf{W})$ is the total flux around the control volume $\Omega_j$.
A solution $\textbf{W}(t)$ at $t=t_n+\Delta t$ with fourth-order time accuracy can be obtained by
\begin{equation}\label{step-hyper-1}
\begin{split}
\textbf{W}^{n+1/2}_j&=\textbf{W}^n_j+\frac{1}{2}\Delta t\mathcal
{L}_j(\textbf{W}^n)+\frac{1}{8}\Delta t^2\frac{\partial}{\partial
	t}\mathcal{L}_j(\textbf{W}^n),\\
\textbf{W}^{n+1}_j&=\textbf{W}^n_j+\Delta t\mathcal
{L}_j(\textbf{W}^n)+\frac{1}{6}\Delta t^2\big(\frac{\partial}{\partial
	t}\mathcal{L}_j(\textbf{W}^n)+2\frac{\partial}{\partial
	t}\mathcal{L}_j(\textbf{W}^{n+1/2})\big),
\end{split}
\end{equation}
where $\mathcal{L}_j$ and $\partial \mathcal{L}_j/\partial t$ are related to the fluxes and the time derivatives of the fluxes which are given in Eq.(\ref{f-F-numer}). The middle state $\textbf{W}_j^{n+1/2}$ is obtained at time $t^{n+1/2} = t^n + \Delta t /2$.

\subsection{Update of the cell-averaged derivatives}
Taking moments $\pmb\psi$ on the same $\hat{f}(t)$ and $\hat{f}_t(t)$ in Eq.(\ref{f-linearize-1}) at $t=0$, the flow variables and their time derivatives at cell interface can be obtained,
\begin{align}\label{f-W}
\begin{split}
&\textbf{W}^n(\textbf{x})=\int f^n \pmb{\psi} \mathrm{d}\Xi,\\
&\textbf{W}^n_{t}(\textbf{x})=\int f^n_t \pmb{\psi} \mathrm{d}\Xi.
\end{split}
\end{align}
Similar to the updates of the cell averages, the update of flow variables at cell interface has a middle state $\textbf{W}^{n+1/2}(\textbf{x})$ at $t^{n+1/2}=t^n + \Delta t /2$ obtained by
\begin{equation}\label{S2O3-1}
\textbf{W}^{n+1/2}(\textbf{x})=\textbf{W}^n(\textbf{x})+\frac{1}{2}\Delta t \textbf{W}_{t}^n(\textbf{x}),
\end{equation}
which gives a second-order approximation to $\textbf{W}(\textbf{x},t^{n+1/2})$ as
\begin{align*}
\textbf{W}^{n+1/2}(\textbf{x}) -\textbf{W}(\textbf{x},t^{n+1/2})= -\frac{1}{8} \Delta t^2 \textbf{W}_{tt}(\textbf{x},t^{n}) + O(\Delta t^3).
\end{align*}
In order to get $\textbf{W}^{n+1}(\textbf{x})$ with a third order of accuracy, a two-stage third-order (S2O3) time marching method can be adopted, and $\textbf{W}^{n+1}(\textbf{x})$ can be obtained as
\begin{equation}\label{S2O3-2}
\textbf{W}^{n+1}(\textbf{x})=\textbf{W}^n(\textbf{x})+\Delta t \textbf{W}_{t}^{n+1/2}(\textbf{x}),
\end{equation}
where $\textbf{W}_{t}^{n+1/2}(\textbf{x})$ is given in the same way as $\textbf{W}_{t}^{n}(\textbf{x})$ in the stage from $t^{n+1/2}$ to $t^{n+1}$.
Thus, the conservative variables $\textbf{W}^{n+1}(\textbf{x})$ at cell interface are obtained with a third-order accuracy from the same distribution function for fluxes in Eq.(\ref{f-linearize-1}) through Eq.(\ref{S2O3-1}) and Eq.(\ref{S2O3-2}), and it has third-order accuracy
\begin{equation}\label{S2O3-2-1}
\textbf{W}^{n+1}(\textbf{x}) -\textbf{W}(\textbf{x},t^{n+1})= -\frac{1}{6} \Delta t^3 \textbf{W}_{ttt}(\textbf{x},t^{n}) + O(\Delta t^4).
\end{equation}
The cell-averaged derivatives of flow variables are calculated directly from $\textbf{W}^{n+1}(\textbf{x})$ by Gauss's theorem.
\begin{align}\label{Gauss-theorem-1}
\begin{split}
\textbf{W}_{j,x}&\equiv \frac{1}{\big| \Omega_j \big|}\iint_{\Omega_j} \frac{\partial \textbf{W}(x,y)}{\partial x} dxdy=\frac{1}{\big| \Omega_j \big|}\int_{\partial \Omega_{j}} \textbf{W}(x,y)dy,\\
\textbf{W}_{j,y}&\equiv \frac{1}{\big| \Omega_j \big|}\iint_{\Omega_j} \frac{\partial \textbf{W}(x,y)}{\partial y} dxdy=-\frac{1}{\big| \Omega_j \big|}\int_{\partial \Omega_{j}} \textbf{W}(x,y)dx.
\end{split}
\end{align}
The above integration can be discretized by Gaussian quadrature, where the same quadrature points are used to calculate the numerical fluxes at the interface in Eq.(\ref{semifvs-rhs}). Thus the cell-averaged derivatives at $t^{n+1}$ can be updated as
\begin{align}\label{Gauss-theorem-2}
\begin{split}
\textbf{W}^{n+1}_{j,x}&=\frac{1}{\big| \Omega_j \big|}\sum_{l=1}^{l_0} \big(|\Gamma_l|\cdot n_{l,x} \sum _{k=1}^q \omega_k\textbf{W}^{n+1}(\textbf{x}_k) \big),\\
\textbf{W}^{n+1}_{j,y}&=\frac{1}{\big| \Omega_j \big|}\sum_{l=1}^{l_0} \big(|\Gamma_l|\cdot n_{l,y} \sum _{k=1}^q \omega_k\textbf{W}^{n+1}(\textbf{x}_k) \big).
\end{split}
\end{align}

Next, it will be demonstrated that the cell-averaged derivatives in Eq.(\ref{Gauss-theorem-2}) have a third order of accuracy. Without loss of generality, the cell-averaged x derivative will be used to give the proof. Substituting Eq.(\ref{S2O3-2-1}) into the first equation of Eq.(\ref{Gauss-theorem-2}), we have
\begin{align}\label{der-average-accuracy}
\begin{split}
\textbf{W}^{n+1}_{j,x}&=\frac{1}{\big| \Omega_j \big|}\sum_{l=1}^{l_0} \big(|\Gamma_l|\cdot n_{l,x} \sum _{k=1}^q \omega_k \big( \textbf{W}(\textbf{x}_k,t^{n+1}) - \frac{\Delta t^3}{6} \textbf{W}_{ttt}(\textbf{x}_k,t^{n}) + O(\Delta t^4)\big) \big)\\
                &=\frac{1}{\big| \Omega_j \big|}\sum_{l=1}^{l_0} \big(|\Gamma_l|\cdot n_{l,x} \sum _{k=1}^q \omega_k \textbf{W}(\textbf{x}_k,t^{n+1}) \big) \\
                &-\frac{1}{\big| \Omega_j \big|}\sum_{l=1}^{l_0} \big(\frac{\Delta t^3}{6}|\Gamma_l|\cdot n_{l,x} \sum _{k=1}^q \omega_k \textbf{W}_{ttt}(\textbf{x}_k,t^{n}) \big) + O(\Delta t^4/h),
\end{split}
\end{align}
where $h$ is a characteristic quantity representing the mesh size, and the magnitude of $h$ is $h \sim O\big( \big| \Omega_j \big|^{1/2} \big) \sim O( |\Gamma_l| )$. Due to the closed cell interfaces of a control volume, the following identity holds
\begin{align}\label{polygon-equation-side}
\sum_{l=1}^{l_0} |\Gamma_l|\cdot n_{l,x} =0.
\end{align}
Then the last two terms on  RHS of Eq.(\ref{der-average-accuracy}) becomes
\begin{equation}\label{der-average-accuracy-1}
\begin{split}
T.E.(\textbf{W}_{j,x}) &\equiv -\frac{1}{\big| \Omega_j \big|}\sum_{l=1}^{l_0} \frac{\Delta t^3}{6}|\Gamma_l|\cdot n_{l,x} \big( \sum _{k=1}^q \omega_k \textbf{W}_{ttt}(\textbf x_{l,k},t^{n}) \big) + O(\Delta t^4/h)  \\
&=-\frac{1}{\big| \Omega_j \big|}\sum_{l=2}^{l_0} \frac{\Delta t^3}{6} |\Gamma_l|\cdot n_{l,x} \big( \sum _{k=1}^q \omega_k  \big( \textbf{W}_{ttt}(\textbf x_{l,k},t^{n})- \textbf{W}_{ttt}(\textbf x_{1,k},t^{n}) \big) \big) + O(\Delta t^4/h) .
\end{split}
\end{equation}
Suppose $\textbf{W}_{ttt}(\textbf x,t^{n+1})$ in Eq.(\ref{der-average-accuracy-1}) is locally bounded and Lipschitz continuous in the domain $\{\textbf x| \textbf x \in \Omega_j+\partial \Omega_j \}$, i.e.,
\begin{align}\label{W-Lipschitz-continuous}
\big| \textbf{W}_{ttt}(\textbf x_{l,q},t^{n+1})- \textbf{W}_{ttt}(\textbf x_{1,q},t^{n+1}) \big| \leq L \cdot h,
\end{align}
where L is a constant.
Substituting Eq.(\ref{W-Lipschitz-continuous}) into Eq.(\ref{der-average-accuracy-1}), the truncation error becomes
\begin{equation}\label{der-average-accuracy-2}
\begin{split}
T.E.(\textbf{W}_{j,x}) &=-\frac{1}{\big| \Omega_j \big|}\sum_{l=2}^{l_0} \frac{\Delta t^3}{6}|\Gamma_l|\cdot n_{l,x} \big(\sum _{k=1}^q \omega_k  L \cdot h \big) + O(\Delta t^4/h ) \\
&=O( \Delta t^3 + \Delta t^4/h).
\end{split}
\end{equation}
Thus the cell-averaged derivative in Eq.(\ref{Gauss-theorem-2}) approximates the exact one with a third order of accuracy. The cell-averaged derivative is used to get the compact reconstruction in the next section.

\section{Compact linear reconstruction}
The GKS needs reconstruction of macroscopic variables for both initial non-equilibrium state and evolved equilibrium one.
The linear reconstruction is the foundation of the nonlinear one.
In this section, the linear reconstruction will be firstly introduced.

In this paper only the triangular mesh is considered. The reconstruction will be conducted in the conforming space.
In GKS, besides the cell average $\textbf{W}_j$, the cell-averaged derivatives $\textbf{W}_{j,x}$ and $\textbf{W}_{j,y}$ of conservative variables are available in each cell, as presented in the Section 3. The components (i.e. components of conservative or characteristic variables) of cell averages and their derivatives are denoted as $Q_j$, $Q_{j,x}$ and $Q_{j,y}$.
To achieve a high-order spatial discretisation, the high-order polynomial will be reconstructed from the given cell averages and cell-averaged derivatives.
With the consideration of possible large aspect ratio of the edges of the triangle, a reference coordinate system $\xi-\eta$ is introduced first to do the reconstruction.  For the cell $\Omega_0$, the transformation from the physical coordinate system $x-y$ into the reference coordinate system $\xi-\eta$ is defined as
\begin{align}\label{refer-physical}
\left(
\begin{array}{c}
x\\
y\\
\end{array}
\right)
=\left(
\begin{array}{c}
x_1\\
y_1\\
\end{array}
\right) +
\textbf{J}
\left(
\begin{array}{c}
\xi\\
\eta\\
\end{array}
\right),
\end{align}
where $\textbf{J}$ is the Jacobian matrix, and $(x_k,y_k), k=1,2,3$ is the nodes' coordinates of $\Omega_0$ in the physical coordinate system $x-y$. With the transformation, the cell $\Omega_0$ is transformed to $\widetilde{\Omega}_0$. The Jacobian matrix takes
\begin{equation}
\textbf{J}={
	\left( \begin{array}{cc}
	x_2-x_1 & x_3-x_1 \\
	y_2-y_1 & y_3-y_1
	\end{array}
	\right ).}
\end{equation}

The reconstruction polynomial for cell $\widetilde{\Omega}_0$ can be expanded over polynomial basis functions expressed as $\varphi_k (\bm{\xi})$ in the reference coordinate system. The polynomial and basis functions are given by
\begin{align}\label{Recons-poly}
\begin{split}
&P^r(\bm{\xi})=\sum_{k=0}^{N} a_k\cdot \varphi_k(\bm{\xi}),\\
&\varphi_k(\bm{\xi})=\frac{1}{k!}\xi^{k_1}\eta^{k_2},~~ k_1+k_2=k, k_1,k_2=0,\cdots,k,
\end{split}
\end{align}
where $\bm{\xi}=(\xi,\eta)$ are the coordinates in the reference coordinate system. $a_k$ are the DOFs in the summation of expansion. $N$ is related to the order of the polynomial by $N=\frac{1}{2}(r+1)(r+2)$ for the current two-dimensional case.
To fully determine the reconstructed polynomials, the given averages in all cells of a stencil $\widetilde{S}_0$ should be recovered, i.e.,
\begin{align}\label{Recons-eqs}
\begin{split}
&\big(\frac{\big|\textbf{J}\big|}{\big|\Omega_l \big|} \int_{ \widetilde{\Omega}_l} \varphi_k(\bm{\xi}) d\xi d\eta \big) a_k=Q_l, \\
&\big(\frac{\big|\textbf{J}\big|}{\big|\Omega_l \big|} \int_{ \widetilde{\Omega}_l} \varphi_{k,\xi}(\bm{\xi}) d\xi d\eta \big) a_k=Q_{l,\xi}, \\
&\big(\frac{\big|\textbf{J}\big|}{\big|\Omega_l \big|} \int_{ \widetilde{\Omega}_l} \varphi_{k,\eta}(\bm{\xi}) d\xi d\eta \big) a_k=Q_{l,\eta}, \widetilde{\Omega}_l \in \widetilde{S}_0,
\end{split}
\end{align}
where $\textbf{J}$ is the determinant of the Jacobian matrix $\textbf{J}$, and $\big|\Omega_l \big|$ is the area of cell $\Omega_l$. There is $\big|\Omega_l \big|/\textbf{J}=\big|\widetilde{\Omega}_l \big|$ for $l=0$.
The $Q_{l,\xi}$ and $Q_{l,\eta}$ are the cell-averaged derivatives in $\widetilde{\Omega}_l$ in the reference coordinate system, and both  derivatives are obtained from $Q_{l,x}$ and $Q_{l,y}$,
\begin{align}\label{Der-relation}
\left(
\begin{array}{c}
Q_{\xi}\\
Q_{\eta}\\
\end{array}
\right)
=\textbf{J}^T
\left(
\begin{array}{c}
Q_{x}\\
Q_{y}\\
\end{array}
\right).
\end{align}
Based on Eq.(\ref{Recons-eqs}), the error in the reconstruction from the cell-averaged derivatives becomes $O(h)\cdot T.E.(Q_{0,x}) \sim O( h\Delta t^3 + \Delta t^4)$. On the other hand, the point-wise value given by the reconstruction polynomial will be a linear combination of $Q_{l}$ and $hQ_{l,s}, s=x,y$. Thus, the fourth-order spatial accuracy can be maintained in the current reconstruction for the fourth-order compact GKS.
In order to improve the order of accuracy of the scheme which is related to the temporal discretisation, a third-order gas distribution function and higher-order MSMD methods can be applied \cite{ji1}.
In fact, the total order of accuracy of the scheme is limited by both the order of spatial and temporal discretisation.
Here the fourth-order compact GKS will be our main target, and the second-order gas evolution model and S2O4 method presented in section 3 are
sufficient for such a purpose.

\begin{figure}[!htb]
\centering
\includegraphics[width=0.4\textwidth]{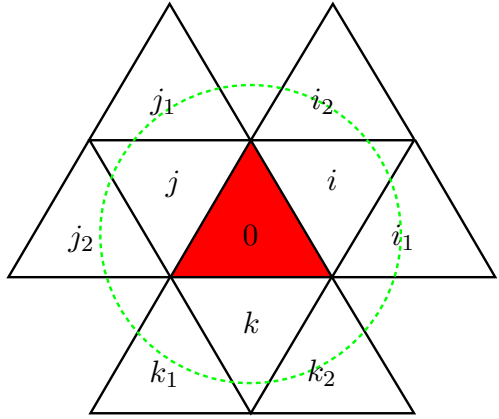}
\caption{\label{0-stencil} A typical stencil of compact high-order GKS from fourth order to sixth order on triangular meshes. Notice that some of the neighbor's neighbors ($\Omega_{i_1}$,$\Omega_{i_2}$,$\Omega_{j_1}$,...) may coincide for arbitrary triangular meshes. The green dotted circle represents the physical domain of dependence.}
\end{figure}

In order to determine the polynomial in Eq.(\ref{Recons-eqs}) for arbitrary geometrical triangular mesh, the number of equations $M$ in Eq.(\ref{Recons-eqs}) should be greater than $N$ to avoid an ill-conditioned system.
In the compact GKS, the reconstruction stencil in Fig. \ref{0-stencil} is used.
The high-order reconstruction will be obtained for the central cell $\Omega_0$.
Compared with the stencils of second-order scheme, the neighbors' neighbors ($\Omega_{i_1}$,$\Omega_{i_2}$,$\Omega_{j_1}$,...,$\Omega_{k_2}$) are needed for the compact fourth- to sixth-order GKS.
The dotted circle represents the physical domain of dependence, which indicates the physically compactness of the scheme.
Even for the regular isotropic mesh shown in Fig. \ref{0-stencil}, the physical domain of dependence of $\Omega_0$ contains neighboring cells of $\Omega_0$'s neighbors, i.e., $\Omega_{i_1},\Omega_{i_2},\cdots, \Omega_{k_2}$. Thus, the current compactness considers the compatibility from the physical domain of dependence in a numerical algorithm.

Since the cell averages and cell-averaged derivatives in each cell are provided, at most $30$ datum can be used in the  reconstruction.
For a series of compact schemes from fourth- up to sixth-order reconstruction, the following data sets are used to get the reconstructed polynomials,
\begin{align*}
S^{p^3}&=\{ Q_{l_1},Q_{l_2,x},Q_{l_2,y} \}; l_1=0,i,j,k,i_1,i_2,\cdots,k_1,k_2; l_2=0,i,j,k; \\
S^{p^4}&=\{ Q_{l_1},Q_{l_1,x},Q_{l_1,y} \}; l_1=0,i,j,k,i_1,i_2,\cdots,k_1,k_2;\\
S^{p^5}&=\{ Q_{l_1},Q_{l_1,x},Q_{l_1,y} \};
l_1=0,i,j,k,i_1,i_2,\cdots,k_1,k_2.
\end{align*}
Generally the least-square method can be adopted to solve the system (\ref{Recons-eqs}) for $a_k$. In the solution process,
the conservation condition is strictly satisfied, i.e., the first equation in Eq.(\ref{Recons-eqs}) for $l=0$, and the others hold in the least squares sense.
In this paper, a new constrained least square (CLS) method is adopted to solve the system (\ref{Recons-eqs}), where the cell averages of $\widetilde{\Omega}_0$, $\widetilde{\Omega}_i$, $\widetilde{\Omega}_j$ and $\widetilde{\Omega}_k$ are strictly satisfied by the polynomial $P^r(\xi,\eta)$. As a result, more accurate and smooth reconstructions are obtained at cell interfaces of $\Omega_0$.
The CLS solution for system (\ref{Recons-eqs}) is obtained by taking the extreme value of the function
\begin{align}\label{CLS-function}
\begin{split}
E= &( \widetilde{A}_{l_1(k),m}a_m-\widetilde{Q}_{l_1(k)})(\widetilde{A}_{l_1(k),n}a_n-\widetilde{Q}_{l_1(k)}) + \\
   &\lambda_{k} (A_{l_2(k),m}a_m-Q_{l_2(k)}),
\end{split}
\end{align}
where $\lambda_{k}$ is a Lagrangian factor.
 For the linear reconstruction on the large stencil $l_2=0,i,j,k$ and $l_1=i_1,i_2,\cdots,k_2$,
$A_{l_2,m}$ is the modified cell average of the basis $\varphi_m$ in $\widetilde{\Omega}_{l_2}$, i.e.,
\begin{align}\label{aver-base-1}
A_{l_2,m}=\frac{M_{l_2}}{\big|\widetilde{\Omega}_{l_2} \big|} \int_{ \widetilde{\Omega}_{l_2}} \varphi_m(\bm{\xi}) d\xi d\eta.
\end{align}
where the geometrically dependent coefficient $M_{l_2}=\big|\textbf{J}\big| \big|\widetilde{\Omega}_{l_2} \big|/\big|\Omega_{l_2} \big|$.
And $\widetilde{A}_{l_1(k),m}$ is the generalized modified cell average of the basis. $\widetilde{A}_{l_1(k),m}$ can be the modified cell average of $\varphi_m$ and the modified cell average of $\partial \varphi_m/ \partial s $ in $\widetilde{\Omega}_{l_1}$, such that the modified cell average of $\partial \varphi_m/ \partial s $ in $\widetilde{\Omega}_{l_1}$ is
\begin{align}\label{aver-base-2}
A^{'}_{l_1,m}=\frac{M_{l_1}}{\big|\widetilde{\Omega}_{l_1} \big|} \int_{ \widetilde{\Omega}_{l_1}}  \frac{\partial \varphi_m(\bm{\xi})}{\partial s} d\xi d\eta,
\end{align}
where $s=x,y$.
$\widetilde{Q}_{l_1}$ is the cell average or cell-averaged derivative in $\widetilde{\Omega}_{l_1}$.
Then the extreme value of $E$ is taken with the conditions
\begin{align*}
\frac{ \partial E }{\partial a_m}=0, ~~\frac{ \partial E }{\partial \lambda_{k}}=0,
\end{align*}
and the final linear system of $a_k$ and $\lambda$ becomes
\begin{equation} \label{CLS-system}
\left(
\begin{array}{cc}
2\widetilde{A}_{l_1(k),m} \widetilde{A}_{l_1(k),l} & \widetilde{A}_{l_2(k),m} \\
\widetilde{A}_{l_2(k),l} & 0 \\
\end{array}
\right)
\left(
\begin{array}{c}
a_l \\
\lambda_{k} \\
\end{array}
\right)
=
\left(
\begin{array}{c}
2\widetilde{A}_{l_1(k),m} \widetilde{Q}_{l_1(k)} \\
Q_{l_2(k)} \\
\end{array}
\right) .
\end{equation}

\begin{figure}[!htb]
	\centering
	\includegraphics[width=0.435\textwidth]{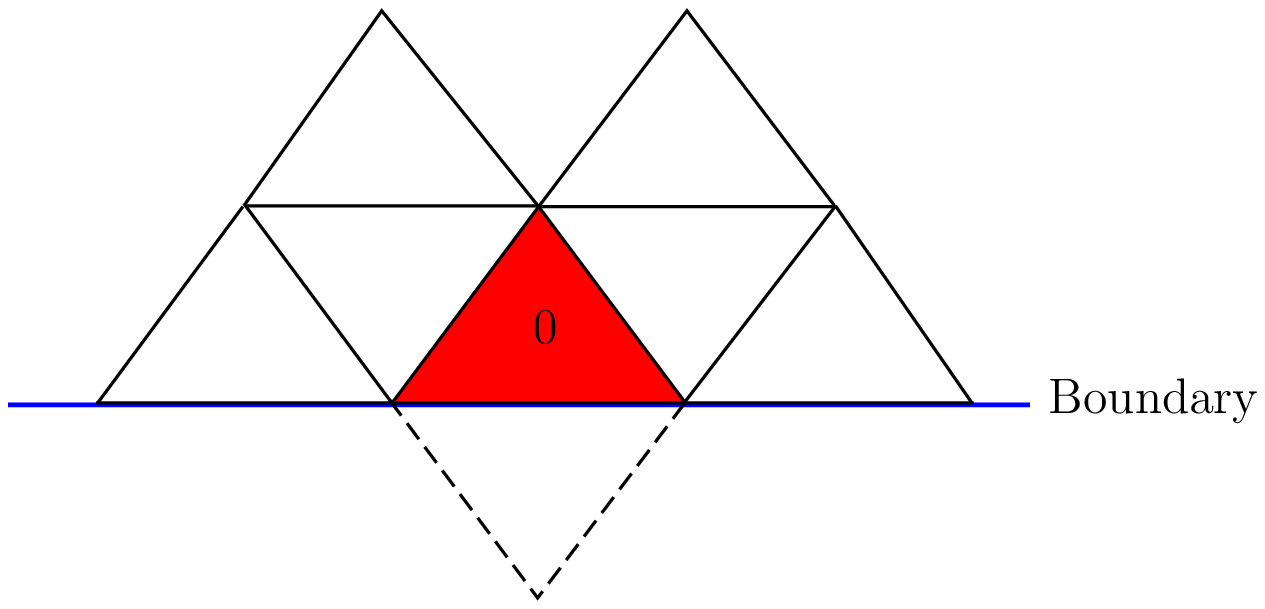}
	\includegraphics[width=0.375\textwidth]{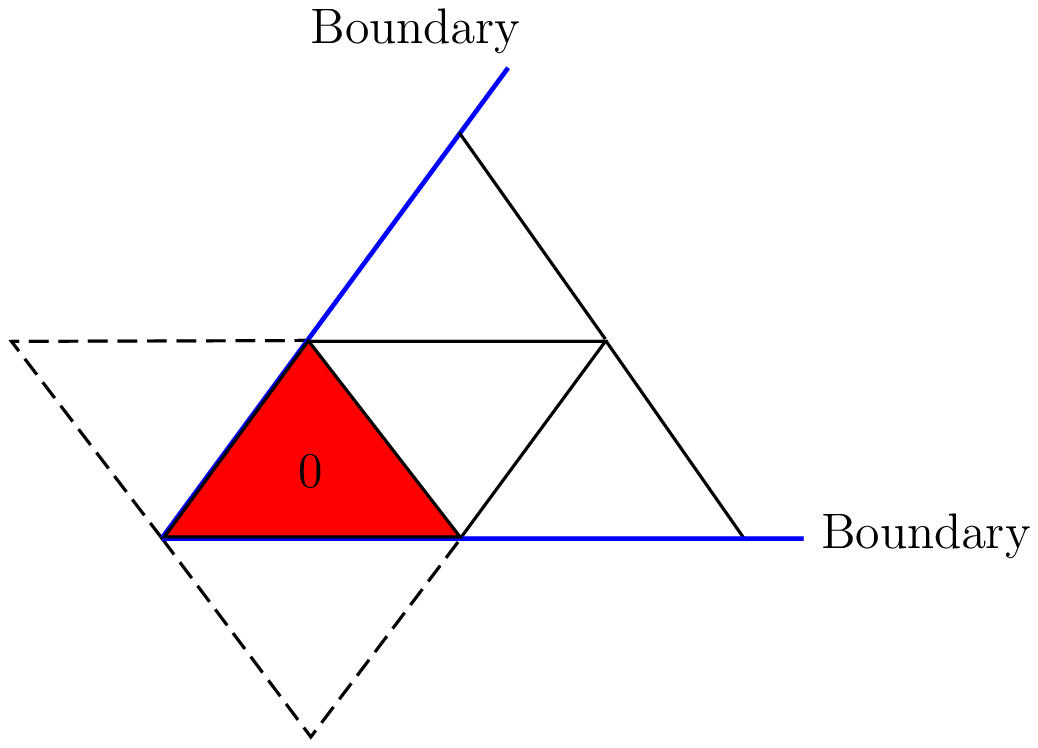}
	\caption{\label{0-stencil-bound} Simplified stencil of compact high-order GKS for a cell adjacent to boundaries of computational domain.}
\end{figure}

For the cell located at the boundary of the computational domain, the simplified stencil is used to get the high-order reconstruction. A graphical illustration is shown in Fig. \ref{0-stencil-bound}. One side or two sides of the cell $\Omega_0$ in Fig. \ref{0-stencil-bound} may coincide with the boundary of computational domain. Along the normal direction of each boundary, only one ghost cell of $\Omega_0$ is constructed, just as in the second-order scheme. A fourth order accuracy consistent with the inner mesh cells can be maintained for the reconstruction in $\Omega_0$. The only difference is that the orders of the candidate polynomials on some sub-stencils are reduced to one, while the orders of candidate polynomials don't affect the accuracy order of final reconstruction, which will be illustrated in the next section.
Numerical test cases demonstrate the same robustness by selecting the simplified stencils for cells adjacent to the boundary.

\section{Compact nonlinear reconstruction}
Nonlinear reconstruction is necessary for shock capturing scheme. The WENO idea is one of the best choices for the nonlinear reconstruction \cite{liu}. On structured mesh \cite{liu,jiang}, the reconstruction is based on the nonlinear convex combination of lower-order candidate polynomials. For the current nonlinear reconstruction on unstructured mesh, a simplified non-oscillatory reconstruction without optimal linear weights is obtained by the combination of the lower-order and high-order polynomials. This simplified reconstruction works well on structured mesh as well, which will not be presented here.

The candidate polynomials in nonlinear reconstruction are $q^r_k(\bm{\xi})$, where $r$ is the order of $q^r_k(\bm{\xi})$, and $k=1,2,\cdots,n$.
In the original WENO method, the value of a high-order polynomial $P^r(\bm{\xi})$ at a given point is a combination of several lower-order polynomials at that point,
\begin{equation}\label{convex-combi-1}
P(\bm{\xi})=\sum_{k=1}^{n}d_k q^r_k(\bm{\xi}),
\end{equation}
where $d_k$ are optimal linear weights and have positive values with the order $O(1)$.
Then, the nonlinear reconstruction is obtained by defining nonlinear weights \cite{jiang}.
In the cases with negative or  very large linear weights in Eq.(\ref{convex-combi-1}),  the scheme becomes unstable.
The modified WENO reconstruction on unstructured mesh is proposed in \cite{zhao2018WENO}.
In the current scheme, a simple non-oscillatory reconstruction is proposed.
The linear combination of lower-order polynomials in Eq.(\ref{convex-combi-1}) can be rewritten as
\begin{equation}\label{convex-combi-2}
P(\bm{\xi})=\sum_{k=1}^{n}\frac{C_k}{1+C} q_k(\bm{\xi}) +  \sum_{k=1}^{n}\frac{CC_k}{1+C}\big(\frac{1+C}{CC_k}d_k -\frac{1}{C}\big) q_k(\bm{\xi}),
\end{equation}
where $C$ is a free parameter, $\sum_{k=1}^n C_k=1$, $C>0$, and $0<C_k<1$.
Here $(1+C)/(CC_k)d_k -1/C=1$ when $C_k=d_k$.
Let's define normalized linear weights as
\begin{equation*}
\overline{d}_k=\frac{C_k}{1+C}, ~~ \overline{d}^{'}_k=\frac{CC_k}{1+C}, k=1,...,n.
\end{equation*}
Based on the linear convex combination of Eq.(\ref{convex-combi-2}), the nonlinear combination can be designed as
\begin{equation}\label{combination-3}
R(\bm{\xi})=\sum_{k=1}^{n}w_k q_k(\bm{\xi}) + \sum_{k=1}^{n} w^{'}_k \big[ \big(\frac{1+C}{C C_k}d_k-\frac{1}{C} \big)q_k(\bm{\xi}) \big],
\end{equation}
where $w_k$ and $w^{'}_k$ are the WENO weights. In this paper, the WENO-Z weight is used,
\begin{equation*}
\begin{split}
&w_{k}=\frac{\widetilde{w}_{k}}{ \sum_{j=1}^{n} \big(\widetilde{w}_{j}+\widetilde{w}^{'}_{j} \big) }, ~~~~~ w^{'}_{k}=\frac{\widetilde{w}^{'}_{k}}{ \sum_{j=1}^{n} \big(\widetilde{w}_{j}+\widetilde{w}^{'}_{j} \big) }, \\
&\widetilde{w}_{k}=\overline{d}_k \big(1+\big(\frac{\tau_Z}{IS_k +\epsilon}\big)^l\big), ~~ \widetilde{w}^{'}_{k}=\overline{d}_k^{'}\big(1+\big(\frac{\tau_Z}{IS_k +\epsilon}\big)^l\big), ~~k=0,...,n,
\end{split}
\end{equation*}
where $\tau_Z$ is a local high-order reference value, and $l$ is a parameter with a value $3$ in WNEO-Z weight.
The smooth indicator $IS_{k}$ is \cite{hu1999WENO}
\begin{equation*}
IS_k=\sum_{|\alpha|=1}^{r_k}|\Omega|^{|\alpha|-1}\iint_{\Omega}\big(D^{\alpha}q(\bm{\xi})\big)^2d\xi d\eta ,
\end{equation*}
where $r_k$ is the order of $q_k(\bm{\xi})$, and $\alpha$ is a multi-index, and $D$ is the derivative operator.
For example, when $\alpha=(1,2)$, then $|\alpha|=3$ and $D^{\alpha}q(\bm{\xi})=\partial^3 q(\bm{\xi}) /\partial \xi^1\partial \eta^2$.
For smooth solution, Eq.(\ref{combination-3}) can give the reconstruction with the same order of accuracy as linear reconstruction.
For discontinuous solution, the ``ENO" property will be  maintained.
The nonlinear weight $w_k$ of $q_k(\bm{\xi})$ when crossing discontinuity becomes $w_k \sim O(h^{2l})$. In smooth case the nonlinear weight $w_k$ for a smooth candidate polynomial $q_k(\bm{\xi})$ is $w_k \sim O(1)$.

A robust reconstruction can be obtained by taking $IS_0=max\{IS_k\}, k=1,2,...,n$, and $IS_0$ is used to determine all  $w^{'}_k$, $k=1,2,...,n$.
As a result, only the first part of the RHS of Eq.(\ref{combination-3}) remains in discontinuous case, while the second part vanishes with $O(h^{2l})$.
Then, Eq.(\ref{combination-3}) can be simplified as
\begin{equation}\label{combination-6}
\begin{split}
R(\bm{\xi})&=\sum_{k=1}^{n}w_k q_k(\bm{\xi}) + w_0 \sum_{k=1}^{n} \big[ \big(\frac{1+C}{C}d_k-\frac{C_k}{C} \big)q_k(\bm{\xi}) \big] \\
    &=\sum_{k=1}^{n}w_k q_k(\bm{\xi}) + w_0 \big( \frac{1+C}{C}P(\bm{\xi}) -\sum_{k=1}^{n} \frac{C_k}{C}q_k(\bm{\xi}) \big),
\end{split}
\end{equation}
where $P(\bm{\xi})$ is the high-order polynomial obtained by the linear reconstruction in Section 4. Thus, the final formula of the nonlinear reconstruction is a combination of the high-order and lower-order polynomials, i.e., the so-called a combination of ENO and WENO methodology.
The nonlinear weights $w_k$ can be simplified as well
\begin{equation}\label{nonlinear-wk}
\begin{split}
&w_{k}=\frac{ \widetilde{w}_{k} }{ \sum_{j=0}^{n} \widetilde{w}_{j}  }, \\
&\widetilde{w}_{k}=\overline{d}_k \big( 1+ \big( \frac{\tau_Z}{IS_k+\epsilon}  \big)^l \big),
\end{split}
\end{equation}
and $\overline{d}_{k}$ is
\begin{equation}\label{nonlinear-w0-1}
\overline{d}_{0}=\frac{C}{1+C}, ~~ \overline{d}_k=\frac{C_k}{1+C},~~ k=1,...,n.
\end{equation}
The compact scheme is insensitive to the values of $C$ and $C_k$. In this paper, we take $C=5$ and $C_k=1/n$.
Generally, $IS_0$ is bounded by the smooth indicator of the high-order polynomial from a large stencil, and $IS_0$ is taken as the smooth indicator of $P(\bm{\xi})$ in this paper.

\begin{figure}[!htb]
	\centering
	\includegraphics[width=0.35\textwidth]{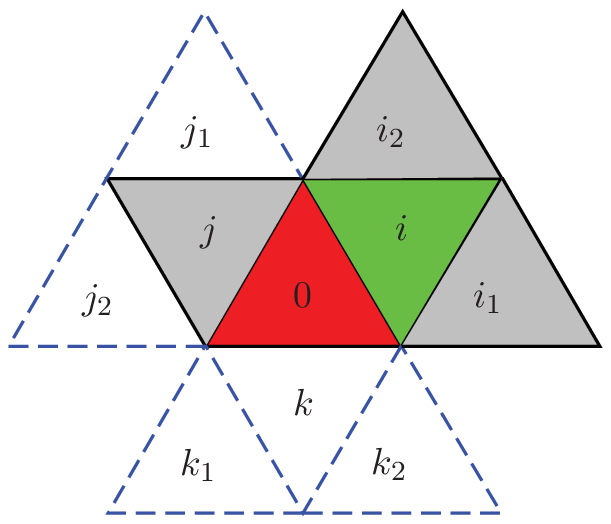}
	\includegraphics[width=0.35\textwidth]{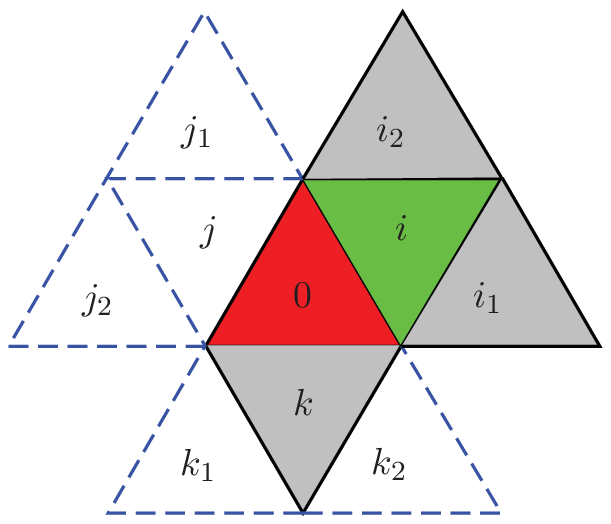}
	\caption{\label{0-stencil-sub} Two sub-stencils in the sub-stencils of compact GKS.}
\end{figure}

The sub-stencils for the nonlinear reconstruction are given by the seven choices,
\begin{align*}
S_1&=\{0,i,i_1,i_2,j\}, ~~~S_2=\{0,i,i_1,i_2,k\}, \\
S_3&=\{0,j,j_1,j_2,i\}, ~~~S_4=\{0,j,j_1,j_2,k\}, \\
S_5&=\{0,k,k_1,k_2,j\}, ~~S_6=\{0,k,k_1,k_2,i\}, \\
S_7&=\{0,i,j,k\}.\\
\end{align*}
The first two sub-stencils are presented in Fig. \ref{0-stencil-sub}. The last four sub-stencils are similar to the presented ones,
but in different directions. And the seventh stencil is a central sub-stencil.
For the sub-stencils $S_1,S_2, \cdots, S_6$, the cell averages of the cells in each sub-stencil and the cell-averaged derivatives of the second cell (highlighted in green in Fig. \ref{0-stencil-sub}) are used to obtain a quadratic polynomial. For the sub-stencil $S_7$, only the four cell averages are used to get a linear polynomial.

Based on the sub-stencils, the same CLS method given in Eq.(\ref{CLS-system}) is adopted to determine the candidate polynomials. In CLS,
in the determination of the candidate polynomials, only the cell average $Q_0$ is strictly satisfied, while the other cell averages are satisfied in the sense of least squares.
When the reconstruction cell is on the boundary of the computational domain as presented in Fig.\ref{0-stencil-bound},
at least two sub-stencils will contain three cells only. Since a quadratic polynomial cannot be fully determined from three cells,
a first-order polynomial is reconstructed. This treatment has no effect on the order of accuracy of the scheme in smooth region.

After the candidate polynomials are determined, the local higher-order reference value $\tau_Z$ in nonlinear weight can be constructed as
\begin{align*}
\tau_Z=\sum_{l=1}^3|2IS_0-IS_{2l-1}-IS_{2l}|.
\end{align*}
Such a definition of $\tau_Z$ can make the nonlinear weight approximate the linear one with $w_k=d_k+O(h^4), k=0,1,2,\cdots,7$. And the nonlinear reconstruction can achieve the formal fourth order up to sixth order of accuracy in the reconstruction.

In conclusion, the high-order reconstruction polynomials in each cell have been obtained based on the compact stencils. The CLS method is adopted to solve linear systems to give the high-order polynomial from a large stencil and lower-order candidate polynomials from the sub-stencils.
The simple WENO is proposed to give the nonlinear reconstruction based on the high-order polynomial and lower-order ones.
At the quadrature point on a cell interface, the conservative variables and their derivatives on both sides of the interface can be obtained by the nonlinear reconstruction. Then, the initial non-equilibrium state of the gas distribution function  and the equilibrium one in GKS are fully determined \cite{xu2}.

\section{Numerical examples }
In this section, we are going to test the compact high-order GKS on unstructured mesh.
The time step is determined by the CFL condition with $CFL=0.4$ in all test cases if not specified. For viscous tests case, the time step is also limited by the viscous term as $\Delta t=0.7h^2/(2\nu)$, where $h$ is the cell size and $\nu$ is the kinematic viscosity coefficient.
In all test cases, the same nonlinear reconstruction is used. There is no any additional ``trouble cell" detection or any additional limiter designed for specific test.
The collision time $\tau$ for inviscid flow at a cell interface is defined by
\begin{align*}
\tau=\varepsilon \Delta t + C_{num}\displaystyle|\frac{p_l-p_r}{p_l+p_r}|\Delta t,
\end{align*}
where $\varepsilon=0.05$, $C_{num}=5$, and $p_l$ and $p_r$ are the pressures at the left and right sides of a cell interface.
For the viscous flow, the collision time is related to the viscosity coefficient,
\begin{align*}
\tau=\frac{\mu}{p} + C_{num} \displaystyle|\frac{p_l-p_r}{p_l+p_r}|\Delta t,
\end{align*}
where  $\mu$ is the dynamic viscosity coefficient and $p$ is the pressure at the cell interface.
In smooth flow regions, it will reduce to $\tau=\mu/p$.
The reason for including pressure jump term in the particle collision time is to increase its value in the
discontinuous region in order to keep the physically consistent non-equilibrium mechanism in the construction of a  dissipative layer with mesh size
thickness.

\subsection{Accuracy test}
\begin{figure}[!htb]
\centering
\includegraphics[width=0.495\textwidth]{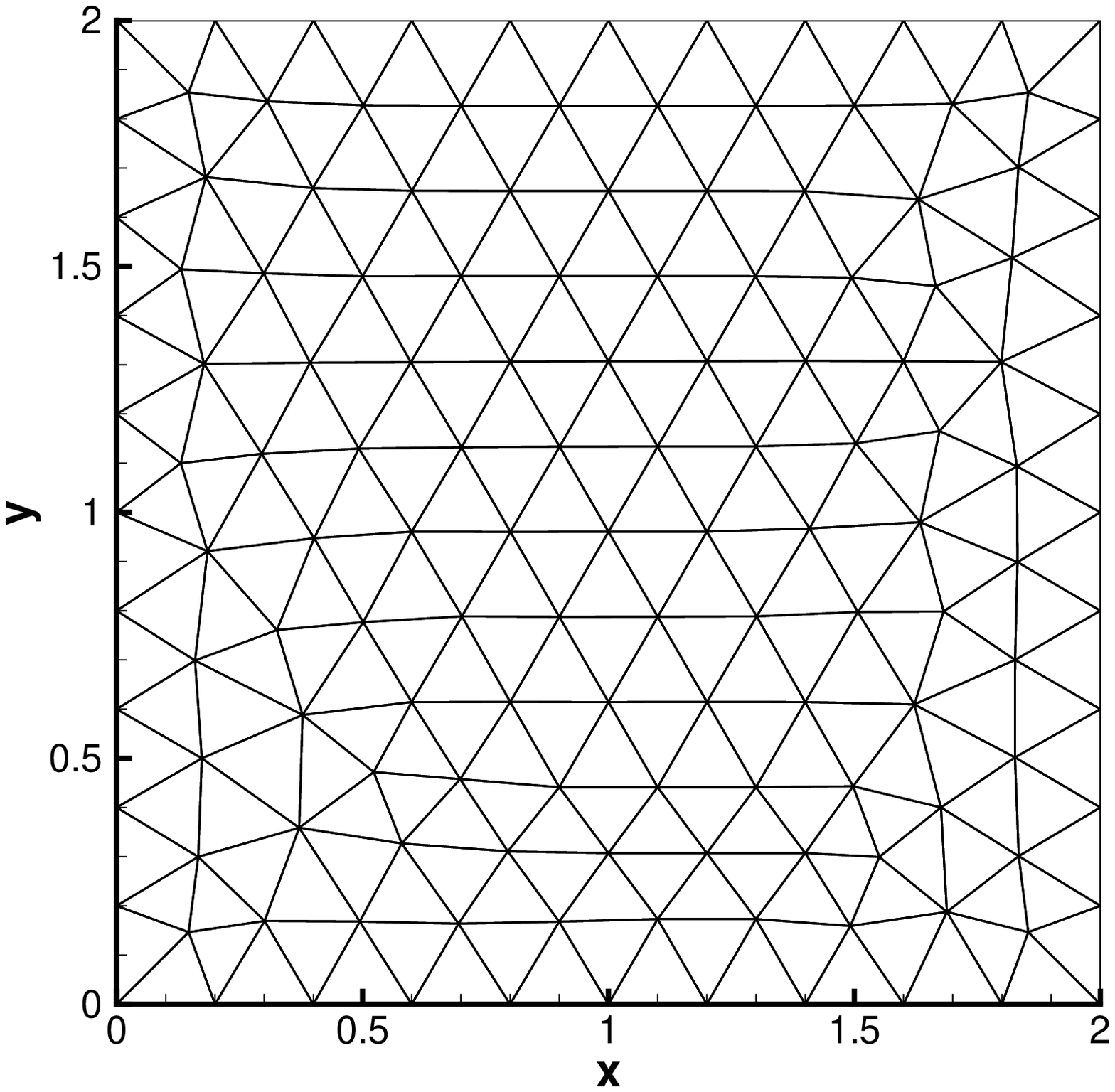}
\includegraphics[width=0.495\textwidth]{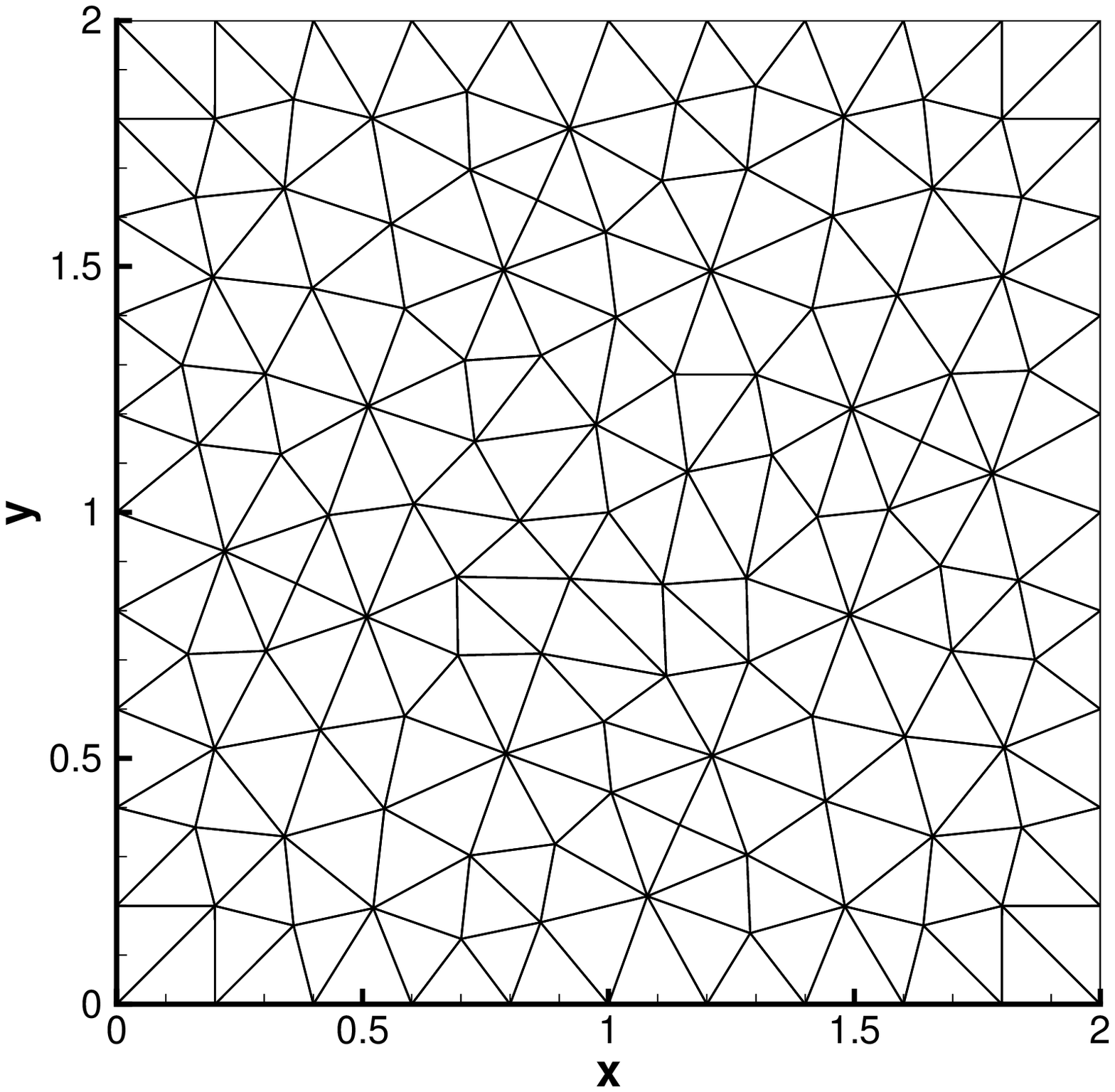}
\caption{\label{Mesh-Accuracy-1} The triangular meshes with cell size $h=1/5$: regular mesh (left) and irregular mesh (right).}
\end{figure}

The two-dimensional advection of density perturbation is used to verify the order of accuracy of compact high-order GKS. The initial condition is given as follows
\begin{align*}
\rho(x,y)=1+0.2\sin(\pi (x+y)),~U(x,y)=1,~V(x,y)=1,~p(x,y)=1.
\end{align*}
The computational domain is $[0,2]\times[0,2]$, and the periodic boundary conditions are applied in both directions.
Two types of triangular meshes are used, i.e., the regular mesh and irregular mesh shown in  Fig.\ref{Mesh-Accuracy-1}.
Except for some small regions, such as those near the boundary of the computational domain, the cells of the regular mesh are very close to equilateral triangles. And the size and shape of triangular cells in irregular mesh are arbitrary.

With the $r$th-order spatial reconstruction and S2O4 temporal discretisation, the leading term of the truncation error is
$O(h^r,\Delta t^4,h\Delta t^3)$, where the third term is from the evolution of conservative variables at cell interface. To verify the accuracy order of reconstruction, $\Delta t \sim \min \big\{h^{r/4},h^{(r-1)/3}\big\}$ is used for the $r$th-order scheme.
The $L^1$ and $L^\infty$ errors and convergence orders at $t=2$ are presented in Table \ref{Accu-sin-p3-linear} to Table \ref{Accu-sin-p5-linear}, respectively.
Due to the non-uniform meshes, the accuracy order of $L^\infty$ cannot reflect the true convergence order of the scheme. From the numerical results listed in Table \ref{Accu-sin-p3-linear} to Table \ref{Accu-sin-p5-linear}, it can be seen that whether for regular or irregular meshes, the accuracy orders of $L^1$ are almost the same as the theoretical ones.
If the time step takes $dt=h^{r/4}$, as that used the Runge-Kutta-based schemes, the same accuracy orders can be achieved, which will not be presented here.

\begin{table}[!h]
	\small
	\begin{center}
		\def\temptablewidth{1.0\textwidth}
		{\rule{\temptablewidth}{0.50pt}}
        \footnotesize
		\begin{tabular*}{\temptablewidth}{@{\extracolsep{\fill}}c|cc|cc|c|cc|cc}
$h_{re}$ &$Error_{L^1}$ &$\mathcal{O}_{L^1}$ &$Error_{L^{\infty}}$ &$\mathcal{O}_{L^{\infty}}$ &$h_{irre}$ &$Error_{L^1}$ &$\mathcal{O}_{L^1}$ &$Error_{L^{\infty}}$ &$\mathcal{O}_{L^{\infty}}$   \\
			\hline
            1/5   & 4.0326e-04    & ~     & 1.0977e-03  & ~      & 1/5     & 7.4212e-04 & ~     & 2.2584e-03   & ~      \\
            1/10  & 2.9859e-05    & 3.76  & 7.7770e-05  & 3.82   & 1/10    & 4.7778e-05 & 3.96  & 1.4986e-04   & 3.91   \\
            1/20  & 1.5647e-06    & 4.25  & 5.7130e-06  & 3.77   & 1/20    & 3.2422e-06 & 3.88  & 1.0439e-05   & 3.84   \\
            1/40  & 9.5781e-08    & 4.03  & 3.1262e-07  & 4.19   & 1/40    & 2.2276e-07 & 3.86  & 1.0642e-06   & 3.29   \\
		\end{tabular*}
		{\rule{\temptablewidth}{0.50pt}}
	\end{center}
	\vspace{-6mm} \caption{\label{Accu-sin-p3-linear} Accuracy test for compact fourth-order GKS: errors and convergence orders of linear compact fourth-order scheme with $P^3(x,y)$ reconstruction on regular and irregular meshes. The time step $\Delta t$ is given by $CFL=1.0$. }
\end{table}

\begin{table}[!h]
	\small
	\begin{center}
		\def\temptablewidth{1.0\textwidth}
		{\rule{\temptablewidth}{0.50pt}}
        \footnotesize
		\begin{tabular*}{\temptablewidth}{@{\extracolsep{\fill}}c|cc|cc|c|cc|cc}
$h_{re}$ &$Error_{L^1}$ &$\mathcal{O}_{L^1}$ &$Error_{L^{\infty}}$ &$\mathcal{O}_{L^{\infty}}$ &$h_{irre}$ &$Error_{L^1}$ &$\mathcal{O}_{L^1}$ &$Error_{L^{\infty}}$ &$\mathcal{O}_{L^{\infty}}$   \\
			\hline
            1/5   & 3.2702e-04    & ~     & 6.1466e-04  & ~      & 1/5     & 3.3221e-04 & ~     & 6.8051e-04        & ~      \\
            1/10  & 1.1291e-05    & 4.86  & 2.1244e-05  & 4.85   & 1/10    & 9.2734e-06 & 5.16  & 2.7177e-05        & 4.65   \\
            1/20  & 3.6287e-07    & 4.96  & 6.9429e-07  & 4.94   & 1/20    & 2.7593e-07 & 5.07  & 7.2777e-07        & 5.22   \\
            1/40  & 8.1121e-09    & 5.48  & 1.8885e-08  & 5.20   & 1/40    & 6.5075e-09 & 5.41  & 3.1871e-08        & 4.51   \\
		\end{tabular*}
		{\rule{\temptablewidth}{0.50pt}}
	\end{center}
	\vspace{-6mm} \caption{\label{Accu-sin-p4-linear} Accuracy test for compact fifth-order GKS: errors and convergence orders of linear compact fifth-order scheme with $P^4(x,y)$ reconstruction on regular and irregular meshes. The time step is determined by $\Delta t=h^{4/3}$. }
\end{table}

\begin{table}[!h]
	\small
	\begin{center}
		\def\temptablewidth{1.0\textwidth}
		{\rule{\temptablewidth}{0.50pt}}
        \footnotesize
		\begin{tabular*}{\temptablewidth}{@{\extracolsep{\fill}}c|cc|cc|c|cc|cc}
			
$h_{re}$ &$Error_{L^1}$ &$\mathcal{O}_{L^1}$ &$Error_{L^{\infty}}$ &$\mathcal{O}_{L^{\infty}}$ &$h_{irre}$ &$Error_{L^1}$ &$\mathcal{O}_{L^1}$ &$Error_{L^{\infty}}$ &$\mathcal{O}_{L^{\infty}}$   \\
			\hline
            2/5  & 1.5914e-03 & ~    & 4.6599e-03 & ~    & 2/5  & 1.1590e-03 & ~    & 3.4237e-03 & ~    \\
            1/5  & 2.5709e-05 & 5.95 & 5.7810e-05 & 6.33 & 1/5  & 2.9413e-05 & 5.30 & 8.2085e-05 & 5.38 \\
            1/10 & 4.5954e-07 & 5.81 & 1.1861e-06 & 5.61 & 1/10 & 4.8447e-07 & 5.92 & 1.6950e-06 & 5.60 \\
            1/20 & 5.2435e-09 & 6.45 & 2.0941e-08 & 5.82 & 1/20 & 6.9893e-09 & 6.12 & 3.5966e-08 & 5.56 \\
		\end{tabular*}
		{\rule{\temptablewidth}{0.50pt}}
	\end{center}
	\vspace{-6mm} \caption{\label{Accu-sin-p5-linear} Accuracy test for compact sixth-order GKS: errors and convergence orders of linear compact sixth-order scheme with $P^5(x,y)$ reconstruction on regular and irregular meshes. The time step is determined by $\Delta t=h^{5/3}$. }
\end{table}

\subsection{One dimensional Riemann problem}
\begin{figure}[!htb]
\centering
\includegraphics[width=0.485\textwidth]{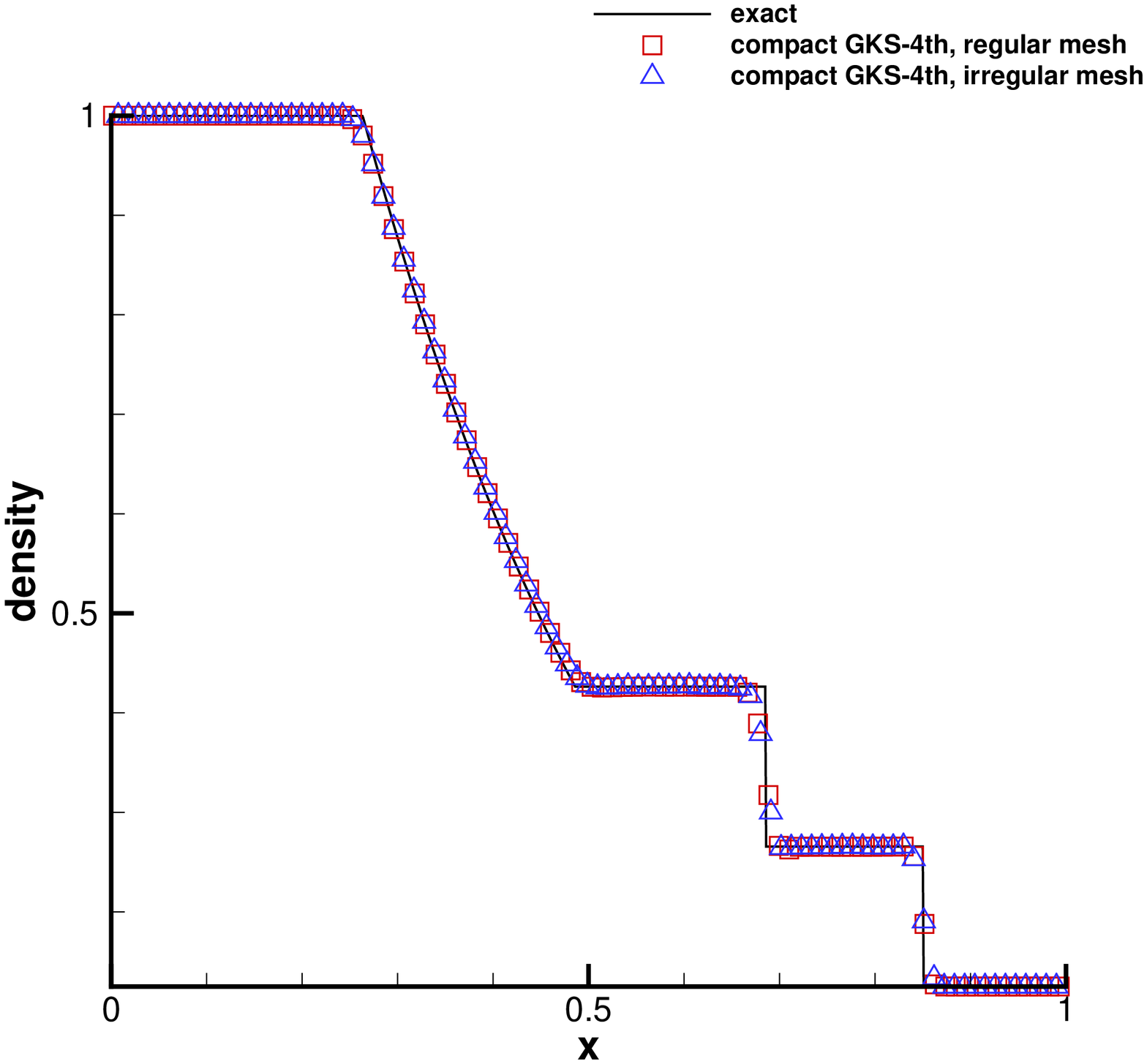}
\includegraphics[width=0.485\textwidth]{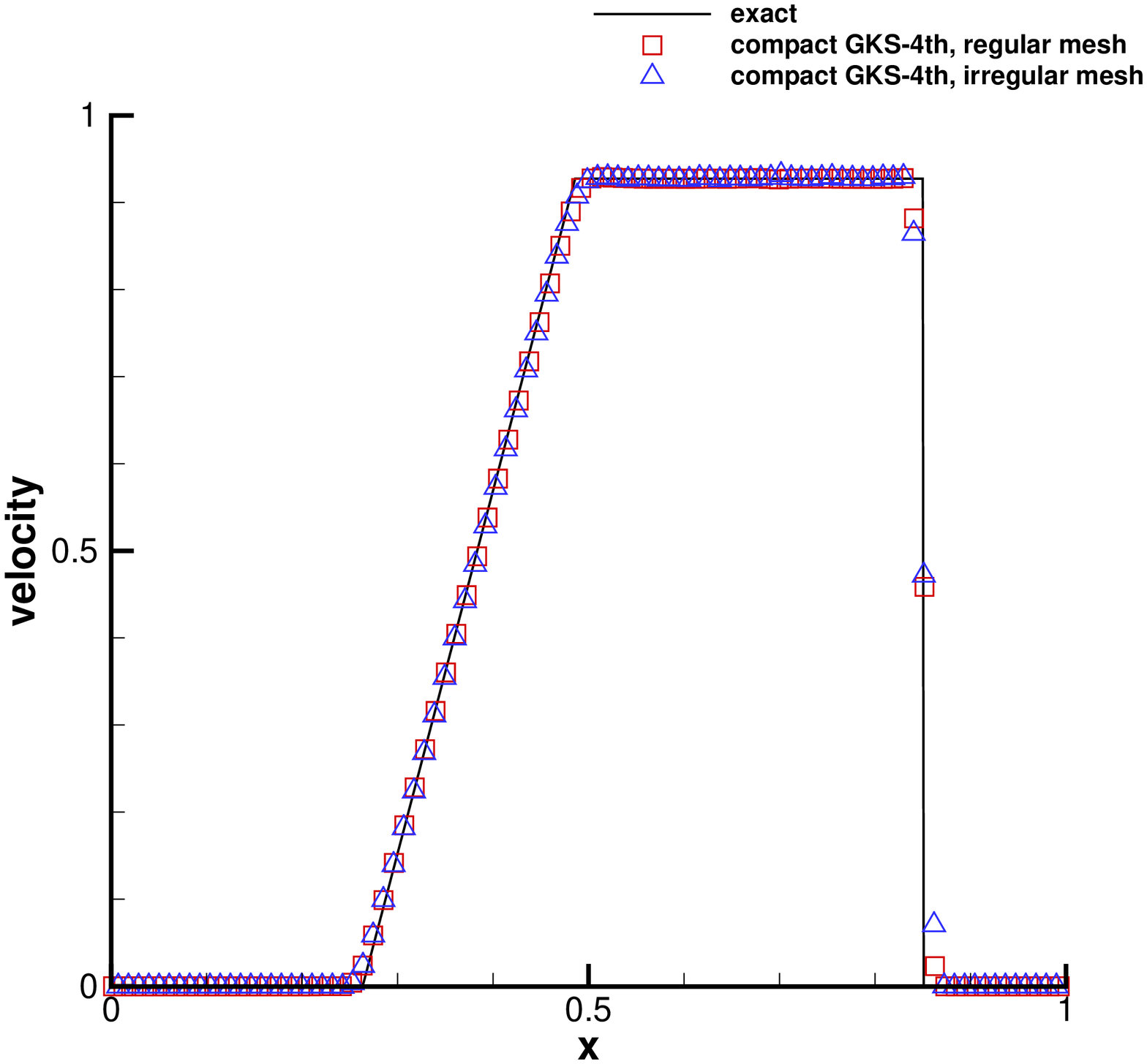}
\caption{\label{1d-shocktube-1} One dimensional Riemann problem: density and velocity distributions $t=0.2$ along the horizontal center line of Sod problem from compact fourth-order GKS. The mesh cell size is $h=1/100$. }
\end{figure}

\begin{figure}[!htb]
\centering
\includegraphics[width=0.485\textwidth]{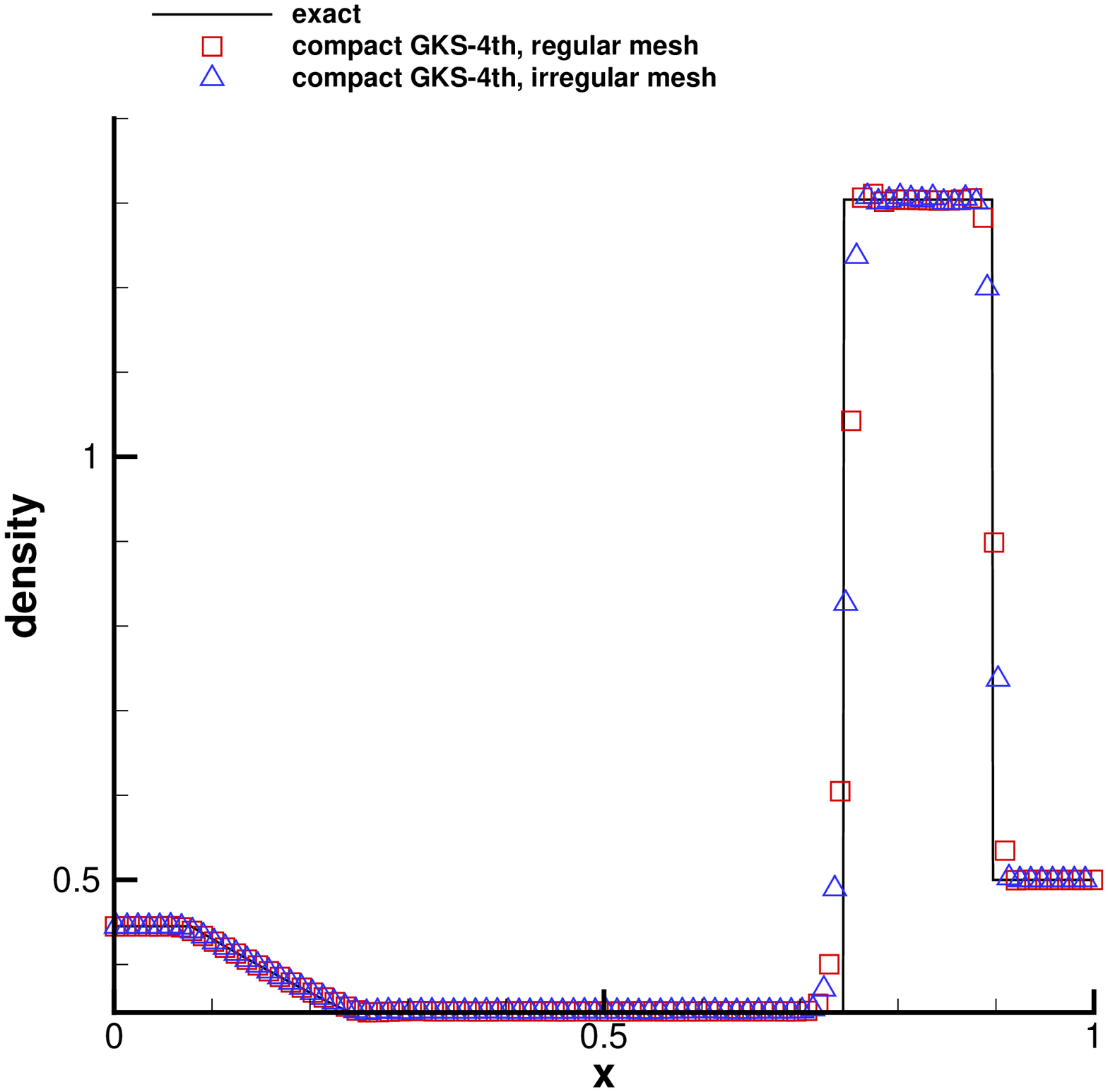}
\includegraphics[width=0.485\textwidth]{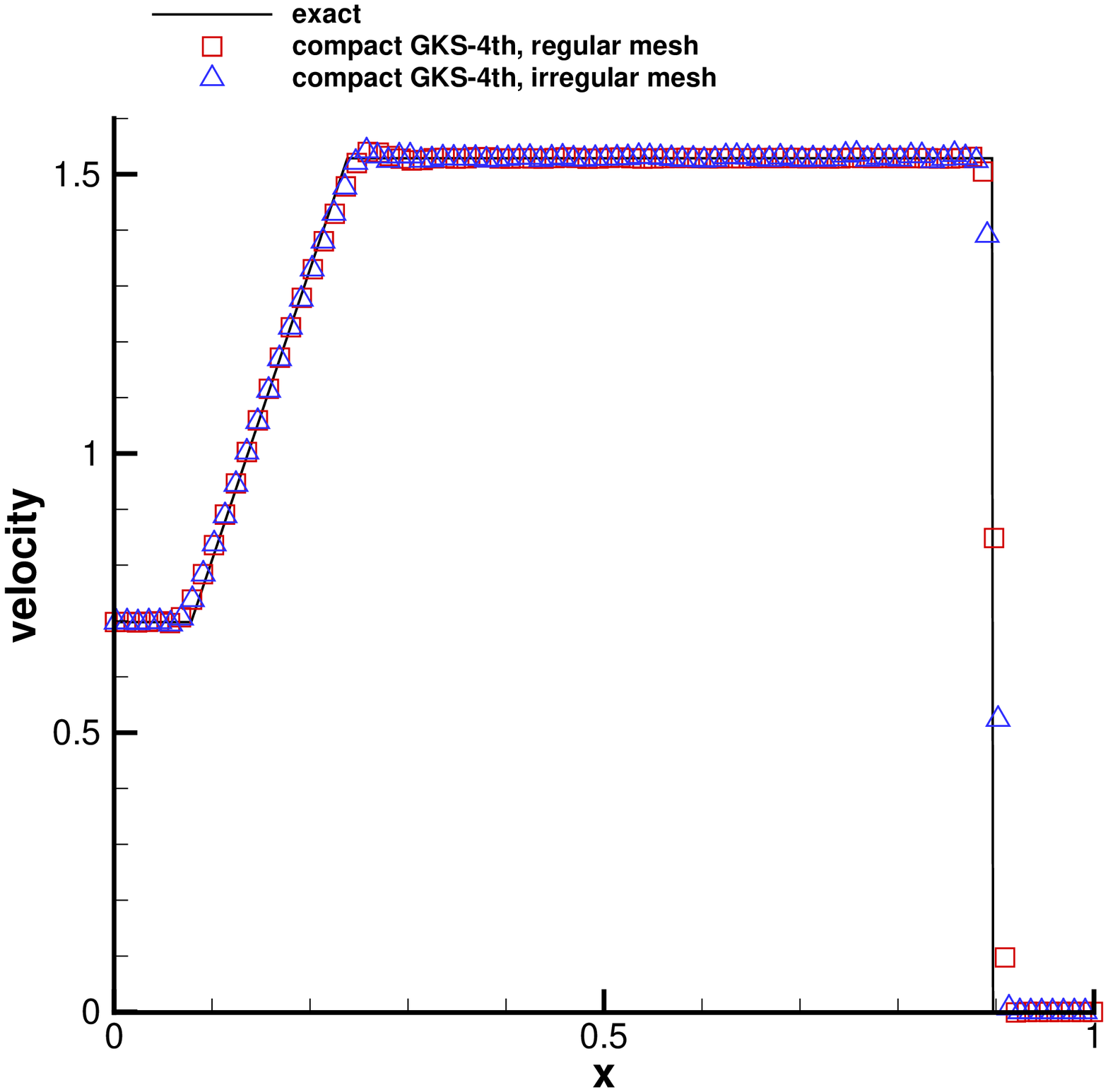}
\caption{\label{1d-shocktube-2} One dimensional Riemann problem: density and velocity distributions $t=0.16$ along the horizontal center line of Lax problem from compact fourth-order GKS. The mesh cell size is $h=1/100$. }
\end{figure}

Three one-dimensional Riemann problems in two-dimensional unstructured mesh are tested here.
The shock tube problem is tested firstly to show the essentially non-oscillatory property for shock waves. The initial condition of Sod problem is
\begin{equation*}
(\rho,U,p) = \begin{cases}
(1,0,1),  0\leq x<0.5,\\
(0.125, 0, 0.1), 0.5\leq x\leq1.
\end{cases}
\end{equation*}
And the initial condition of Lax problem is
\begin{equation*}
(\rho,U,p) = \begin{cases}
(0.445, 0.698, 3.528),  0\leq x<0.5,\\
(0.5, 0, 0.571), 0.5\leq x\leq1.
\end{cases}
\end{equation*}
The computational domain is $[0,1]\times[0,0.5]$. The reflecting boundary condition is applied in the $y$ direction. The regular and irregular triangular meshes are used with cell size $h=1/100$. The density and velocity distributions at $t=0.2$ and $t=0.16$ along the horizontal center line are shown in Fig.\ref{1d-shocktube-1} and Fig.\ref{1d-shocktube-2}, respectively. The compact fourth-order GKS gives the essentially non-oscillatory solutions with high resolution.

\begin{figure}[!htb]
\centering
\includegraphics[width=0.485\textwidth]{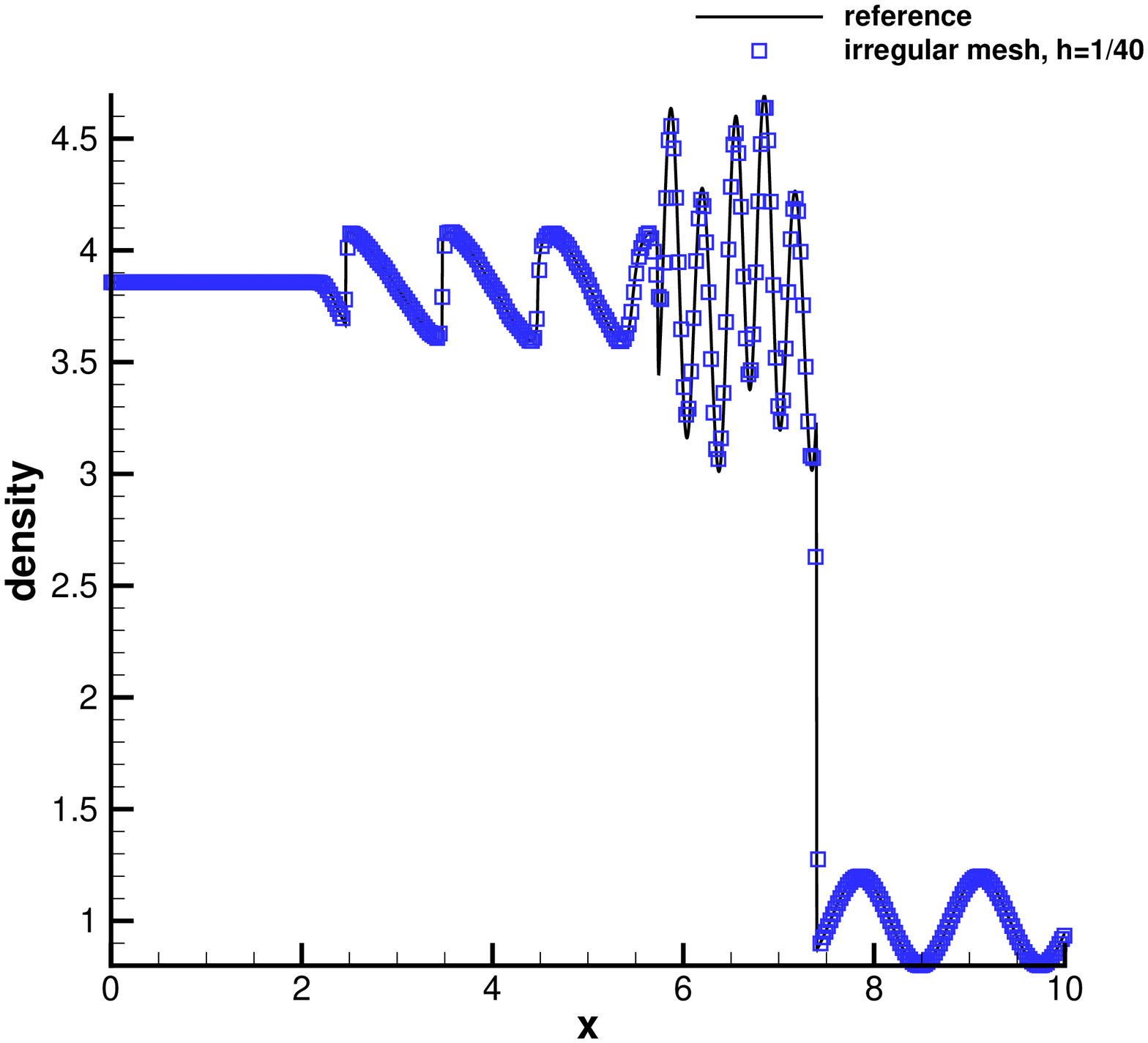}
\includegraphics[width=0.485\textwidth]{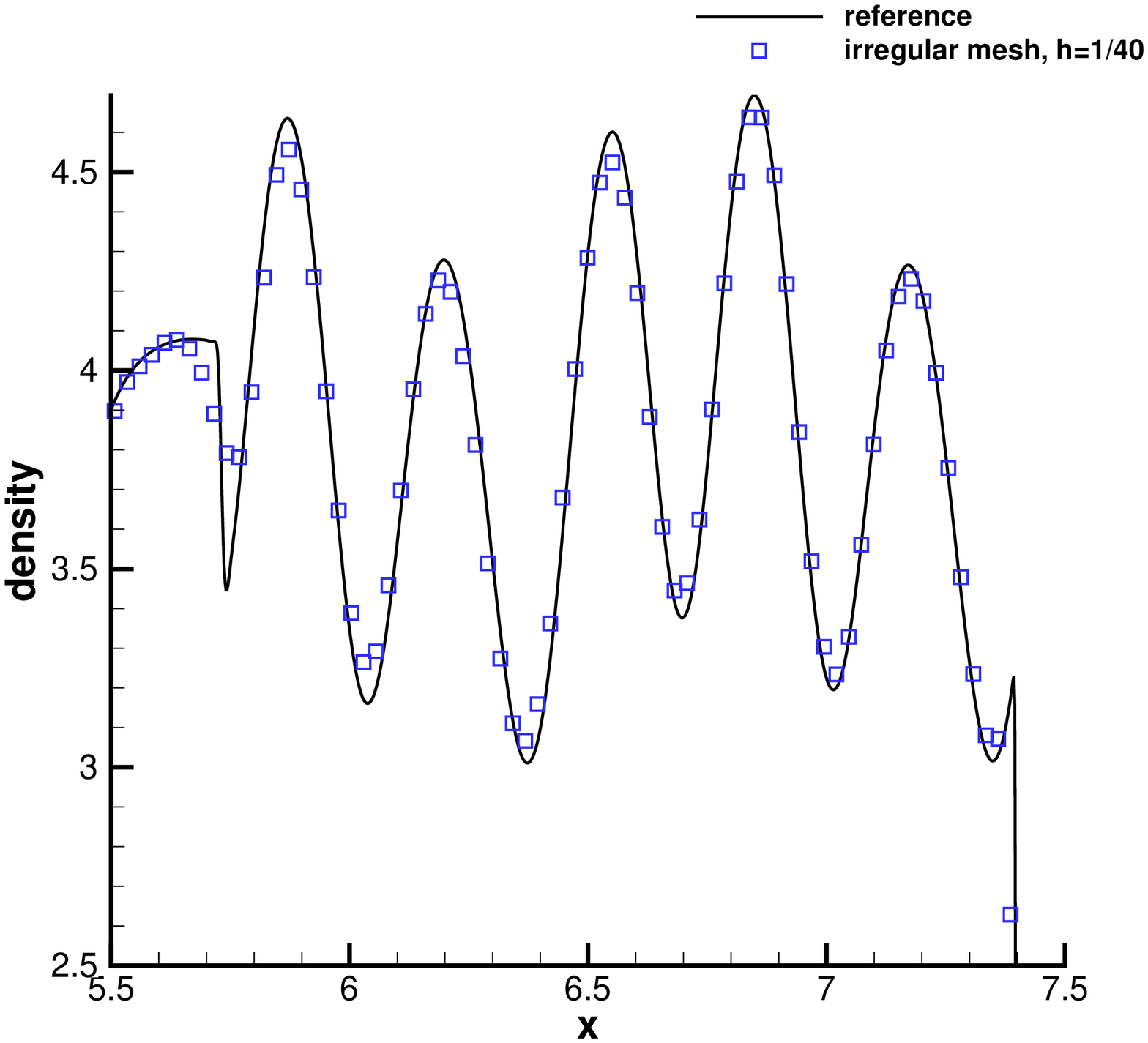}
\caption{\label{1d-shocktube-shu} One dimensional Riemann problem: the density distribution, local enlargement for Shu-Osher problem at $t=1.8$ on irregular triangular mesh with $h=1/40$.}
\end{figure}

To test the performance of capturing high frequency waves, the Shu-Osher problem \cite{Case-Shu-Osher} is tested, which is a case of a shock wave interacting with the density wave. The initial condition is given as follows
\begin{equation*}
(\rho,U,p)=\left\{\begin{array}{ll}
(3.857134, 2.629369, 10.33333),  \ \ \ \ &  x \leq 1,\\
(1 + 0.2\sin (5x), 0, 1),  &  1 <x.
\end{array} \right.
\end{equation*}
The computational domain is $[0, 10]\times[0, 0.25]$ and $h=1/40$ irregular mesh is used. The reflected boundary condition is applied in the $y$ direction. The density distributions and local enlargement at the center line at $t=1.8$ are presented in Fig.\ref{1d-shocktube-shu}. The current compact GKS demonstrates a high resolution for the high frequency wave.

\begin{figure}[!htb]
\centering
\includegraphics[width=0.485\textwidth]{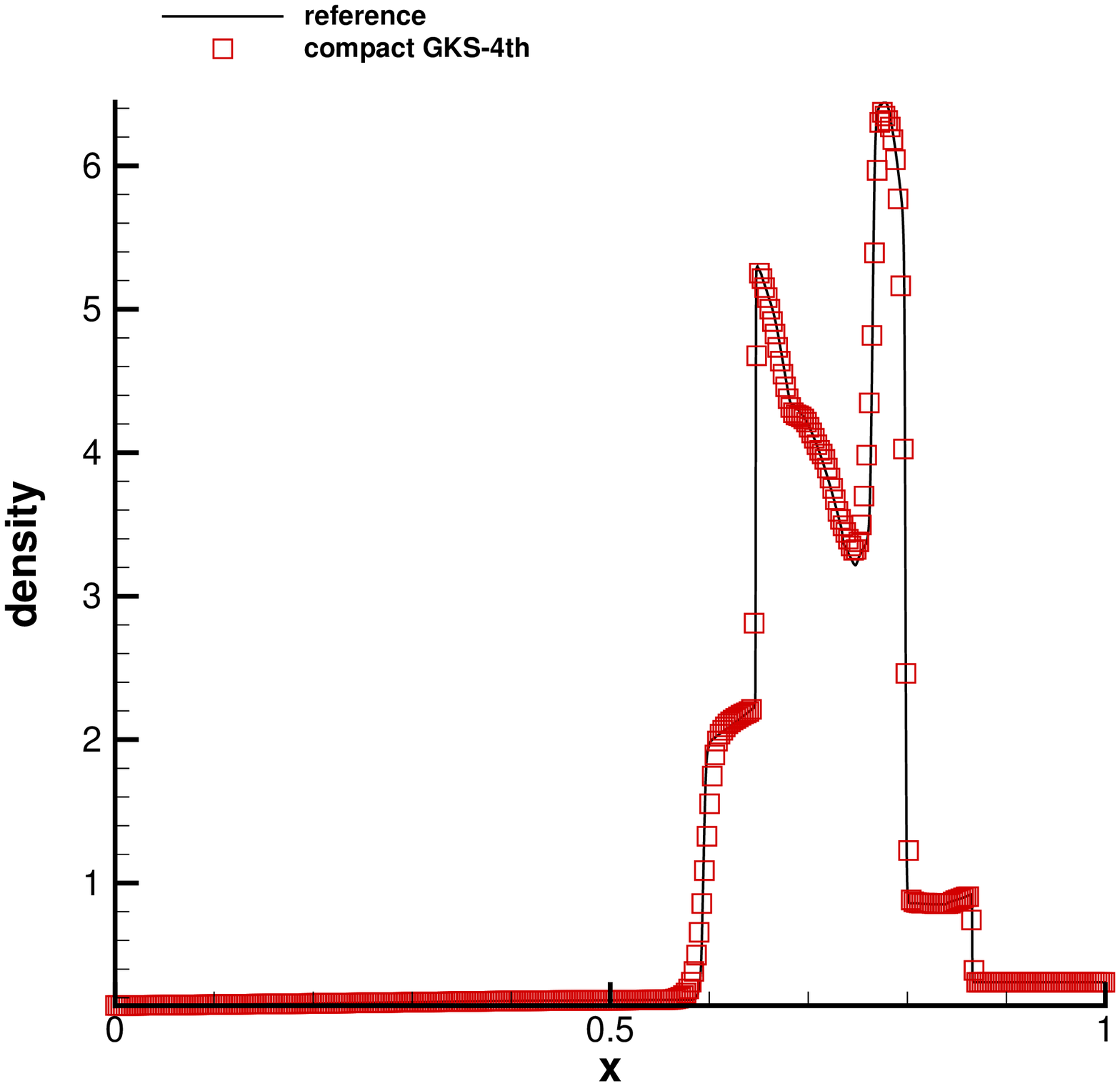}
\includegraphics[width=0.485\textwidth]{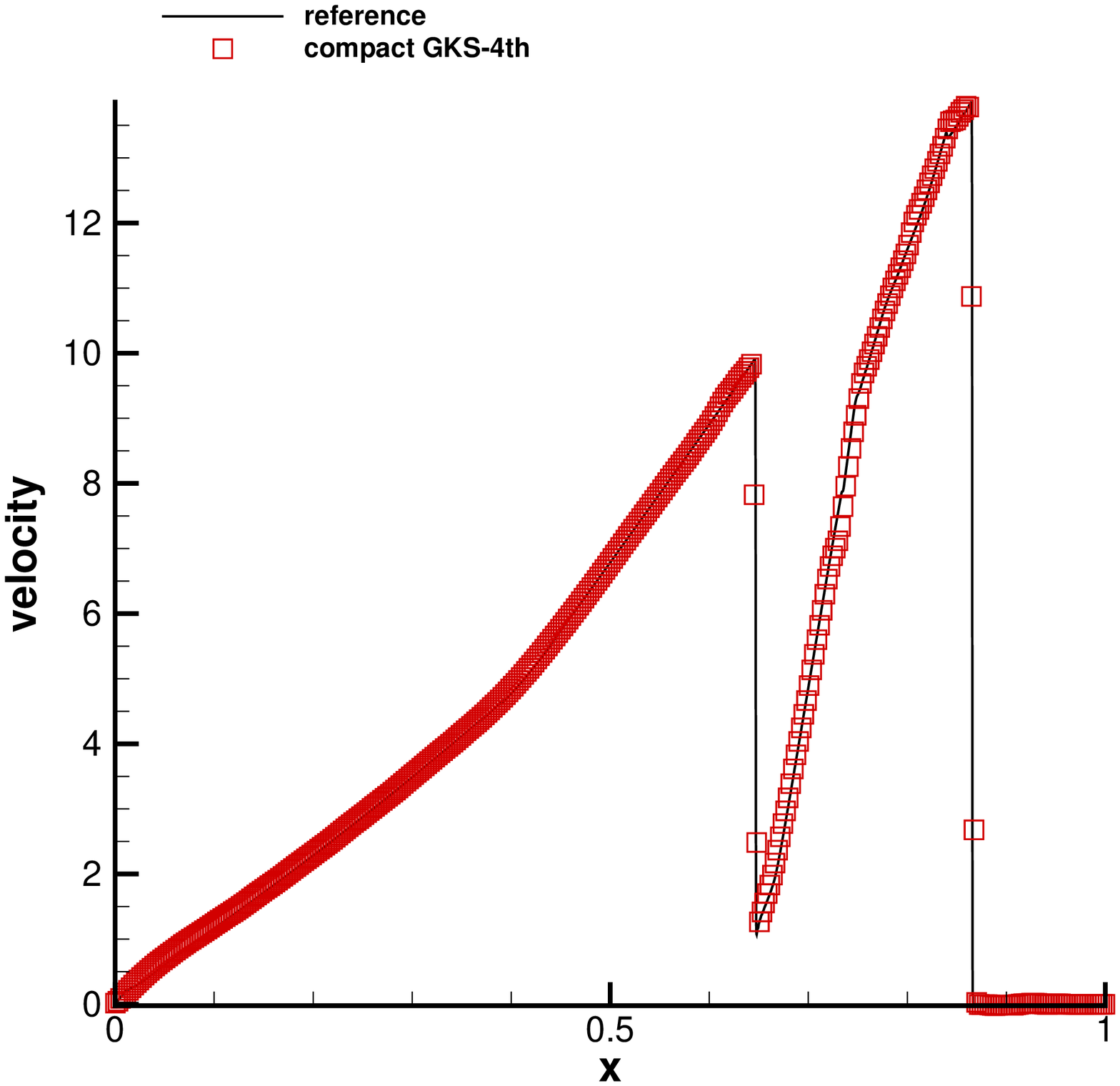}
\caption{\label{1d-blastwave} One dimensional Riemann problem: the density distribution, local enlargement for blast wave problem at $t=0.038$ on regular triangular mesh with $1/400$.}
\end{figure}

The third case is the Woodward-Colella blast wave problem \cite{Case-Woodward}. The test case verifies the robustness of compact GKS for capturing strong shock wave on unstructured mesh. The initial condition is given as follows
\begin{align*}
(\rho,U,p) =\begin{cases}
(1, 0, 1000), & 0\leq x<0.1,\\
(1, 0, 0.01), & 0.1\leq x<0.9,\\
(1, 0, 100),  & 0.9\leq x\leq 1.
\end{cases}
\end{align*}
The computational domain is $[0,1]\times[0,0.25]$, and the reflecting boundary conditions are imposed on both directions. The regular mesh with $h=1/400$ is used. In this test case, the $CFL$ number takes $0.2$. The computed density and velocity profiles at $t=0.038$ are shown in Fig. \ref{1d-blastwave}. The numerical results agree well with the reference solutions. The scheme can resolve the wave profiles well, particularly for the local extreme values. And there are not numerical oscillations near the shock waves.

\subsection{Double Mach reflection}
\begin{figure}[!htb]
\centering
\includegraphics[width=0.80\textwidth]{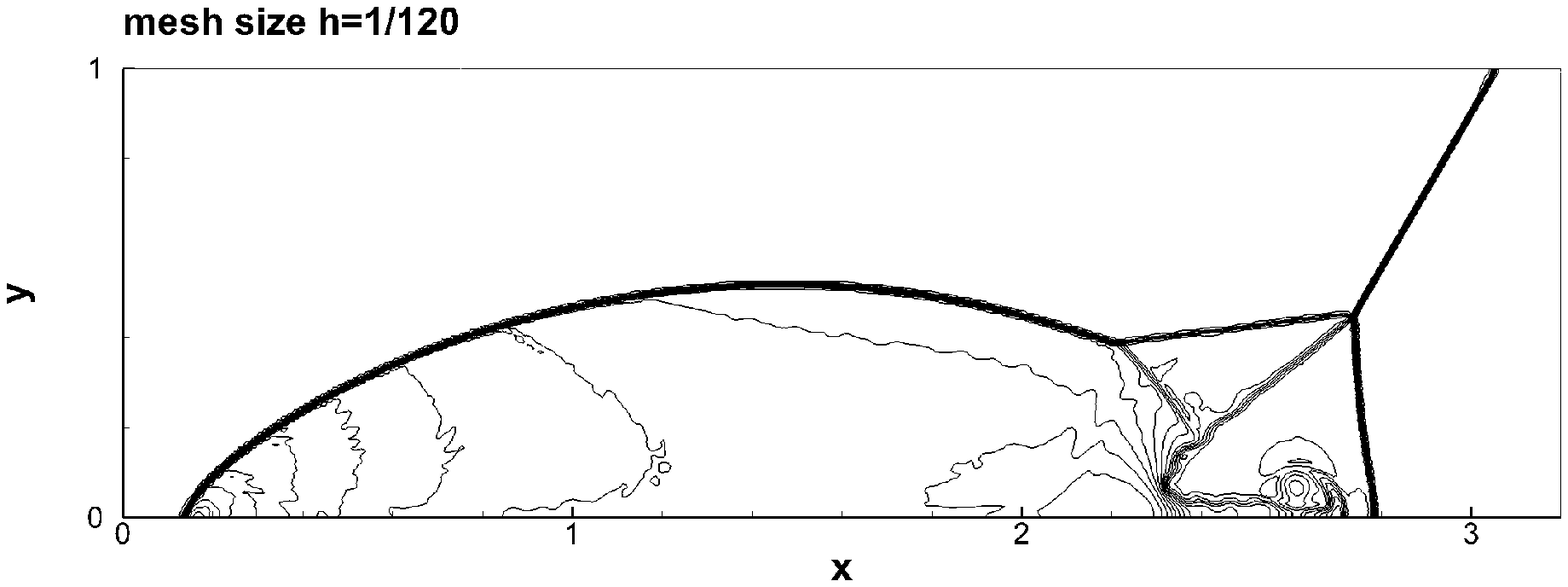}
\includegraphics[width=0.80\textwidth]{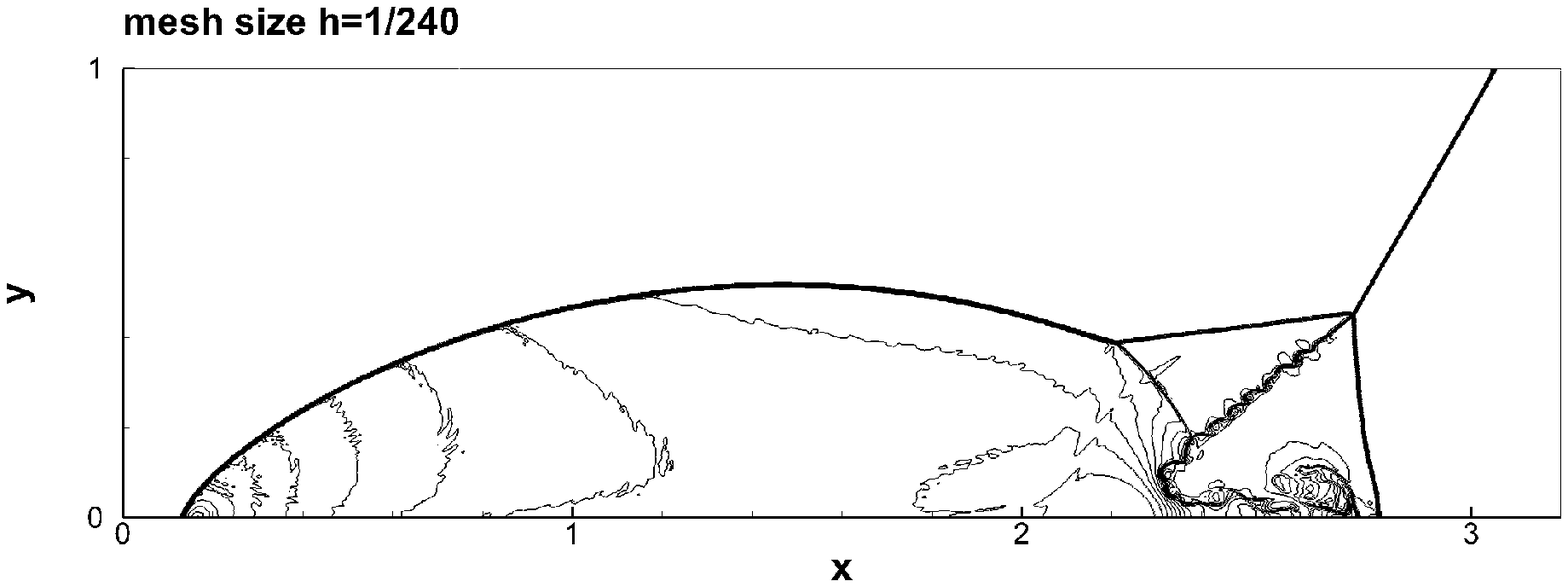}
\caption{\label{2d-doublemach-1} Density distribution on the triangular mesh with $h=1/120$ and $1/240$ at $t=0.2$.}
\end{figure}

\begin{figure}[!htb]
\centering
\includegraphics[width=0.45\textwidth]{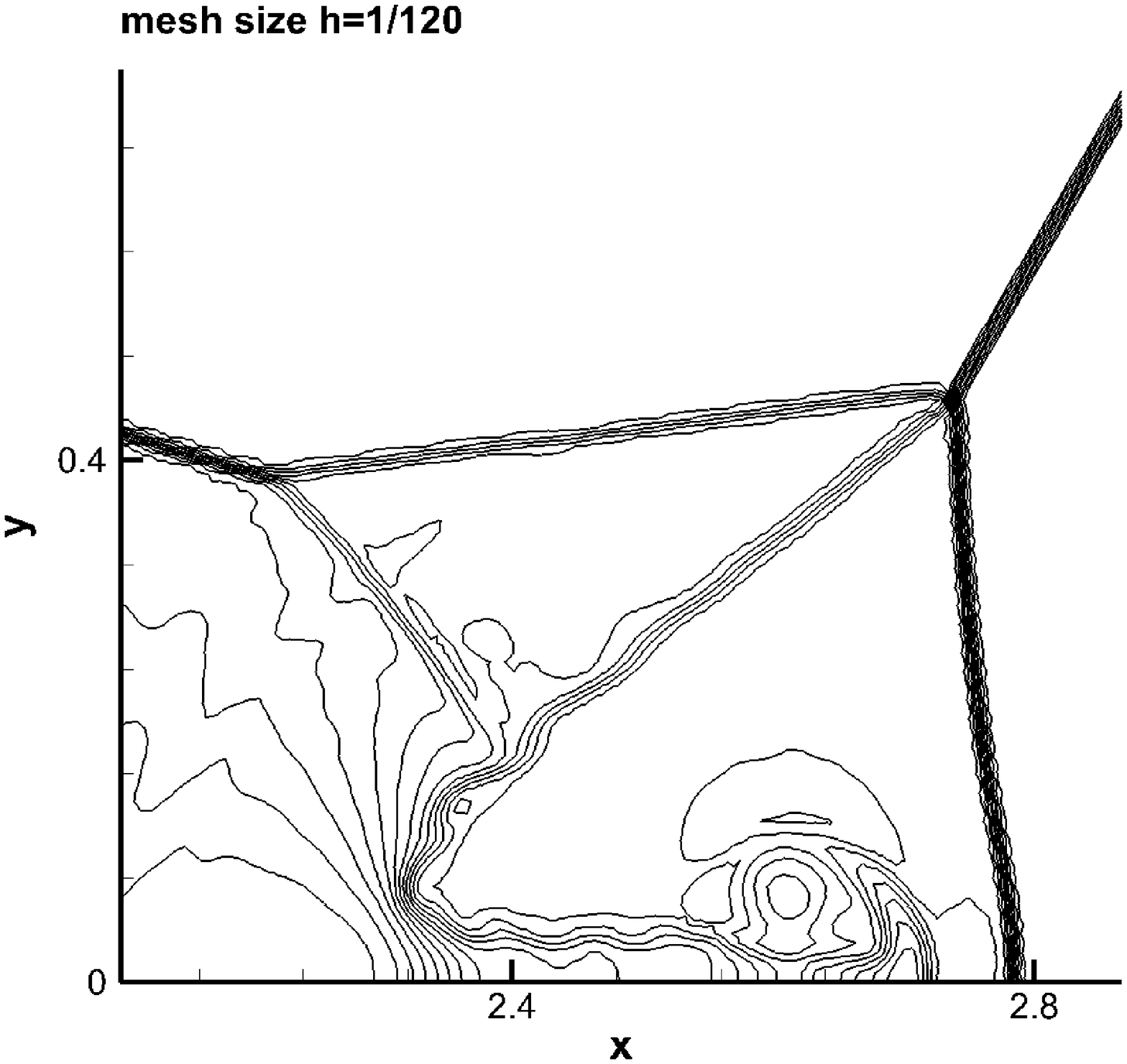}
\includegraphics[width=0.45\textwidth]{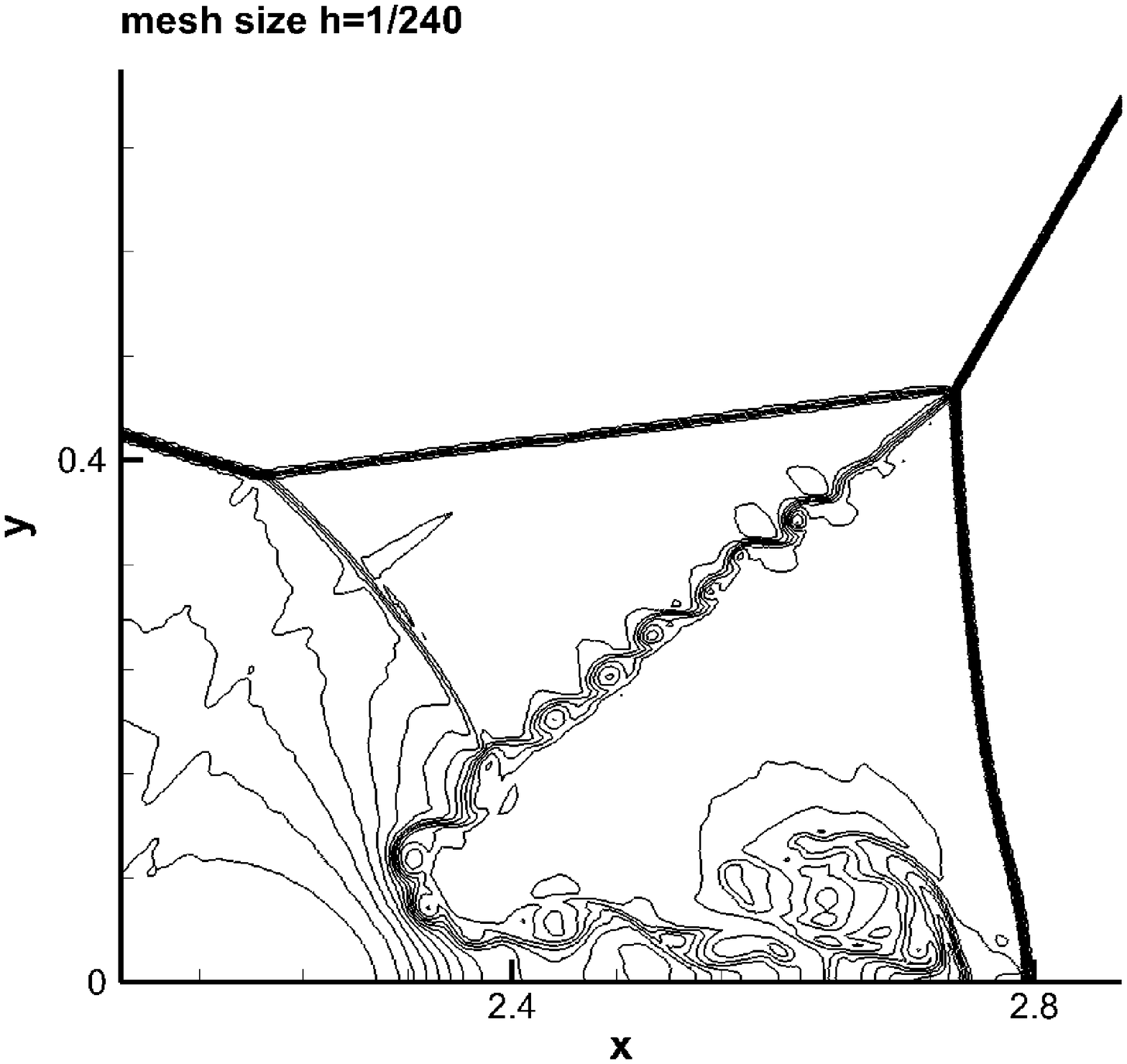}
\caption{\label{2d-doublemach-2} Density enlargement distribution on the triangular mesh with $h=1/120$ and $1/240$ at $t=0.2$.}
\end{figure}

This problem was extensively studied by Woodward and Colella \cite{Case-Woodward} for the inviscid flow. The computational domain is $[0,3.2]\times[0,1]$, and a solid wall lies at the bottom of the computational domain starting from $x =1/6$. Initially a right-moving Mach $10$ shock is positioned at $(x,y)=(1/6, 0)$, and makes a $60^\circ$ angle with the x-axis. The initial pre-shock and post-shock conditions are
\begin{align*}
(\rho, U, V, p)&=(8, 4.125\sqrt{3}, -4.125,
116.5),\\
(\rho, U, V, p)&=(1.4, 0, 0, 1).
\end{align*}
The reflecting boundary condition is used at the wall, while for the rest of bottom boundary, the exact post-shock condition is imposed.
At the top boundary, the flow variables are set to follow  the motion of the Mach $10$ shock.
The density distributions with $h=1/120$ and $1/240$ regular triangular meshes at $t=0.2$ are shown in Fig. \ref{2d-doublemach-1},
and the corresponding local density enlargements are shown in Fig. \ref{2d-doublemach-2}. The current compact fourth-order scheme resolves the flow structure under the triple Mach stem very well.

\subsection{Mach 3 wind tunnel with a step}
\begin{figure}[!htb]
\centering
\includegraphics[width=0.80\textwidth]{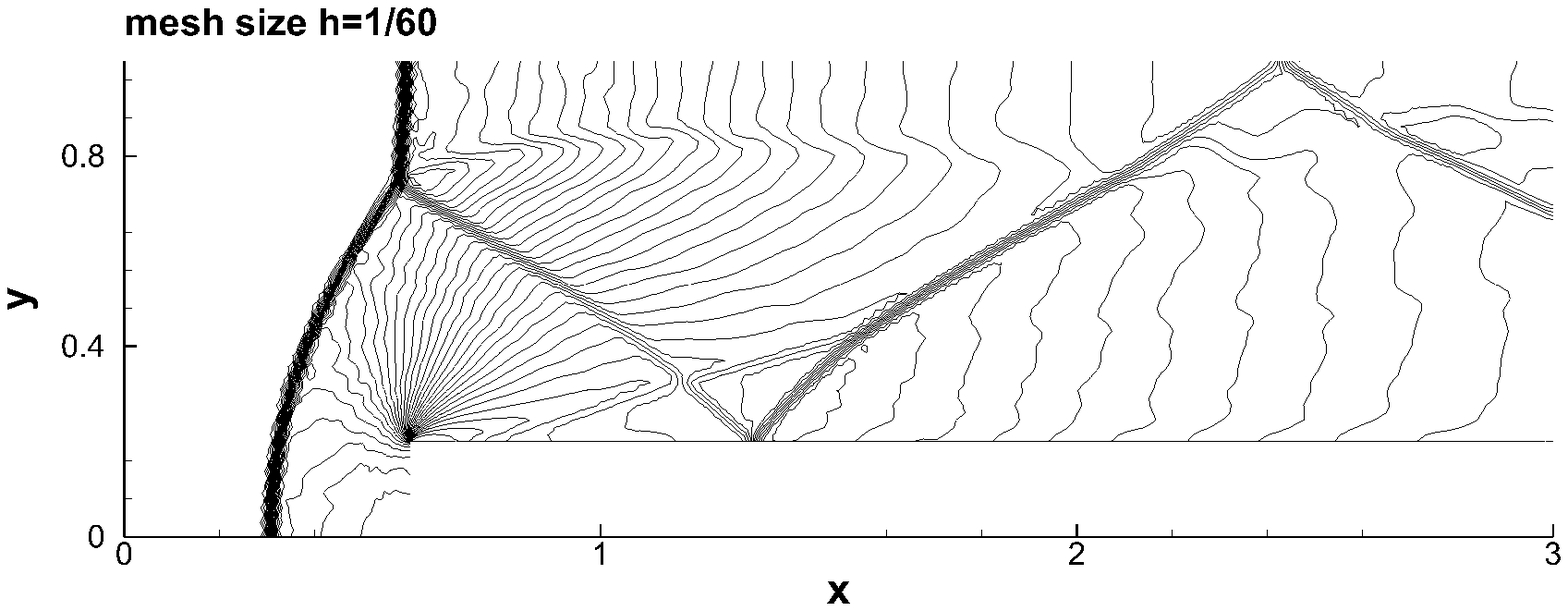}
\includegraphics[width=0.80\textwidth]{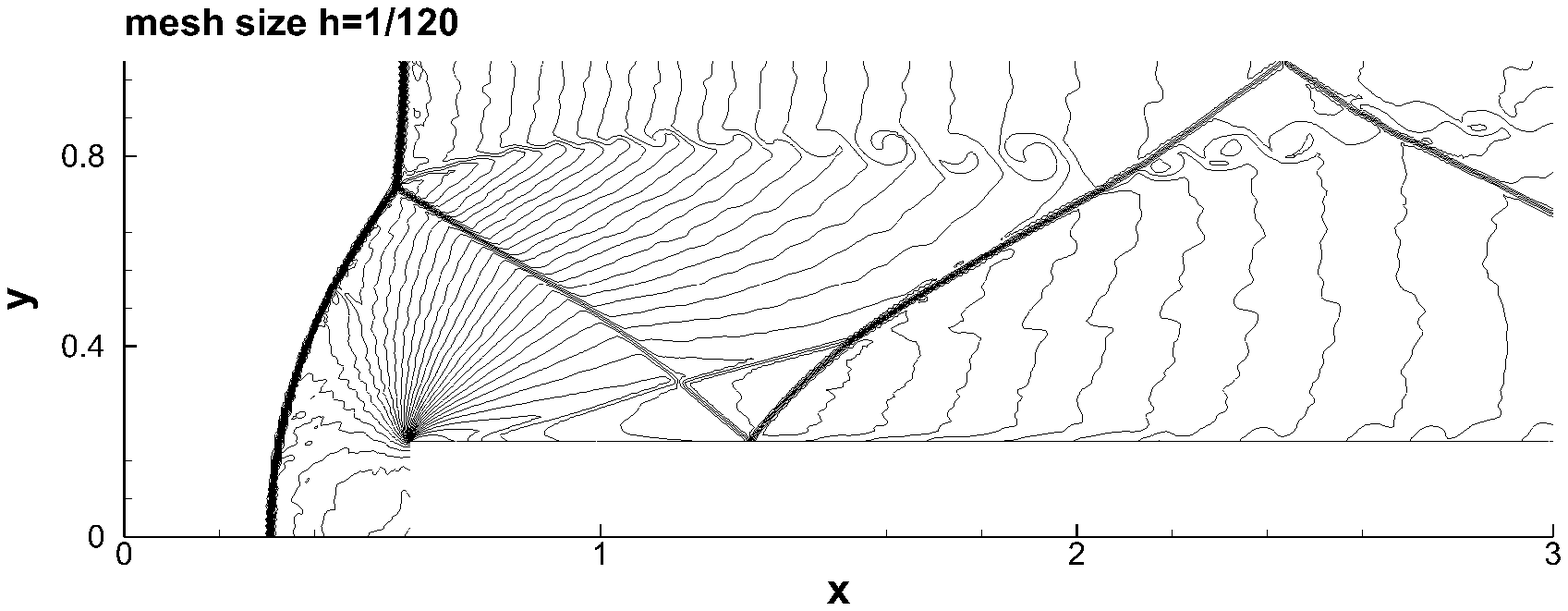}
\caption{\label{2d-tunnel-step} Mach 3 step problem: Density contours at $t=4$ with $h=1/60$ and $h=1/120$ triangular mesh.}
\end{figure}

The step problem was extensively studied in \cite{Case-Woodward} for inviscid flow. The computational domain is $[0,3]\times[0,1] \setminus [0.6,3]\times[0,0.2]$. The height of the wind tunnel is $1$, and the length is $3$. The step is located at $x=0.6$ with height $0.2$ in the tunnel. Initially the tunnel is filled with the gas which elsewhere has $\rho=1, U=3, V=0, p=1/1.4$. The same state is used as the left boundary condition.
 The upper and lower boundaries are wall with slip Euler boundary condition.
The corner of the step is the center of a rarefaction fan. The method of modifying the density and velocity magnitude on the several cells around the corner has not been used in the current computation \cite{Case-Woodward}. The local solution at the corner is properly resolved by the compact GKS on the regular mesh with $h=1/60$ and $1/120$.
The density distribution at time $t=4.0$ is plotted in Fig. \ref{2d-tunnel-step}. One can clearly observe that the current scheme provides a high resolution solution of the physical instability and rolling up of the slip line. This indicates that even with excellent shock capturing capability the scheme has less dissipation in resolving the small scale structure.

\subsection{Hypersonic flow around scramjet configuration}
\begin{figure}[!htb]
\centering
\includegraphics[width=0.85\textwidth]{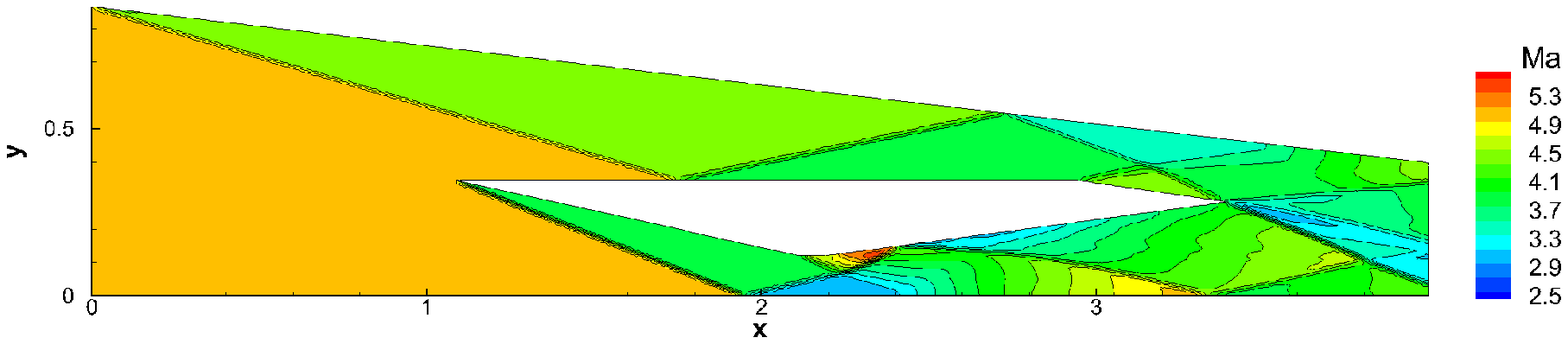}
\includegraphics[width=0.85\textwidth]{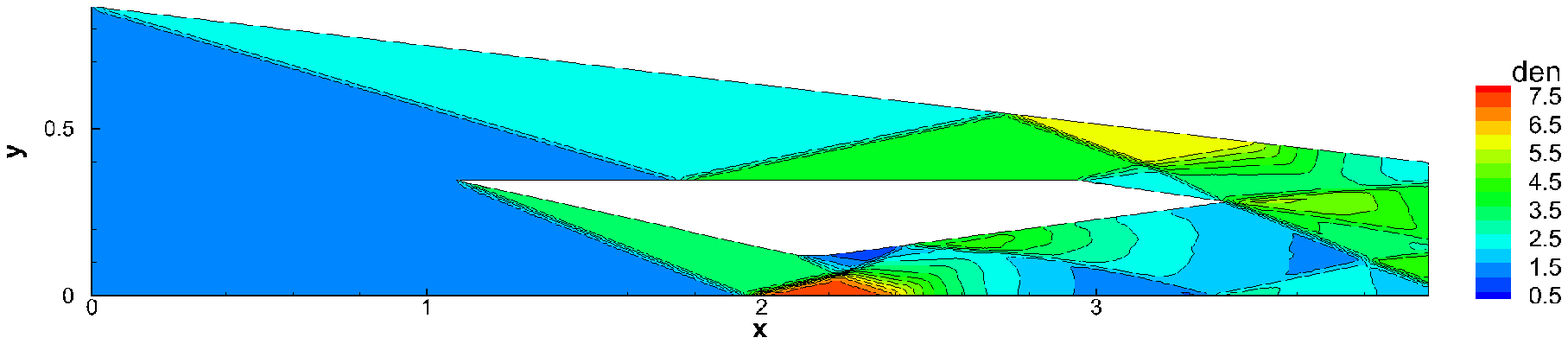}
\includegraphics[width=0.85\textwidth]{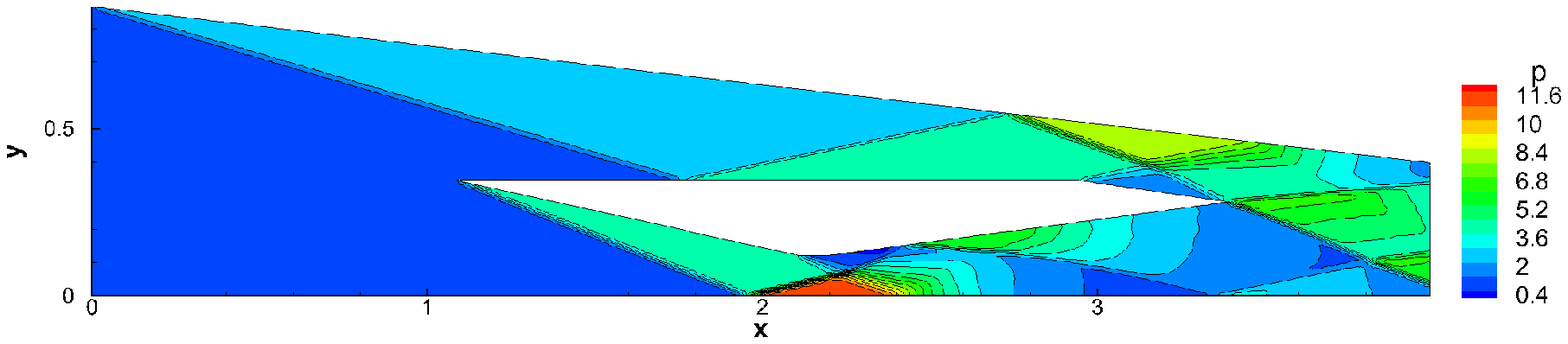}
\caption{\label{2d-scramjet} Hypersonic scramjet flow: Mach number, density and pressure contours with $h=1/100$ triangular mesh.}
\end{figure}

In order to show the robustness of the scheme in capturing discontinuous solution, the scramjet configuration was computed at $Ma=5$ \cite{ollivier1997_scramjet}. A regular triangular mesh with $52204$ cells was generated in the current computation. Fig. \ref{2d-scramjet} shows the solutions of Mach number, density, and pressure from the compact fourth-order GKS. Much of the details of flow field are resolved, such as the rarefaction wave near the throat and the contact slip layer near the outlet.

\subsection{Hypersonic flow around an inviscid cylinder}
\begin{figure}[!htb]
	\centering
	\includegraphics[width=0.235\textwidth]{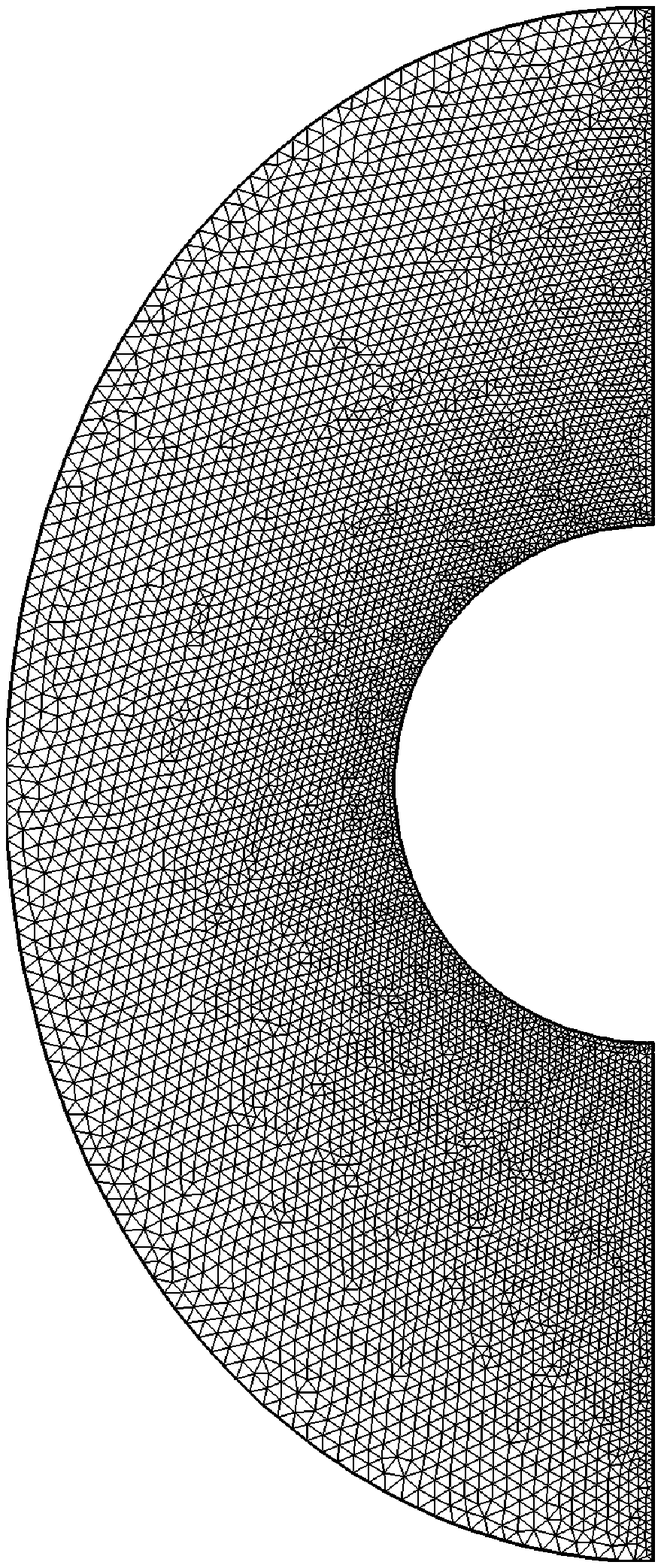}
	\includegraphics[width=0.235\textwidth]{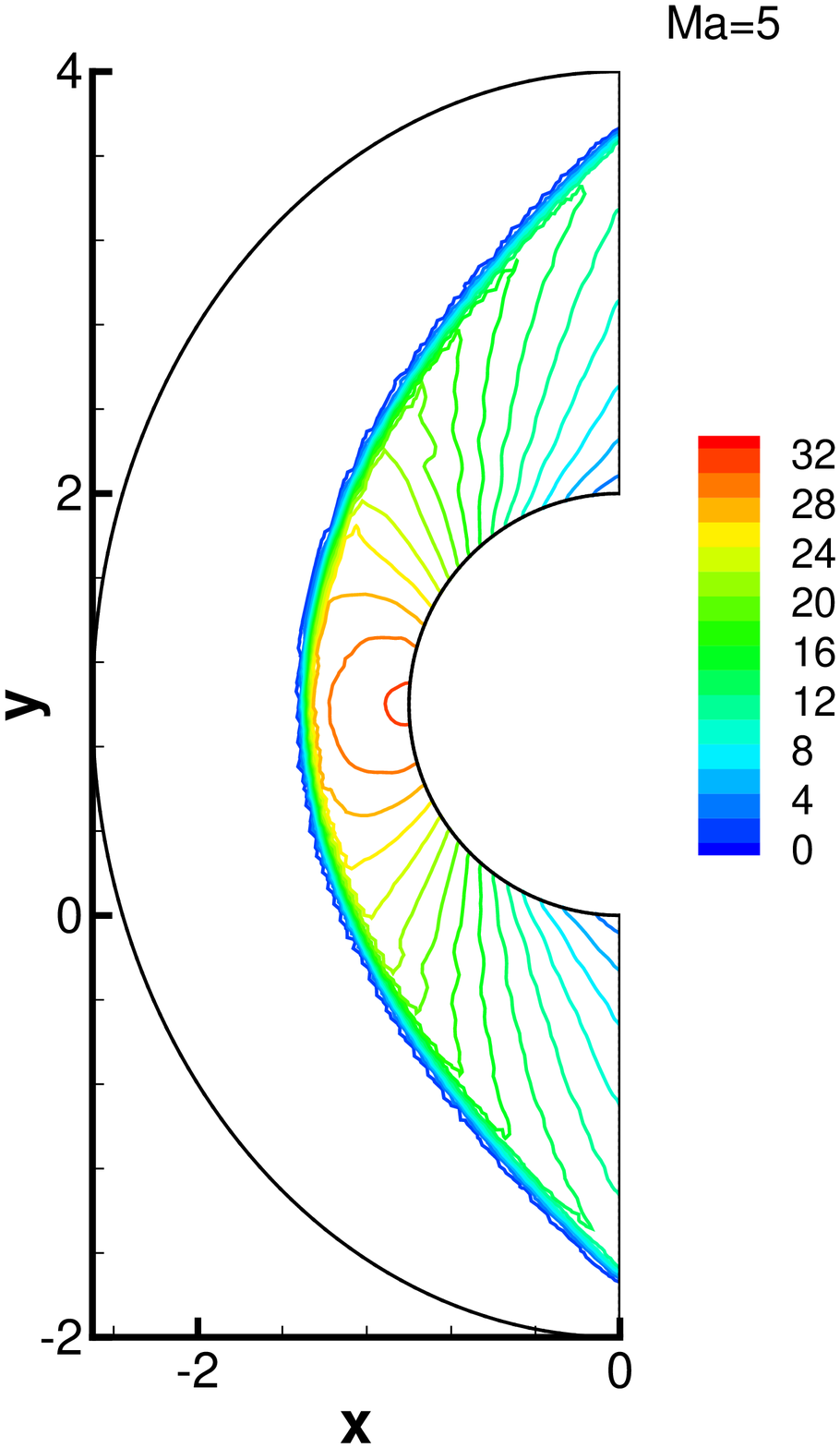}
	\includegraphics[width=0.235\textwidth]{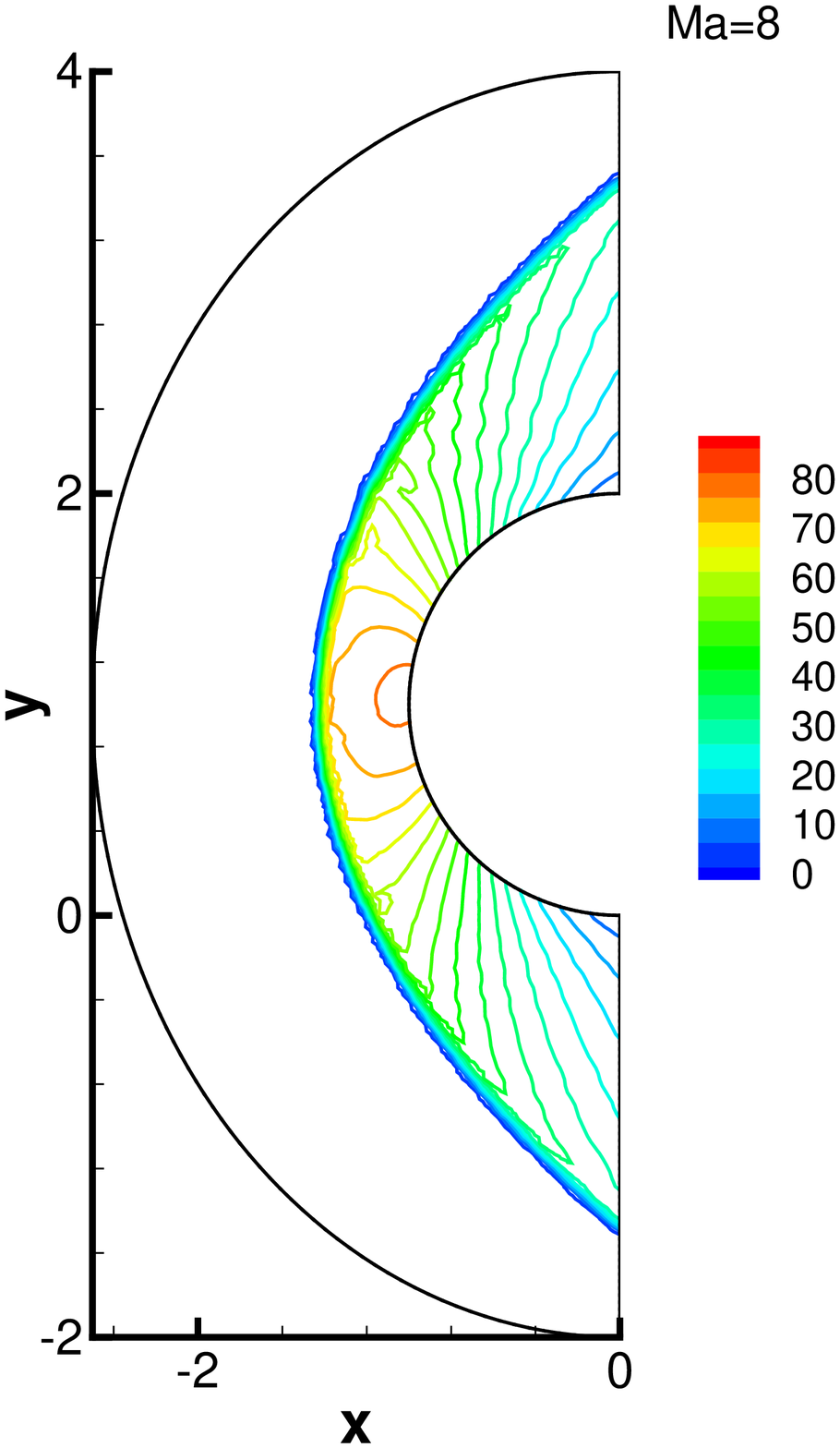}
	\includegraphics[width=0.235\textwidth]{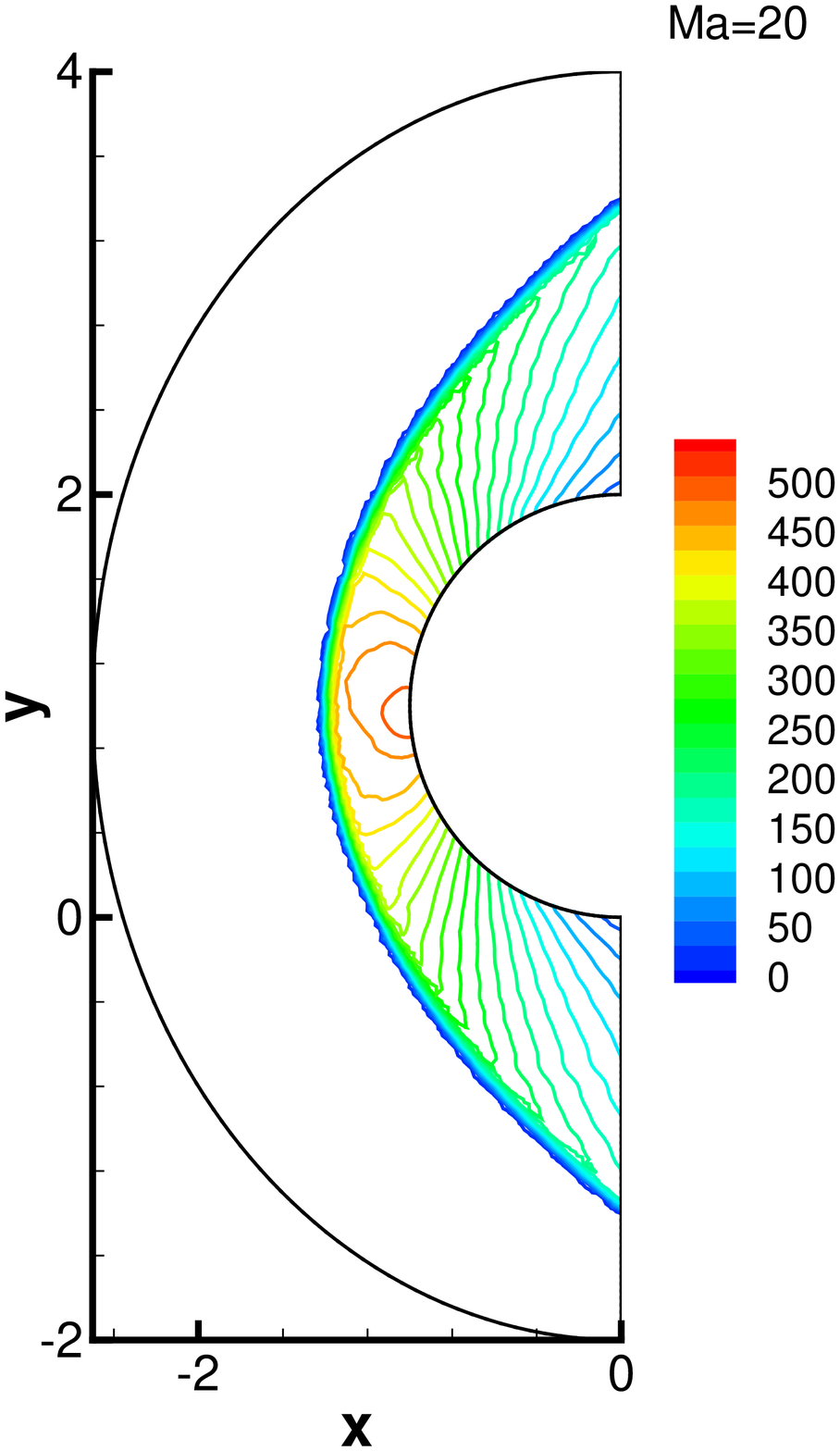}
	\caption{\label{2d-inv-cylinder} Hypersonic inviscid flow passing through a circular cylinder: pressure distributions with Mach number $Ma =5, 8$ and $20$ on triangular mesh.}
\end{figure}

The incoming flow has Mach numbers up to $20$. The adiabatic reflective boundary condition is imposed on the wall of the cylinder while the right boundary is set as outflow boundary condition. The mesh and pressure distributions are presented in Fig. \ref{2d-inv-cylinder}. The regular triangular mesh is used, and the mesh is refined near the cylinder. The results agree well with those calculated on structured mesh by the non-compact high-order GKS. For the case of $Ma=20$, the $CFL$ takes $0.2$.

\subsection{Viscous shock tube}
\begin{figure}[!htb]
	\centering
	\includegraphics[width=0.495\textwidth]{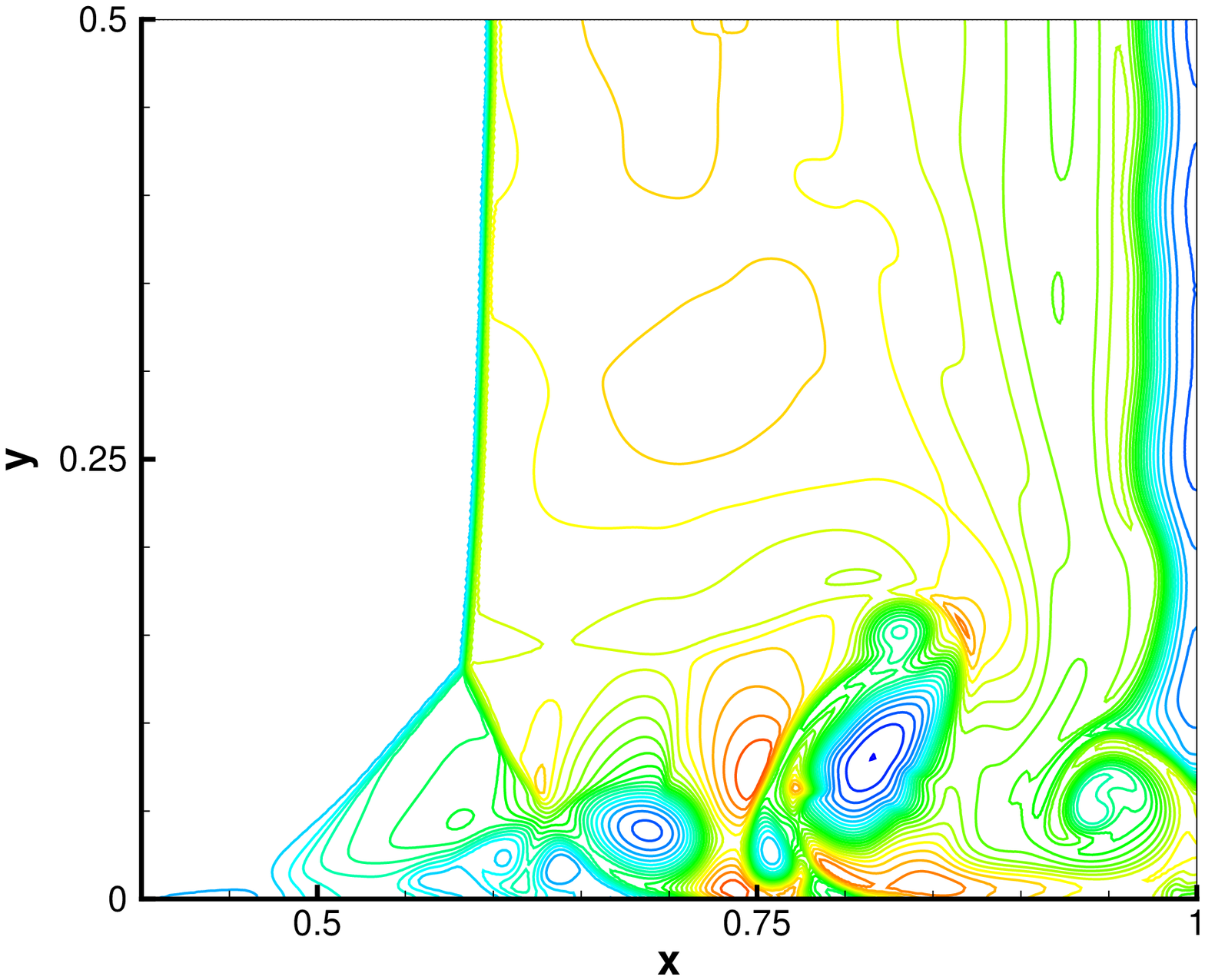}
	\includegraphics[width=0.495\textwidth]{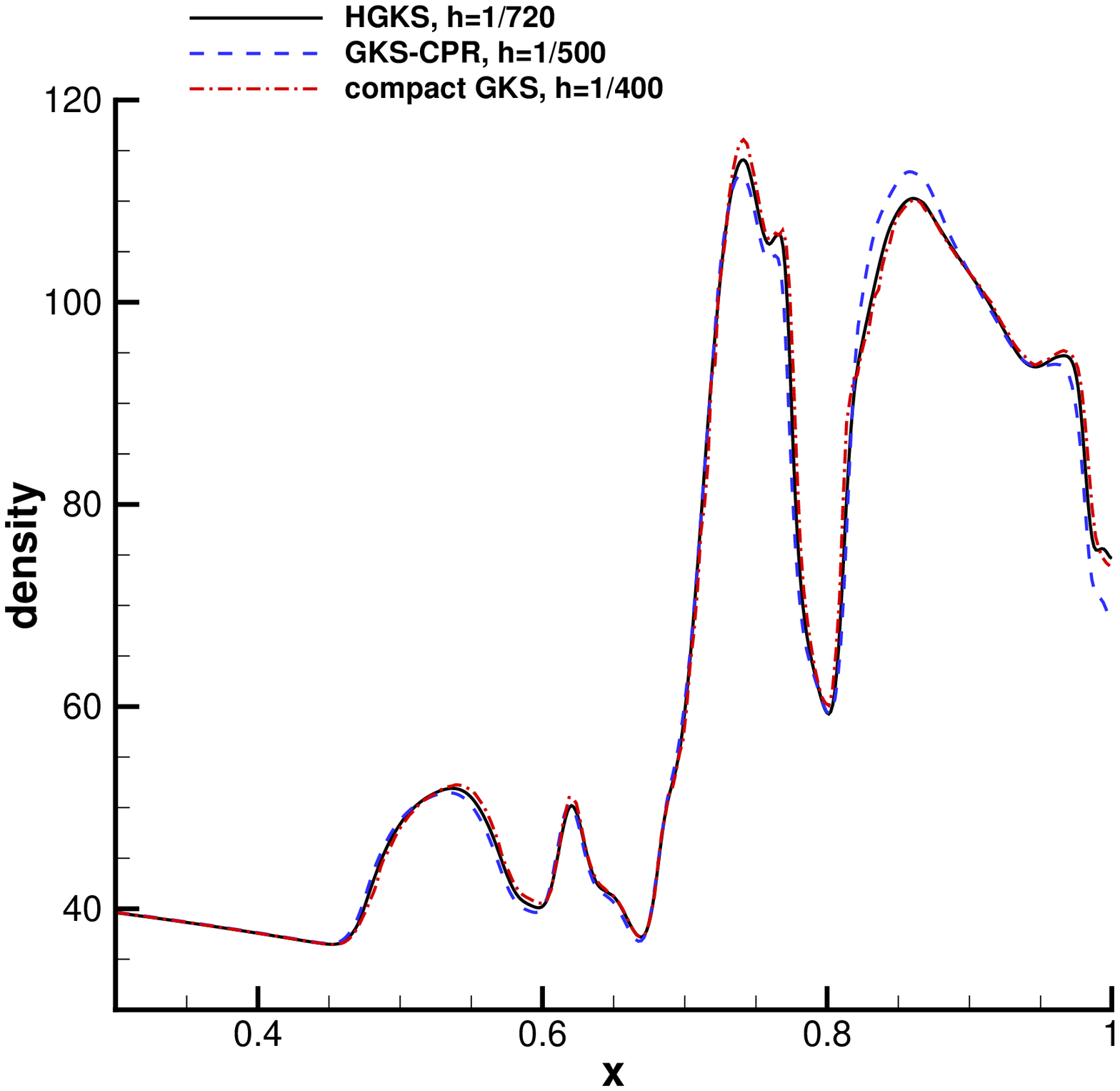}
	\caption{\label{2d-vis-shocktube-2} Viscous shock tube: density field (left) and the density distribution along the lower wall (right) with the mesh $h=1/400$ at $t=1.0$. $24$ uniform contours are presented ranging from $20$ to $120$ for the density field.}
\end{figure}

Viscous shock tube problem is tested to validate the compact fourth-order GKS for the complicated shock wave and boundary layer interaction.
This problem requires not only the robustness of the scheme, but also the accuracy of the numerical method. The flow is bounded in a unit square cavity. The computational domain is set as $[0,1] \times [0,0.5]$ and  symmetrical condition is used on the top boundary. The non-slip and adiabatic wall conditions are imposed on other boundaries. The initial condition is
\begin{equation*}
(\rho,U,V,p) = \begin{cases}
(120,0,0,120/\gamma),  0\leq x<0.5,\\
(1.2,0,0,1.2/\gamma),  0.5\leq x\leq1.
\end{cases}
\end{equation*}
The viscosity coefficient is $\mu=0.005$ with a corresponding Reynolds number $Re=200$.
The Prandtl number in the current computation is set to be $P_r=1$.
Initially, the shock wave, followed by a contact discontinuity, moves towards to the right wall. A thin boundary layer is created above the lower wall. The complex shock and boundary layer interaction occurs and results in a lambda-shape shock pattern after the reflecting shock wave from the right wall.
The density field at $t=1$ is presented in Fig \ref{2d-vis-shocktube-2}. The complex flow structure, including the lambda shock and the vortex configurations, are well resolved by the current compact GKS with a mesh size $h=1/400$. A quantitative verification for the result is also given. The density distribution along the lower wall is presented. The result from non-compact high-order GKS \cite{liQB2010high} with $h=1/720$ structured mesh is used as the reference solution. The result of third-order CPR-GKS scheme \cite{CPR-GKS} with $h=1/500$ triangular mesh is plotted as well. The compact GKS has a better resolution on the coarse mesh. Compared with the third-order CPR-GKS, there are no additional DOFs inside each cell and trouble cell detection in the current compact GKS. This test case demonstrates that the current compact GKS can resolve the flow with shock and boundary layer interaction robustly and accurately.

\begin{figure}[!htb]
	\centering
	\includegraphics[width=0.80\textwidth]{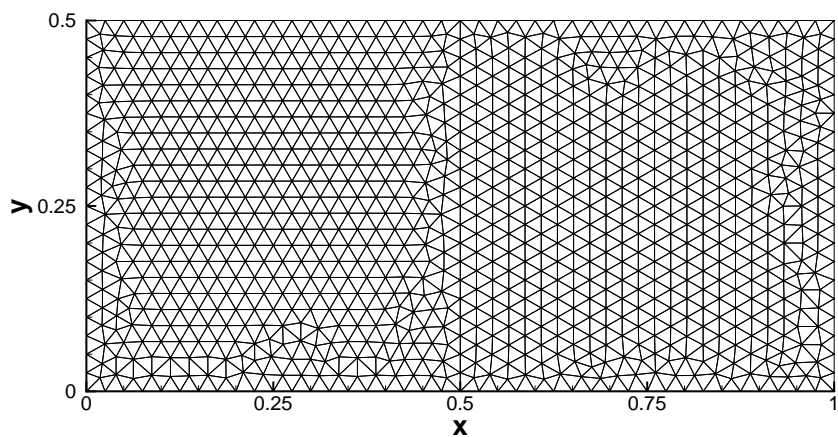}
	\caption{\label{2d-vis-shocktube-1} Viscous shock tube: mesh distribution with $h=1/40$.}
\end{figure}

\subsection{Laminar boundary layer flow}
\begin{figure}[!htb]
\centering
\includegraphics[width=0.670\textwidth]{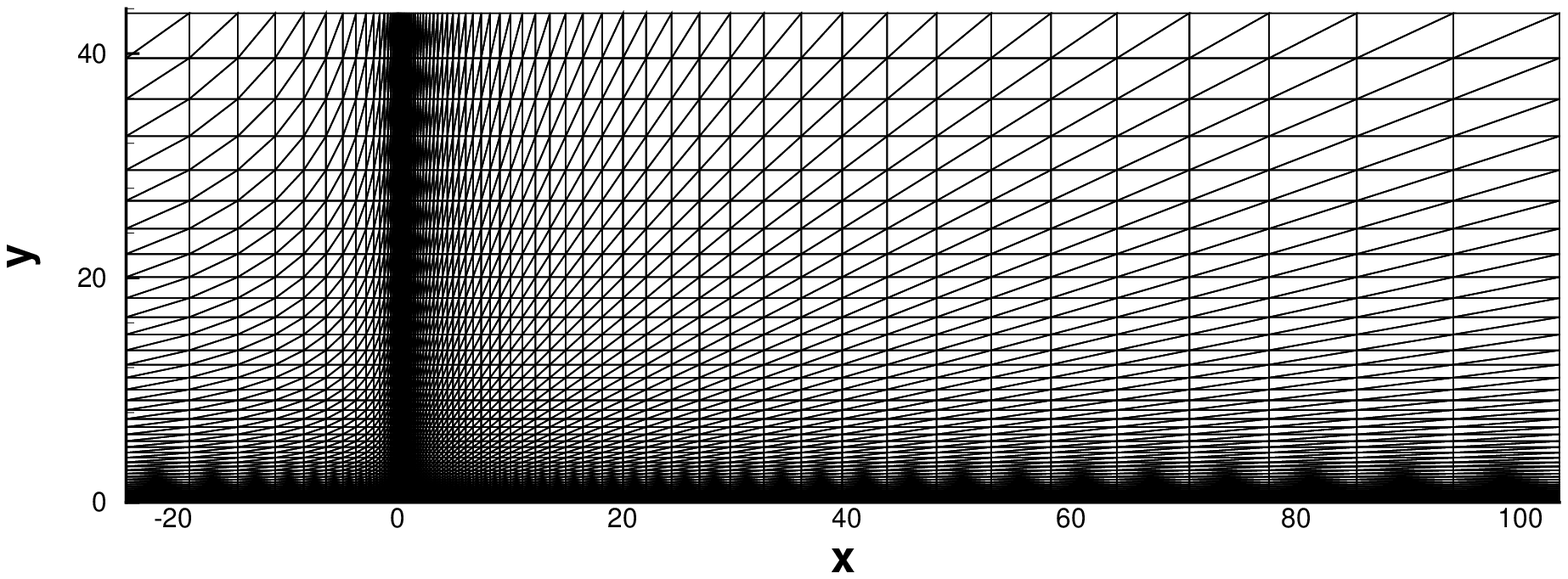}
\includegraphics[width=0.285\textwidth]{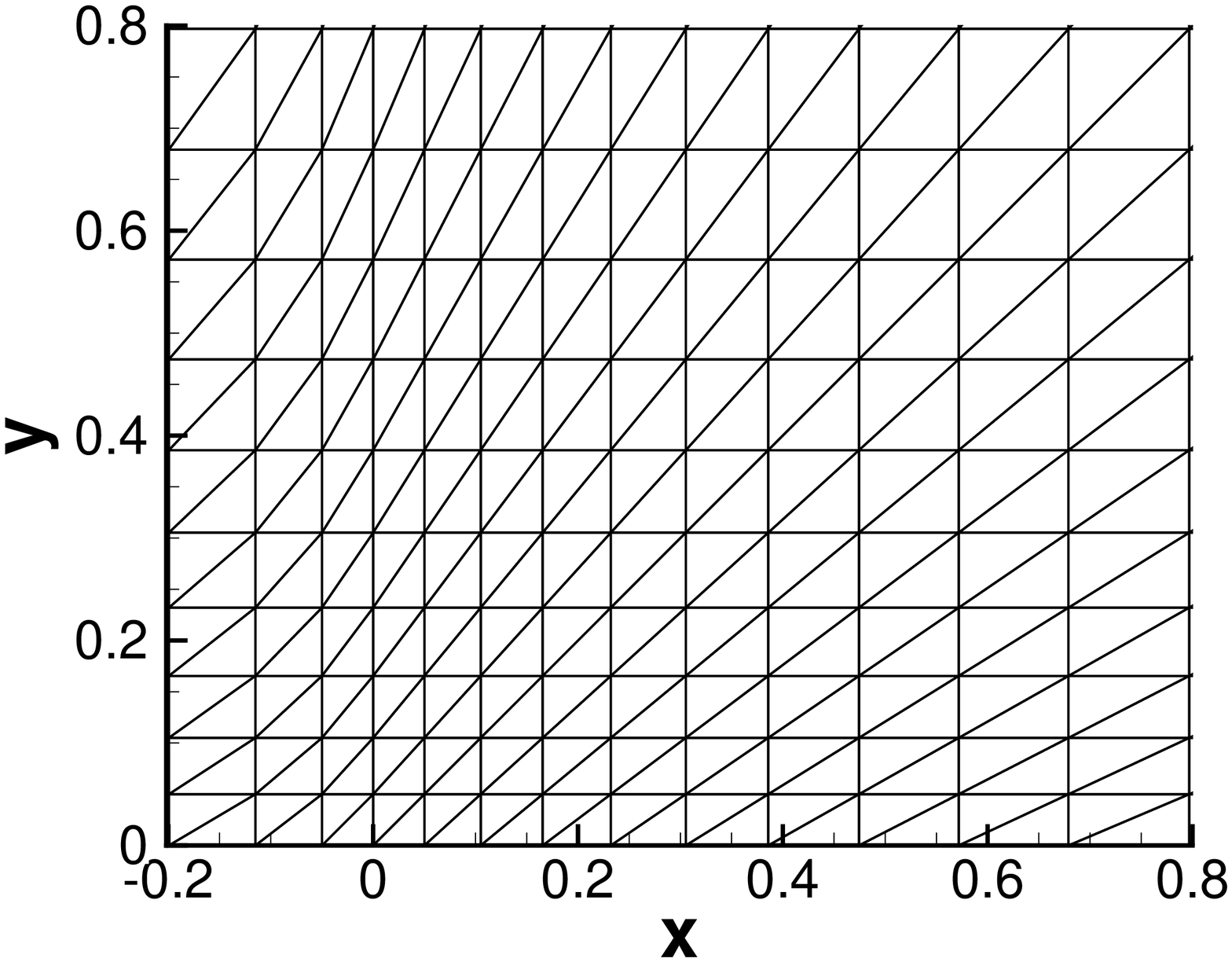}
\caption{\label{2d-boundary-layer-1} Laminar boundary layer flow: the computational domain and mesh. The minimum size of cells along streamwise and transverse directions is $h=0.05$, and the maximum aspect ratio of the triangle in the mesh is about $190$.}
\end{figure}

\begin{figure}[!htb]
\centering
\includegraphics[width=0.45\textwidth]{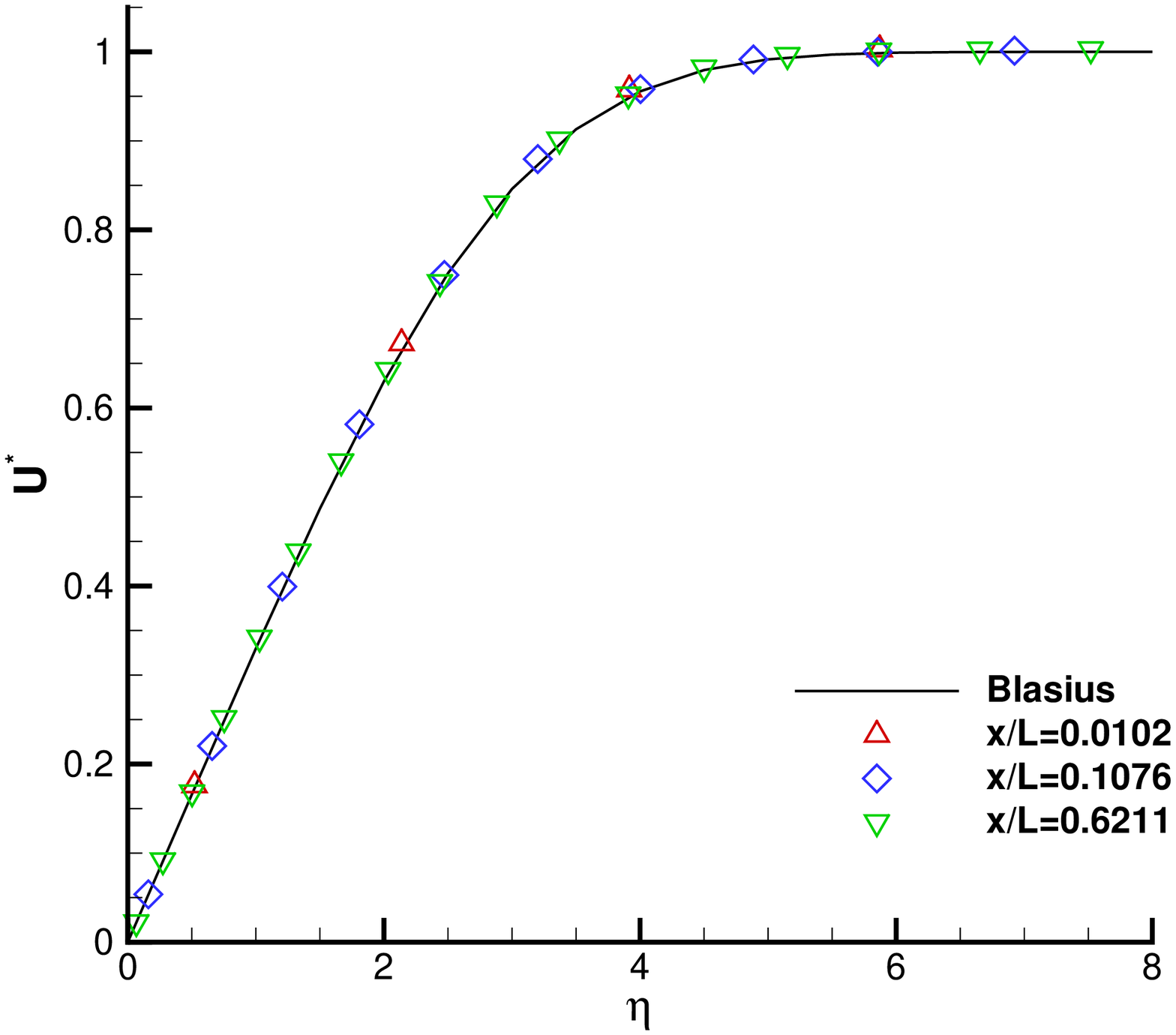}
\includegraphics[width=0.45\textwidth]{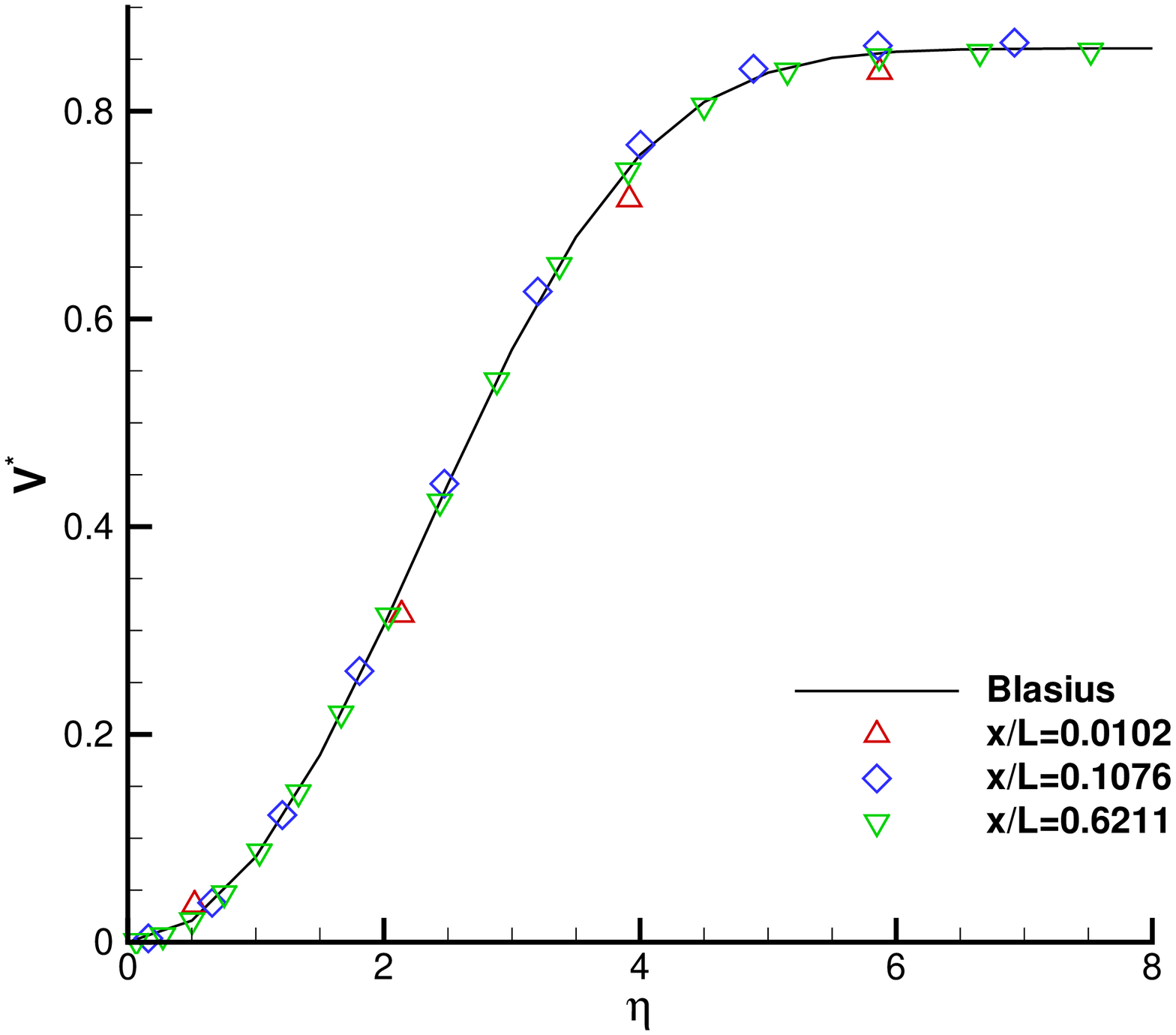}
\caption{\label{2d-boundary-layer-2} Laminar boundary layer flow: streamwise and transverse velocity profiles at different locations in flat plate boundary layer.}
\end{figure}

\begin{figure}[!htb]
	\centering
	\includegraphics[width=0.55\textwidth]{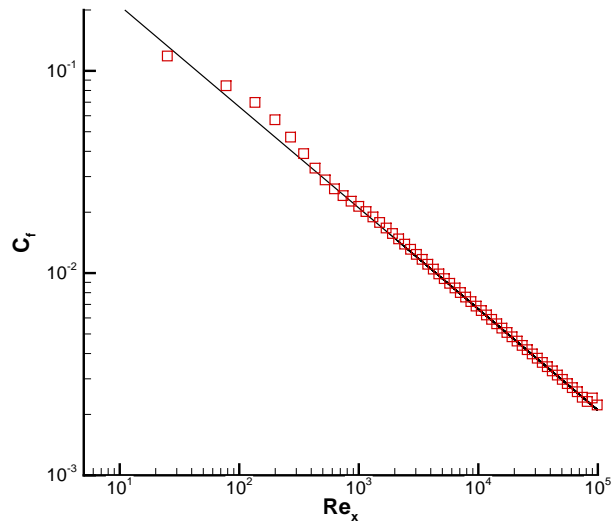}
	\caption{\label{2d-boundary-layer-3} Laminar boundary layer flow: wall friction coefficient distribution along the flat plate.}
\end{figure}

The laminar boundary layer flow problem is tested to validate the compact GKS for viscous flow, especially when the computational mesh is strongly stretched in one direction with a large aspect ratio. The Mach number of the free stream is $Ma_{\infty}=0.15$, and the Reynolds number is $Re_{\infty}=1.0\times10^5$. The Reynolds number is defined as $Re_{\infty}=U_{\infty}L/\nu$, where $L=100$ is the length of plate and $U_{\infty}$ is the velocity of free stream. The density of the free stream is $\rho_{\infty}=1.0$, and the pressure is $p_{\infty}=1.0/1.4$.
The computational domain and mesh are shown in Fig. \ref{2d-boundary-layer-1}. The plate is placed at $y=0$ and $x>0$.
A total $75\times47\times2$ cells are used in the domain. The minimum size of cells along streamwise and transverse directions is $h=0.05$. The mesh is generated from $(0,0)$ with a stretching rate of $1.1$ along the positive x-direction, $1.3$ along the negative x-direction, and $1.1$ along the positive y-direction. At the trailing edge, the maximum aspect ratio of the triangle in the current mesh is about $190$.
The streamwise and transverse velocity profiles at different locations in the flat plate boundary layer are shown in Fig. \ref{2d-boundary-layer-2}. The numerical results agree well with the Blasius solution at the locations near and far from the leading edge. The wall friction coefficient distribution along the flat plate is given in Fig. \ref{2d-boundary-layer-3}. The accurate wall friction coefficient distribution is obtained along the flat plate except for a few mesh points close to the leading edge.

\subsection{Lid-driven cavity flow}
\begin{figure}[!htb]
\centering
\includegraphics[width=0.45\textwidth]{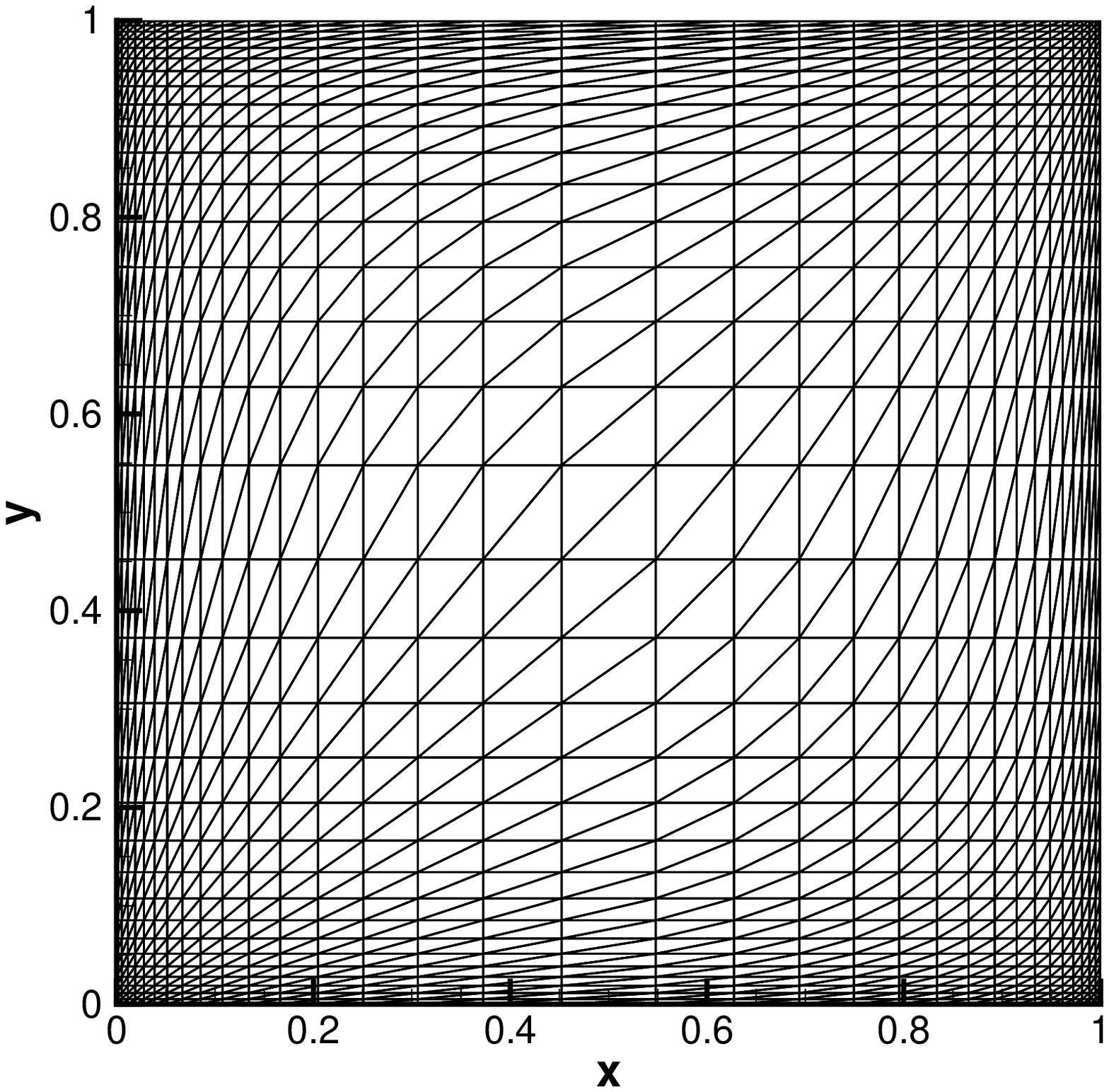}
\includegraphics[width=0.45\textwidth]{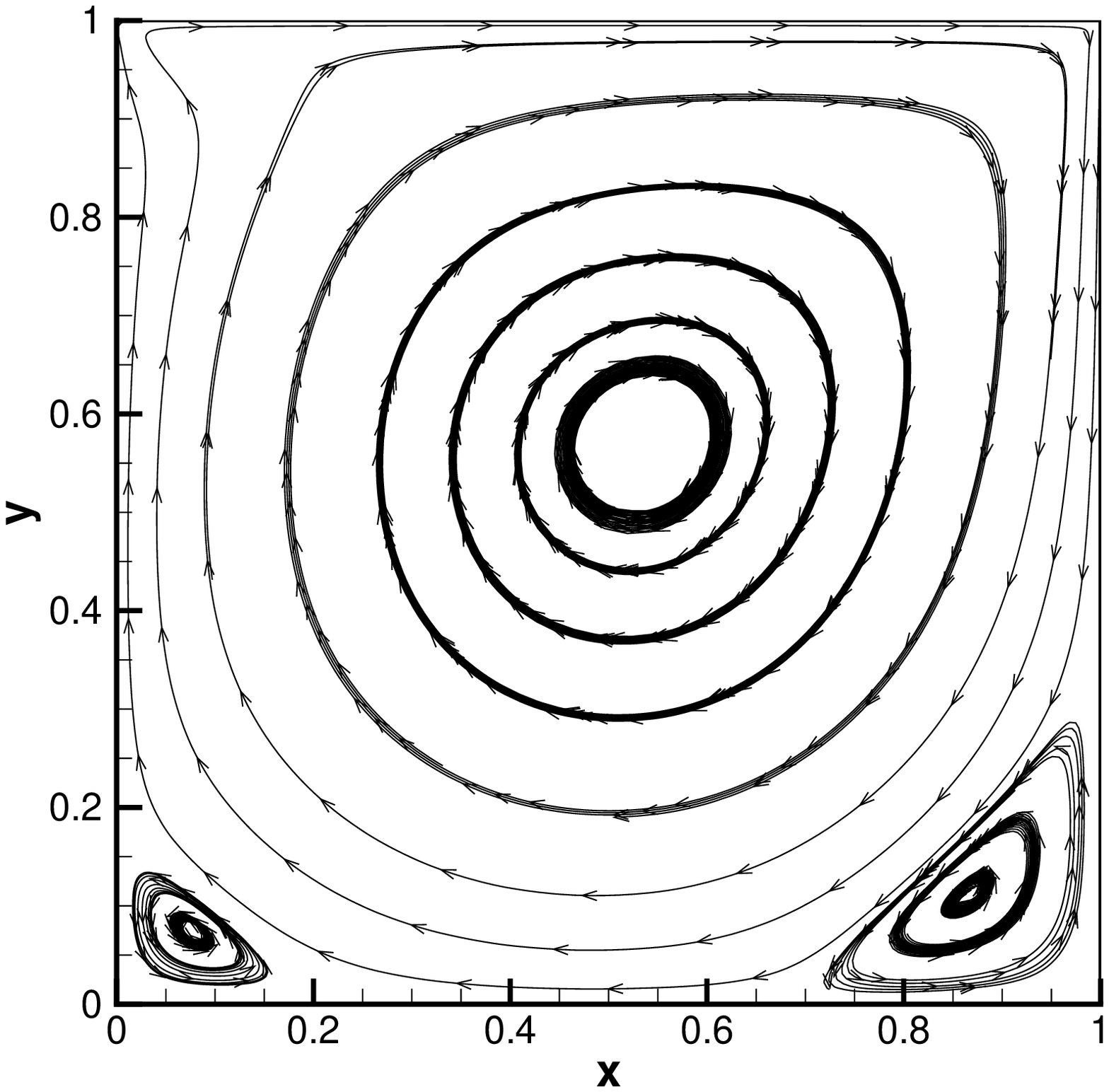}
\caption{\label{2d-cavity-1} Lid-driven cavity flow: the computational mesh (left) with total $33\times33\times 2$ mesh cells and streamlines for $Re=1000$ flow (right). }
\end{figure}

\begin{figure}[!htb]
\centering
\includegraphics[width=0.45\textwidth]{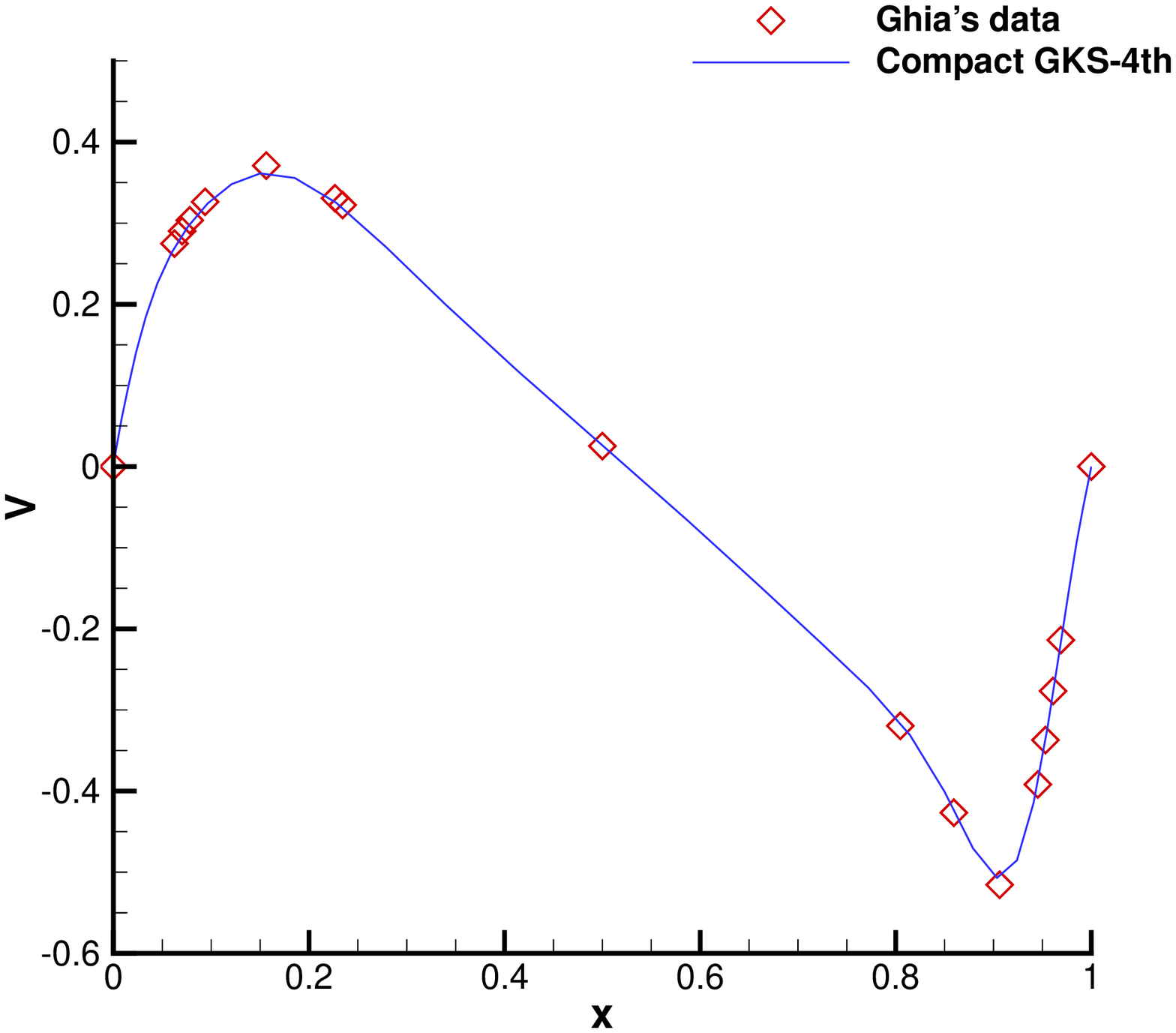}
\includegraphics[width=0.45\textwidth]{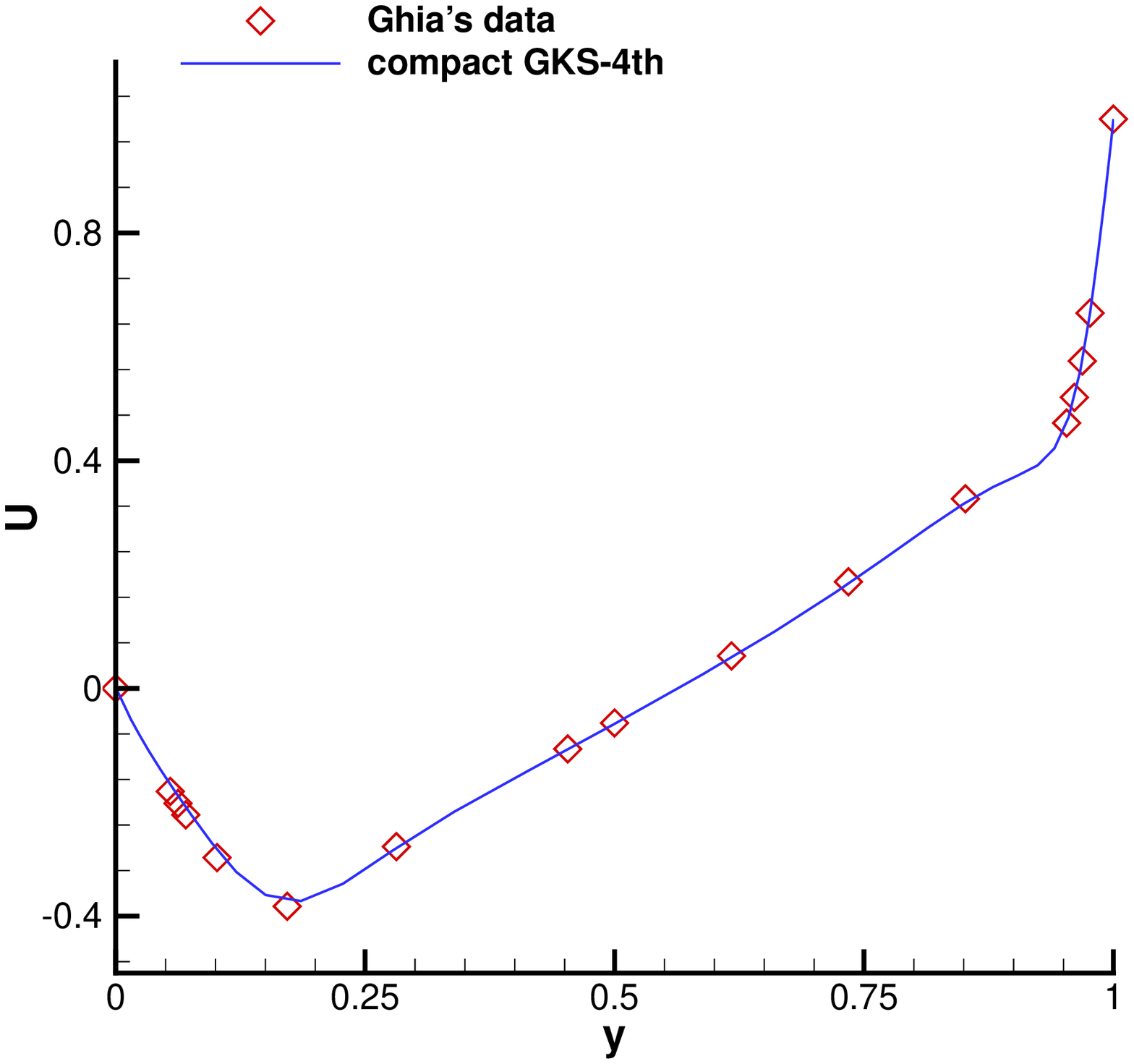}
\caption{\label{2d-cavity-2} Lid-driven cavity flow: V-velocities along the horizontal centerline and U-velocities along the vertical centerline with Re =1000. }
\end{figure}

The lid-driven cavity problem is a test case for incompressible viscous flow. The fluid is bounded in a unit square where the top boundary is  moving with a uniform speed $U_0=1$ and temperature $T_0=1$.
The corresponding Mach number is $Ma=U_0/\sqrt{\gamma R T_0}=0.15$. The non-slip and isothermal boundary conditions are imposed on all boundaries with the wall temperature $T_w=T_0$. The initial flow is stationary with density $\rho_1=\rho_0$ and temperature $T_1=T_0$. The computational domain is $[0,1]\times[0,1]$. The case of $Re=1000$ is tested.
The computational mesh and the streamlines are shown in Fig. \ref{2d-cavity-1}. Total $33\times33\times 2$ mesh cells are used, and the stretching rate in two directions is $1.2$. The mesh is refined close the wall and the minimum mesh size is about $h=0.0052$.
The velocities distribution along horizontal and vertical center lines are shown in Fig. \ref{2d-cavity-2}. The numerical results agree well with the reference solutions.

\subsection{Subsonic laminar flow over a circular cylinder}
\begin{figure}[!htb]
	\centering
	\includegraphics[width=0.495\textwidth]{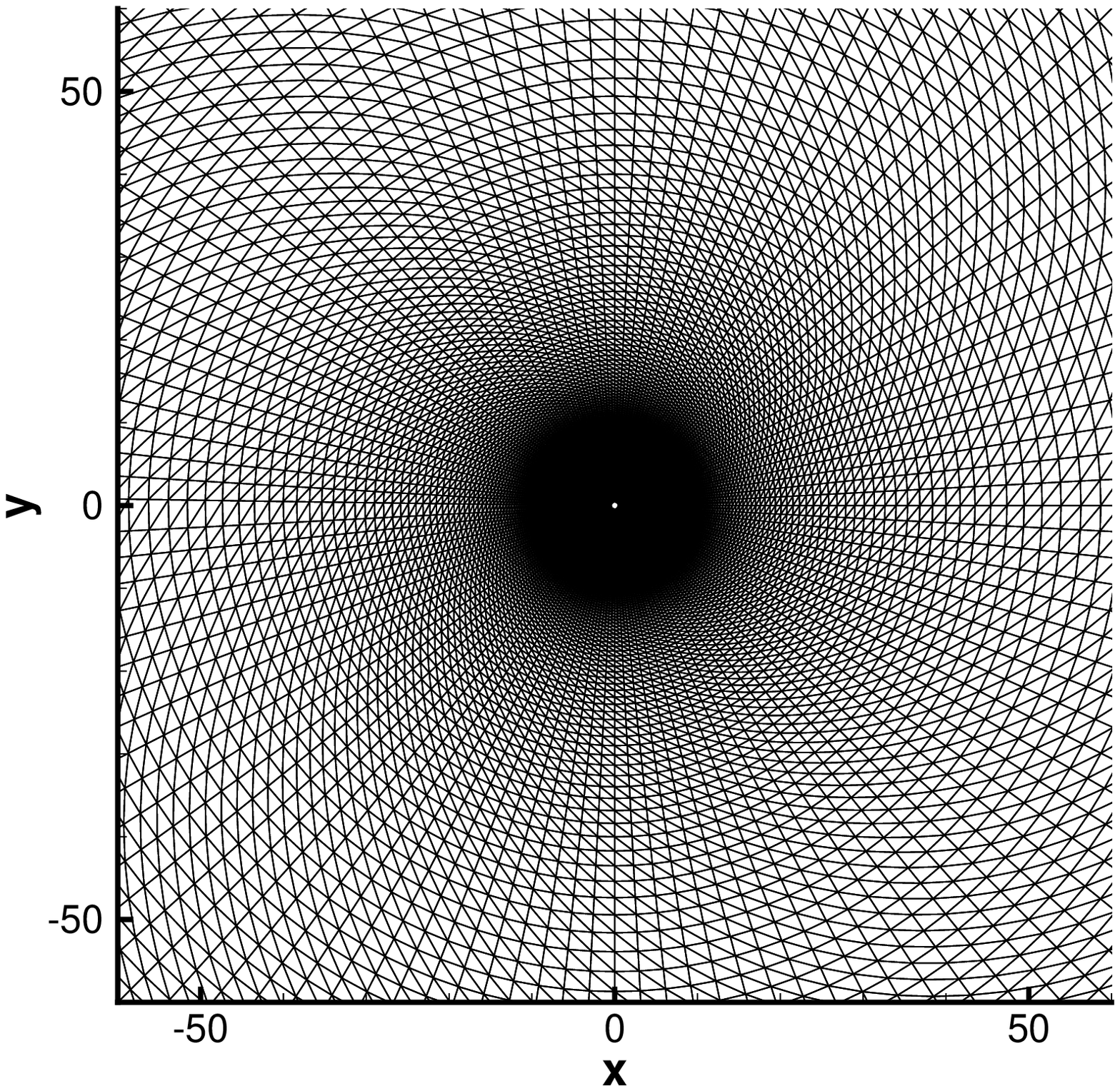}
	\includegraphics[width=0.495\textwidth]{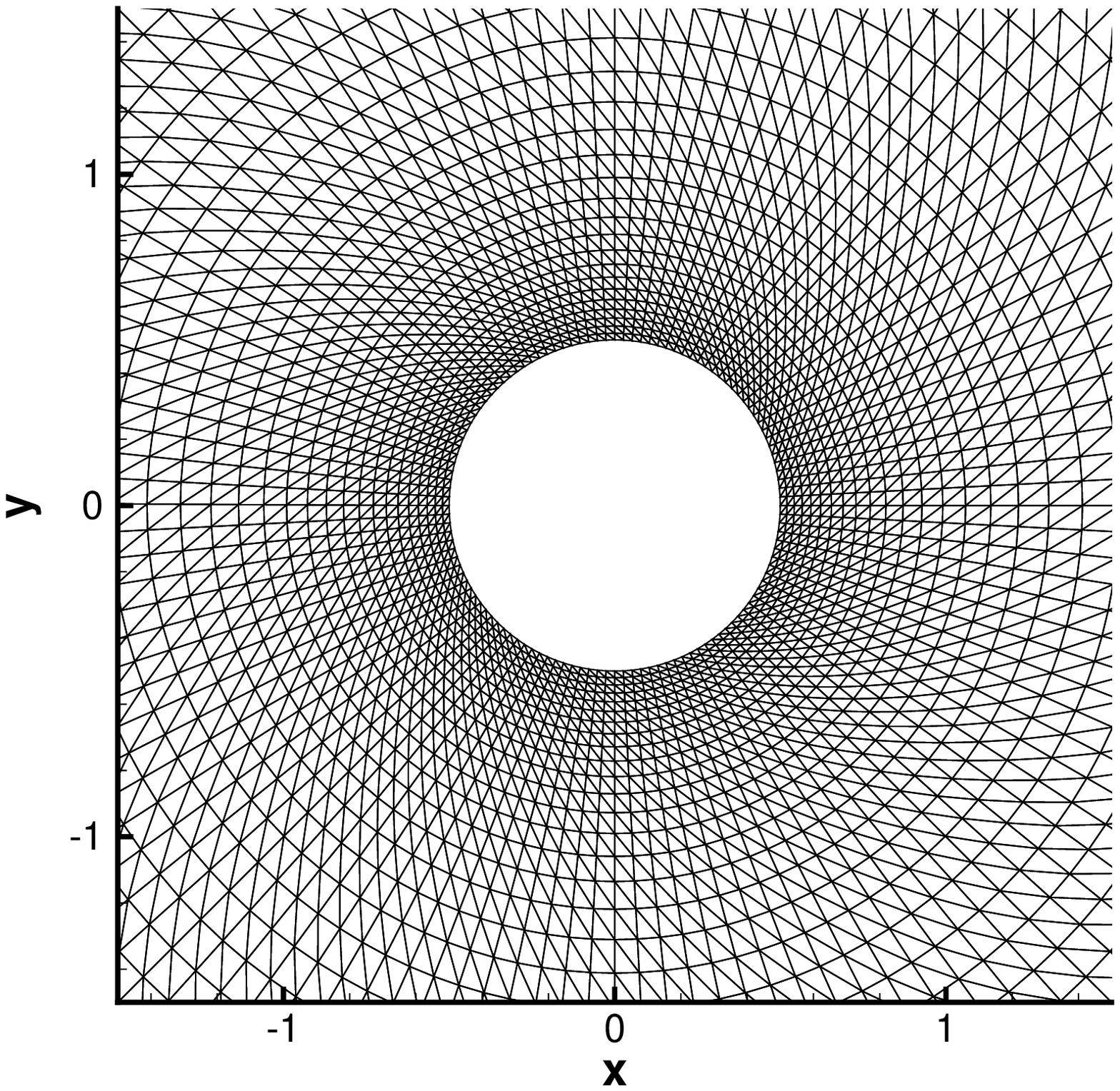}
	\caption{\label{2d-subcylinder-mesh} The computational mesh for the laminar Flow over a circular cylinder. The right figure is the local enlargement of the mesh. The radial mesh size of the first layer above the wall is $2.0 \times 10^{-2}$, and the number of circumferential mesh is $120$. Starting at a distance of $50$ from the center of the cylinder, the mesh size along the radial direction is about $2.5$. }
\end{figure}

The laminar flow over a circular cylinder is studied in this section to validates the compact GKS on triangular mesh for unsteady viscous flow. The incoming free streaming flow is uniform with a Mach number of $Ma=0.2$ and a Reynolds number $Re=150$ based on diameter of the circular cylinder  $D=1$. At this Reynolds number, the flow is essentially two dimensional with periodic shedding vortex pairs from the downstream side of the cylinder. The drag coefficient $C_D$ and Strouhal number $St$ are quantitatively compared.
A triangular mesh in Fig. \ref{2d-subcylinder-mesh} is used. The outer boundary of the computational domain is a circle with diameter of about $346$. The radial size of the first layer mesh on the wall is $2.0 \times 10^{-2}$, and the stretching rate of the cell size in the radial direction is $1.1$ in the region near the cylinder. The number of circumferential mesh points is $120$. The total number of mesh points is $44640$.

\begin{figure}[!htb]
	\centering
	\includegraphics[width=0.45\textwidth]{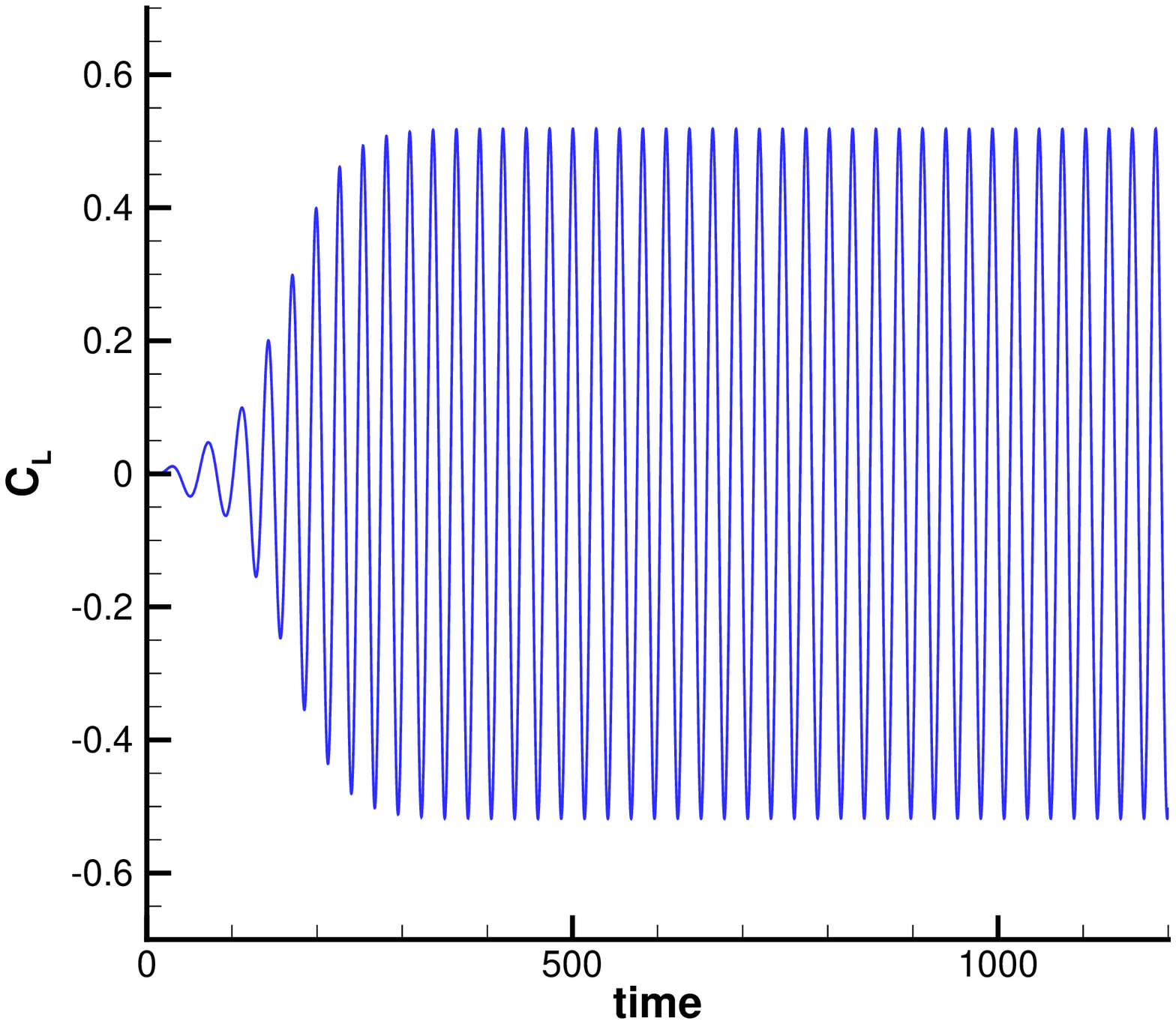}
	\includegraphics[width=0.45\textwidth]{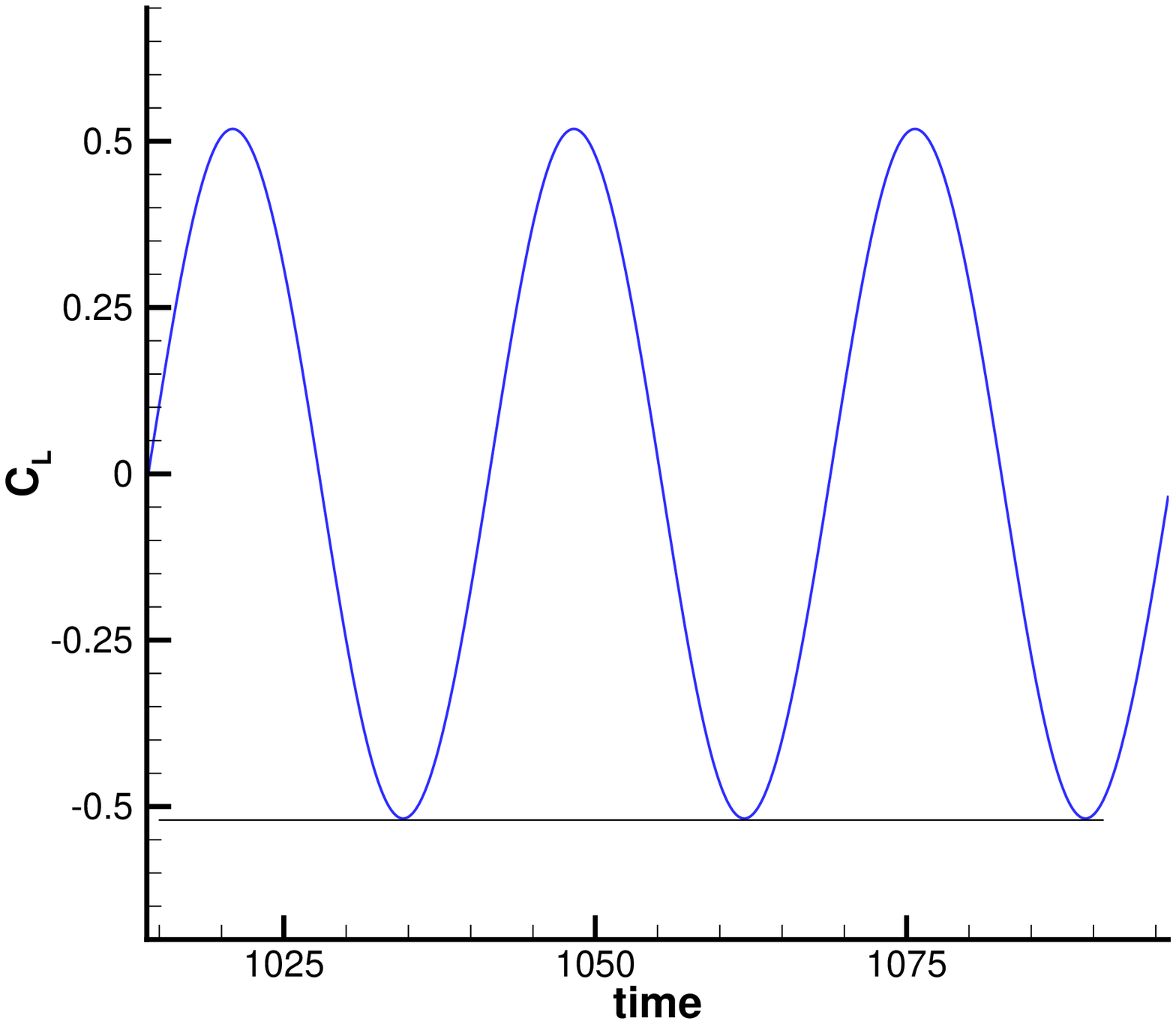}\\
	\includegraphics[width=0.45\textwidth]{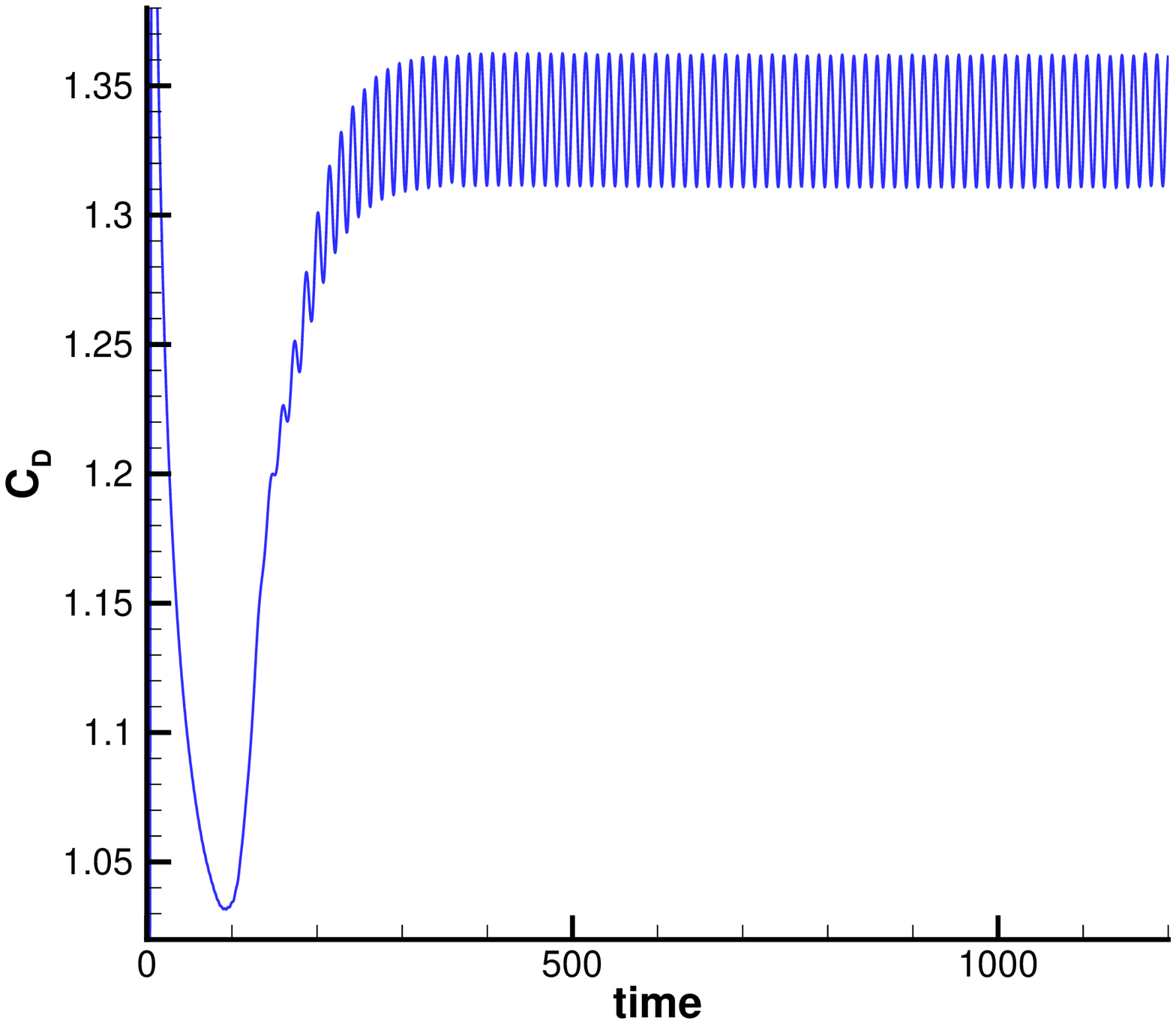}
	\includegraphics[width=0.45\textwidth]{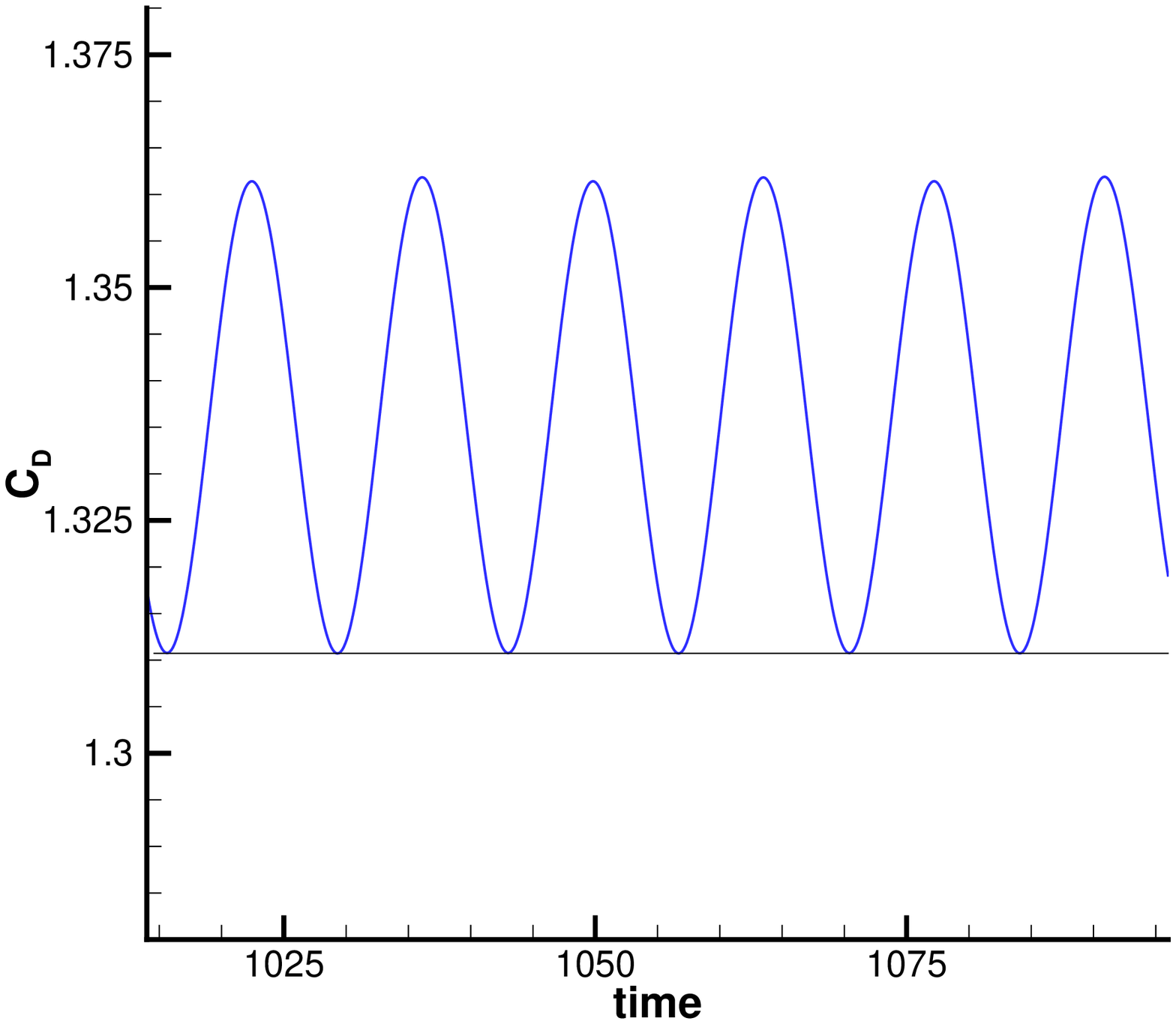}
	\caption{\label{2d-subcylinder-clcd} Total lift coefficient.}
\end{figure}

\begin{figure}[!htb]
	\centering
	\includegraphics[width=0.85\textwidth]{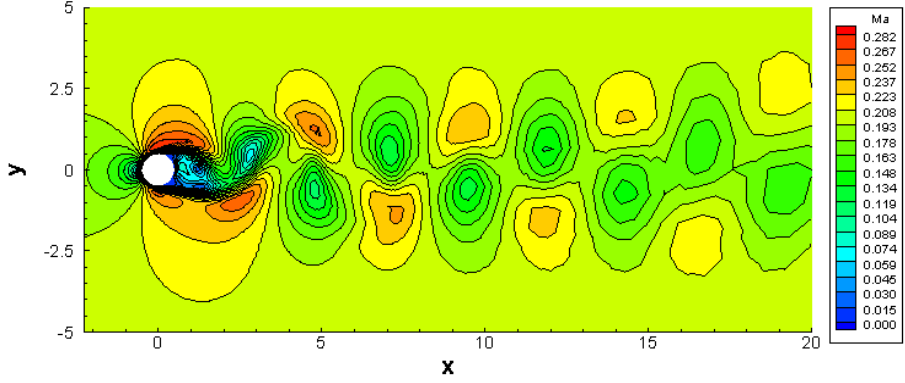}
	\caption{\label{2d-subcylinder-Ma} The contours of Mach number.}
\end{figure}

\begin{figure}[!htb]
	\centering
	\includegraphics[width=0.60\textwidth]{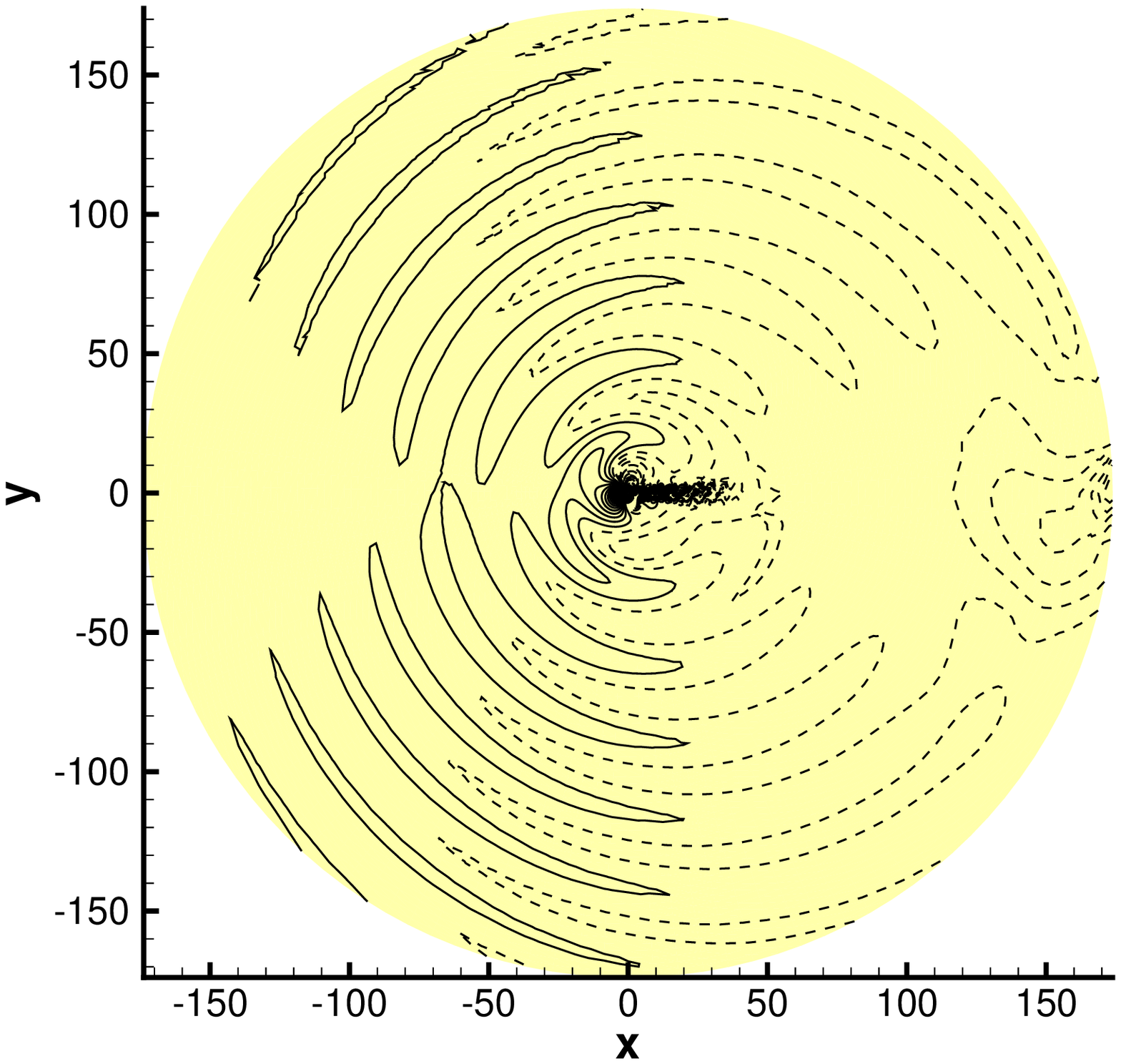}
	\caption{\label{2d-subcylinder-ps} The contours of pressure perturbation $\Delta p$ in the range of $[-0.0285,0.0290]$. The $\Delta p$ is defined as $\Delta p=(p-p_{\infty})/p_{\infty}$. The dash line represents negative value, and the solid line represents positive value}
\end{figure}

The computed time history of the drag and lift coefficients are presented in Fig. \ref{2d-subcylinder-clcd}, in which it can be observed that the flow is almost fully developed after $t=400$. The single-frequency variation of the lift and drag coefficients is obtained. The frequency of lift coefficient can be obtained by taking a Fourier transform of the lift coefficient in the range of $t=[500,1200]$. The frequency of the lift coefficient is $f=0.0362549$.
And from the lift coefficient, we can get the dimensionless frequency, i.e. the Strouhal number
\begin{equation*}
St=f \cdot D/U=0.0362549 \times \frac{1.0}{0.2}=0.181.
\end{equation*}
The experimental value of the Strouhal number $St$ ranges from $0.179 - 0.182$ \cite{cylinder-experiments}.
The average of the drag coefficient can be obtained with the drag coefficient in the range of $t=[500,1200]$. The reference $C_D$ from experiment is $C_D=1.34$ \cite{cylinder-experiments}.
The results from other computations are given in Table \ref{2d-subcylinder-St}.
The contours of Mach number of the fully developed flow are shown in Fig. \ref{2d-subcylinder-Ma}. The vortex street is captured clearly. The pressure perturbation is shown in Fig. \ref{2d-subcylinder-ps}, where the pressure perturbation is defined as $\Delta p=(p-p_{\infty})/p_{\infty}$. The dashed line represents negative value, and the solid line represents positive value.

\begin{table}[!h]
	\small
	\begin{center}
		\def\temptablewidth{0.5\textwidth}
		{\rule{\temptablewidth}{0.50pt}}
		\footnotesize
		\begin{tabular*}{\temptablewidth}{@{\extracolsep{\fill}}c|c|c}
			Scheme 		    				& $St$         		&$C_D$   \\
			\hline
			Compact 4th-order GKS			& 0.181      	 	&1.337   \\
			P3P5 scheme\cite{dumbserpnpm}   & 0.182   	 		&1.33	 \\
			FD-6th\cite{muller-subcylinder} & 0.183   	 		&1.34	 \\
			DG-4th\cite{zhang-DGFV} 		& 0.183				&1.349   \\
		\end{tabular*}
		{\rule{\temptablewidth}{0.50pt}}
	\end{center}
	\vspace{-6mm} \caption{\label{2d-subcylinder-St} Strouhal number and average drag coefficient obtained by different schemes. }
\end{table}

\subsection{Propagation of sound, entropy and vorticity waves}
\begin{figure}[!htb]
	\centering
	\includegraphics[width=0.325\textwidth]{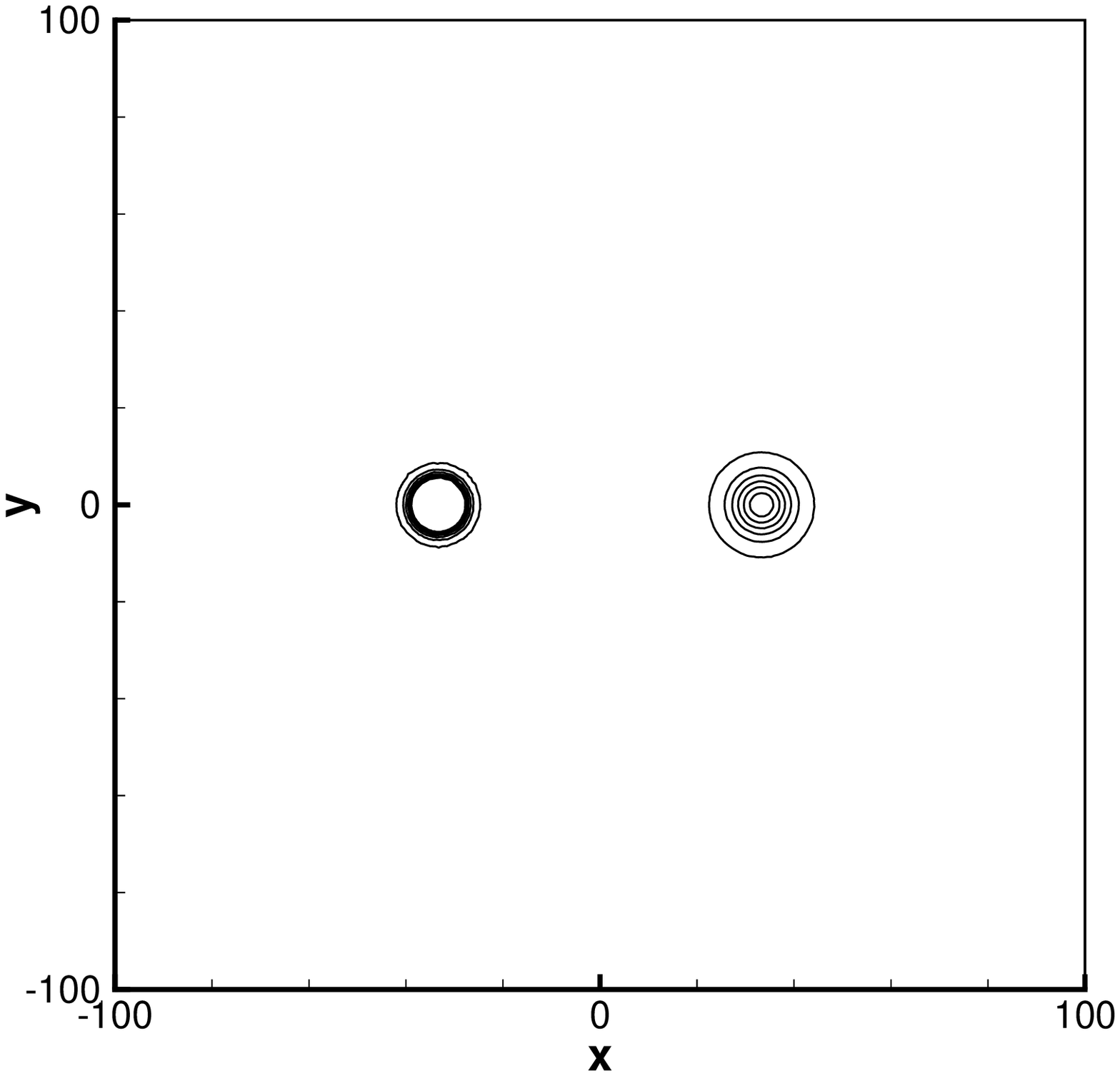}
	\includegraphics[width=0.325\textwidth]{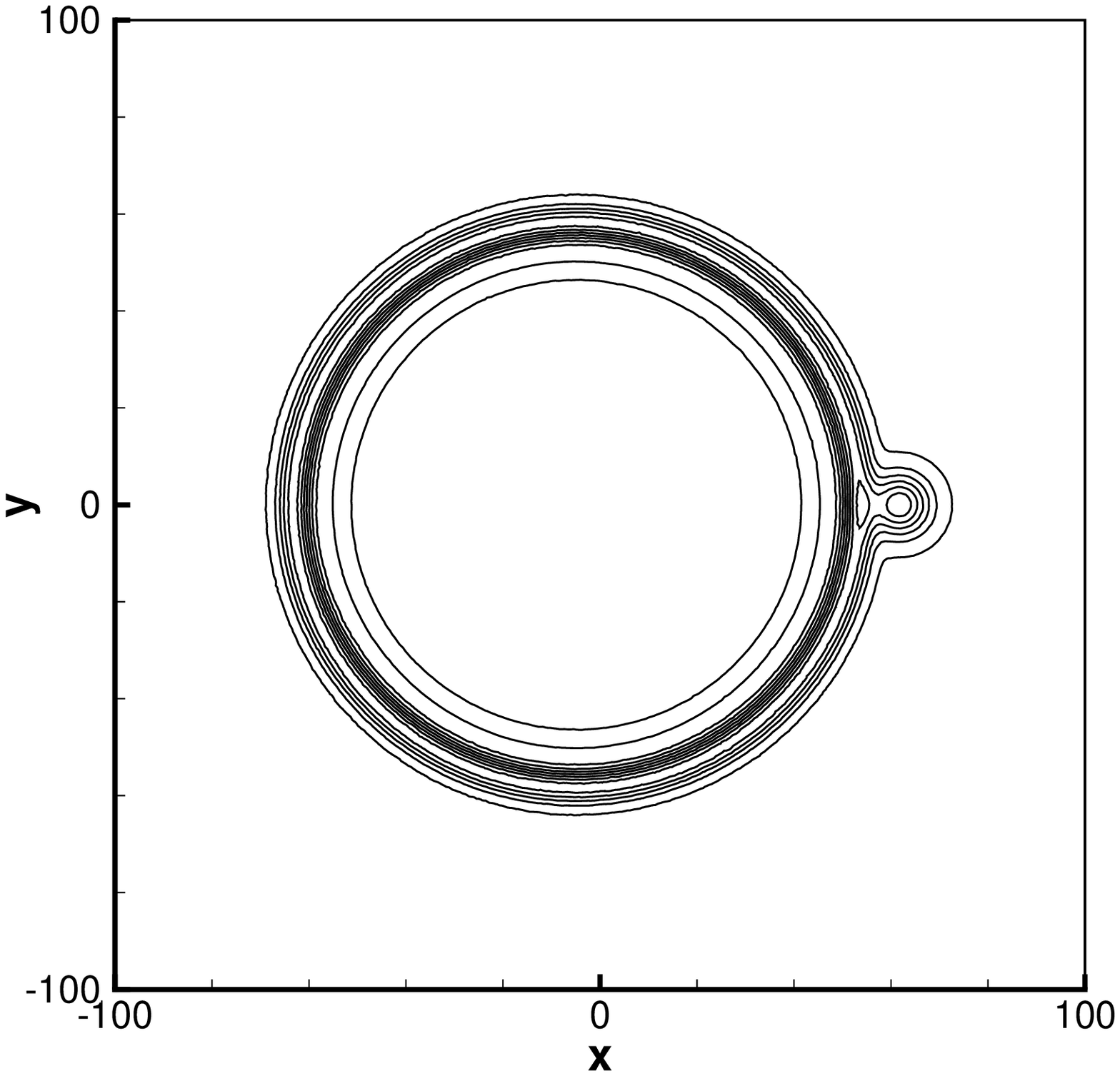}
	\includegraphics[width=0.325\textwidth]{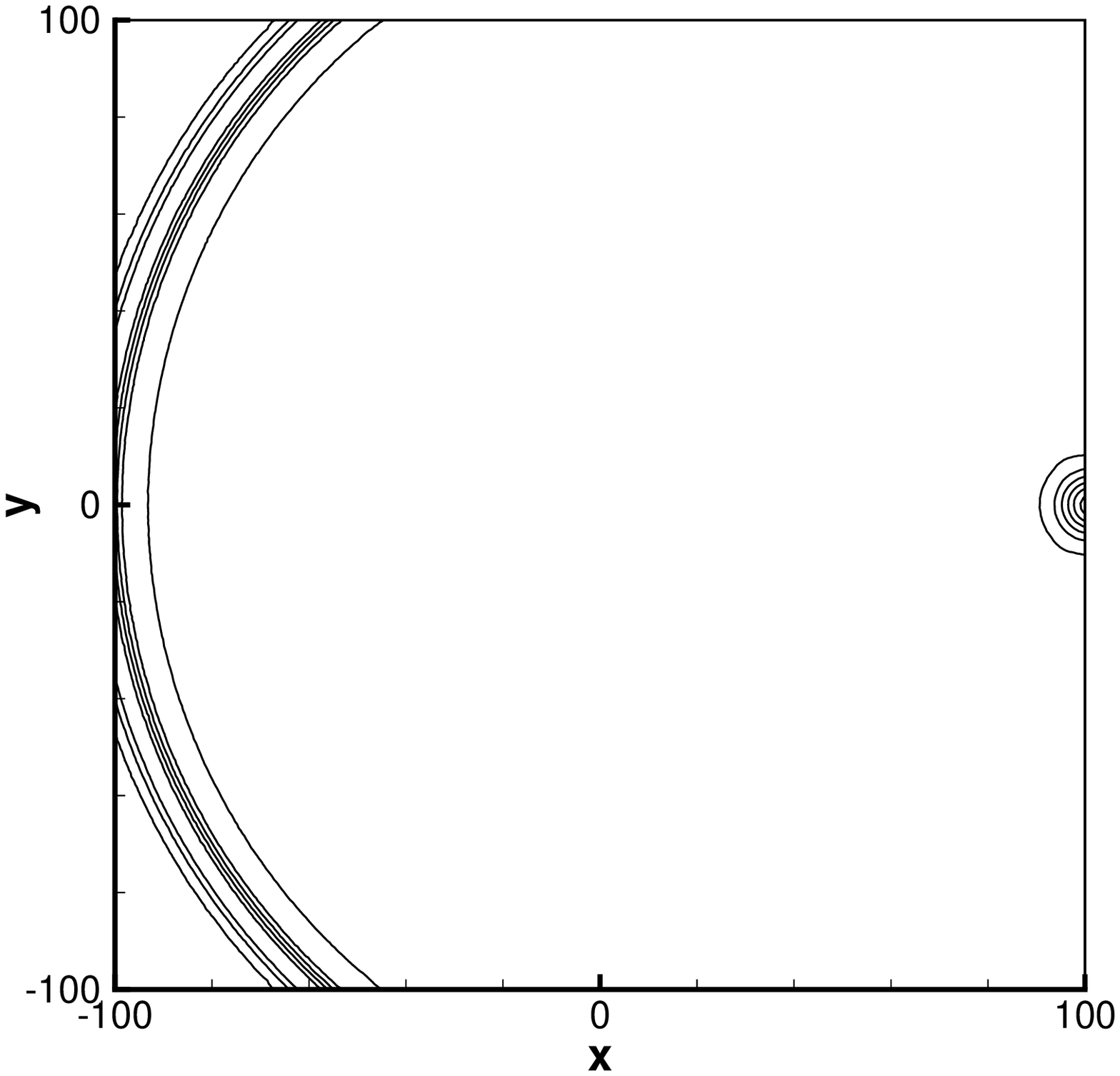}\\
	\includegraphics[width=0.325\textwidth]{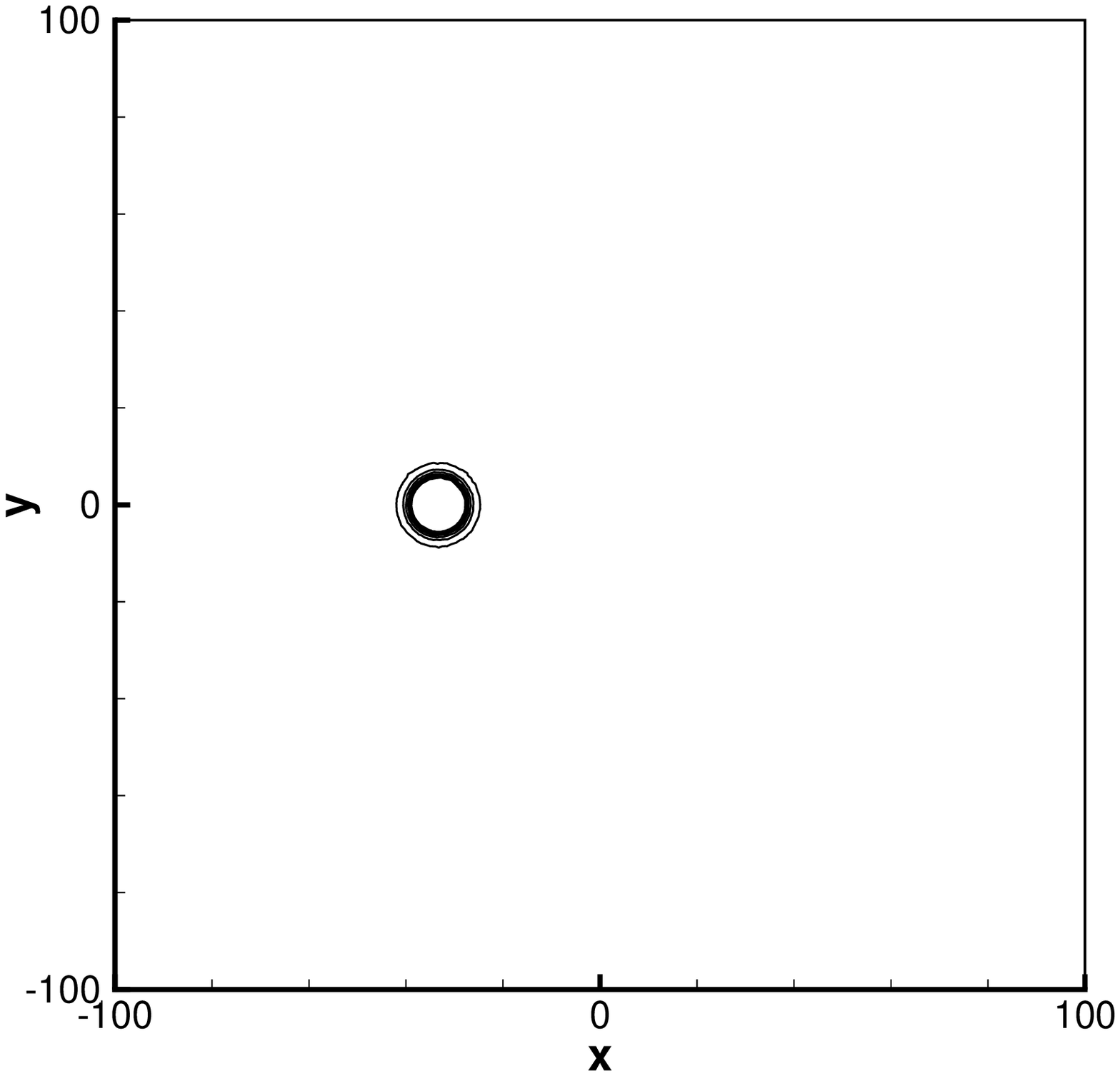}
	\includegraphics[width=0.325\textwidth]{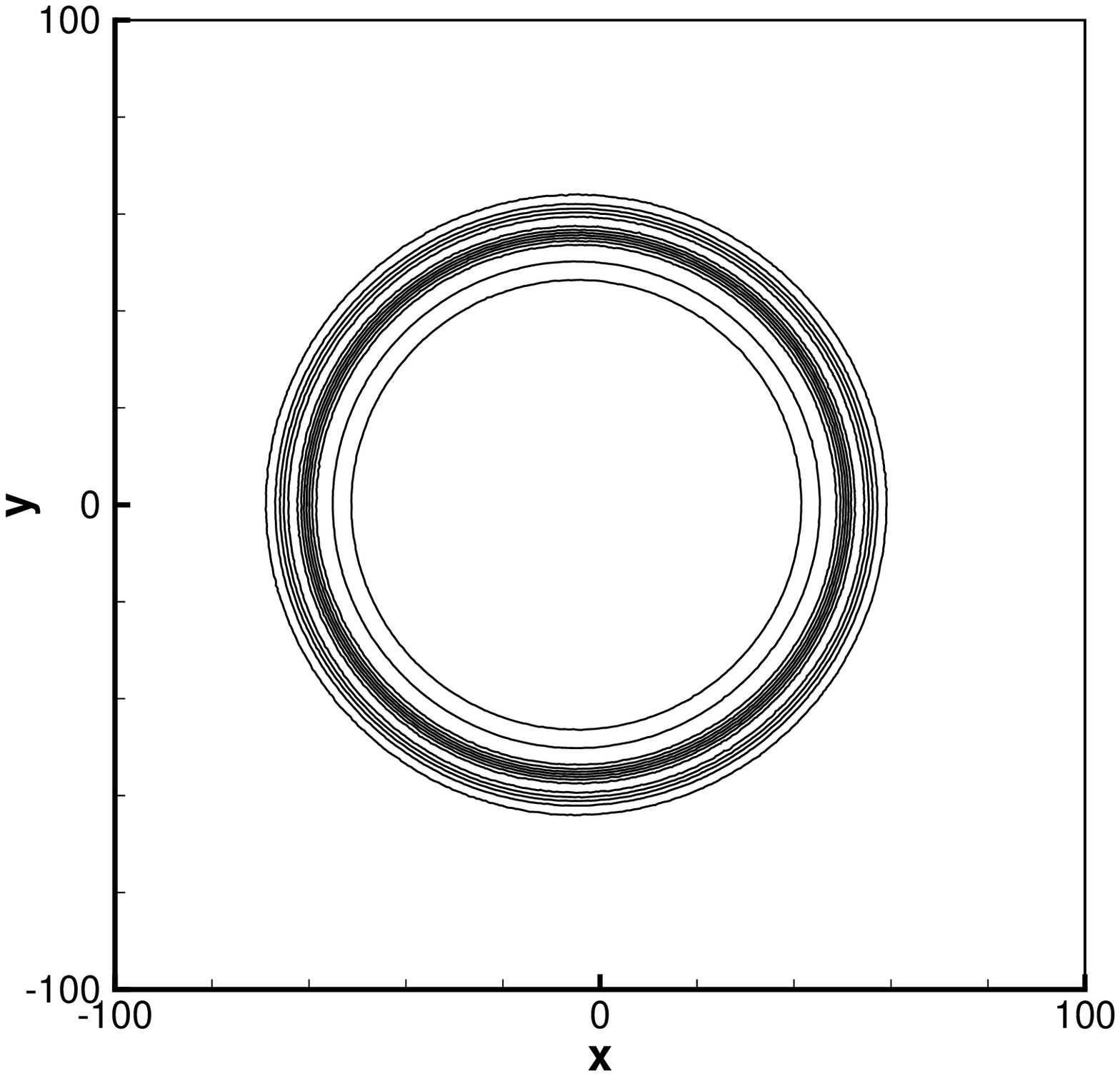}
	\includegraphics[width=0.325\textwidth]{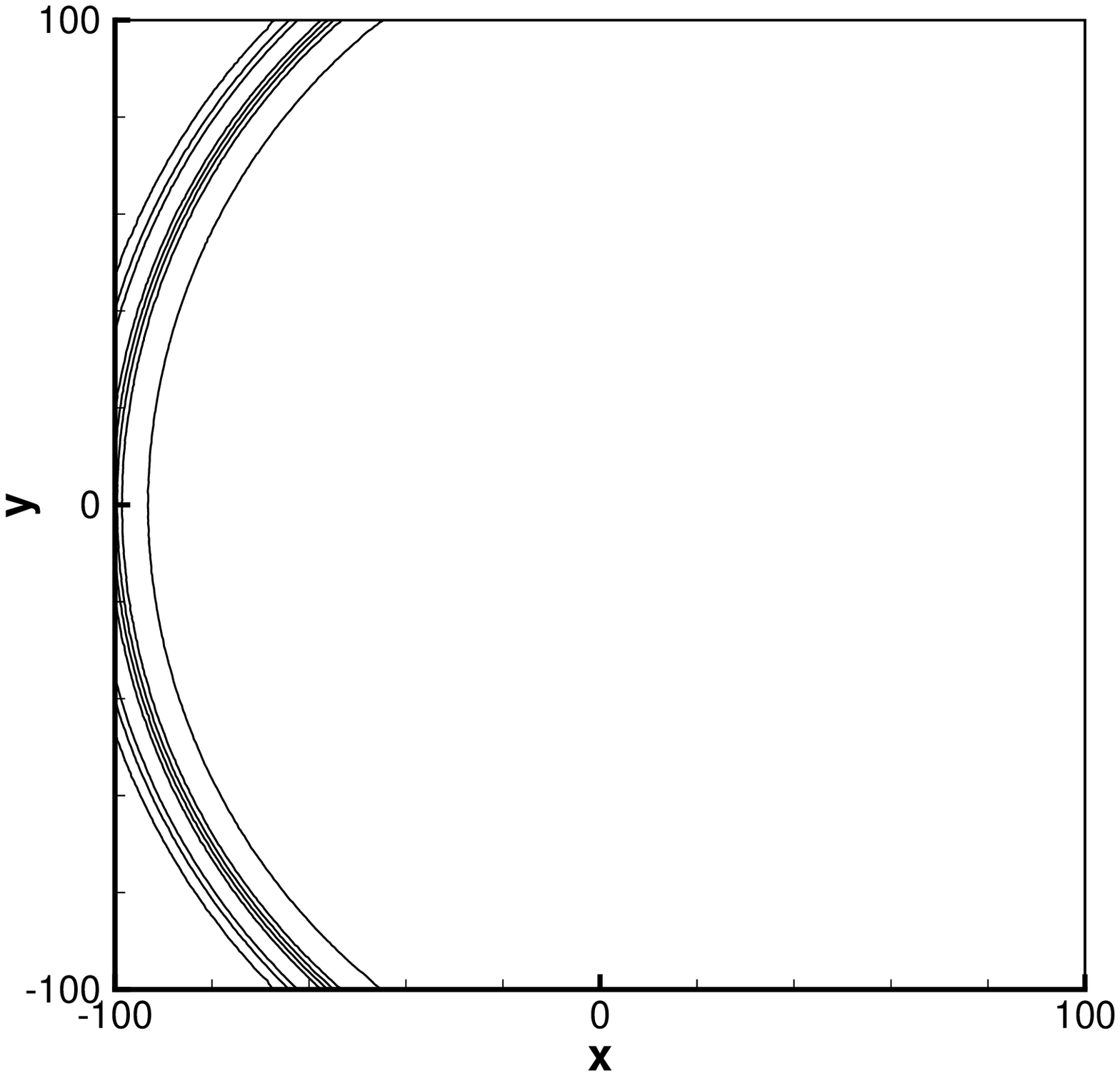}
	\caption{\label{2d-acoustic-1} Sound, entropy and vorticity wave test: density perturbation contours (upper) and pressure perturbation contours (lower) in the range $[0.0006, 0.001]$ with $0.00016$ increment obtained by compact fourth-order GKS at $t=0$, $t = 56.9$ and $t=136$. The irregular triangular mesh with $h=1$ is used. }
\end{figure}

The case was studied by Tam and Webb \cite{tam1993} to demonstrate the high resolution of finite difference schemes.
Initially, an acoustic, entropy, and vorticity pulses are set on a uniform mean flow. The wavefront of acoustic pulse expands radially, and the wave pattern is convected downstream with the mean flow. The entropy and vorticity pulses are convected downstream with the mean flow without any distortion.
The initial condition of the uniform mean flow with three pulses is the following,
\begin{align*}
&\rho= \rho_{\infty} +\epsilon_1 e^{-\alpha_1\cdot r_1^2} + \epsilon_2 e^{-\alpha_2\cdot r_2^2}, \\
& p=p_{\infty} +\epsilon_1 e^{-\alpha_1\cdot r_1^2},  \\
& U=U_{\infty} +\epsilon_3 e^{-\alpha_3\cdot r_3^2}(y-y_3), \\
& V=V_{\infty} -\epsilon_3 e^{-\alpha_3\cdot r_3^2}(x-x_3),
\end{align*}
where $\rho_{\infty}=1.0$, $U_{\infty}=0.5$ and $V_{\infty}=0$. The Mach number is $Ma=0.5$. $\alpha_l$ $ (l=1,2,3)$ is $\alpha_l=\ln 2/b_l^2$ and $b_l$ is the half-width of the Gaussian perturbation. The parameters of these initial pulses are $\epsilon_1=1\times 10^{-2}$, $\epsilon_2=1\times 10^{-3}$, $\epsilon_3=4\times 10^{-4}$, $b_1=3$, $b_2=b_3=5$. $r_l $ $ (l=1,2,3)$ is $r_l=\sqrt{(x-x_l)^2+(y-y_l)^2}$, where $(x_1,y_1)=(-100/3,0)$ and $(x_2,y_2)=(x_3,y_3)=(100/3,0)$. The computational domain is $[-125,125]\times[-125,125]$, where the region  $[-125,125]\times[-125,125]/[-100,100]\times[-100,100]$ is the buffer region. The irregular triangular mesh is used in the current computation. The cell size is $h=1$ in $[-100,100]\times[-100,100]$, and the cell size of the outer boundary of the buffer zone is $10$. The contours of density and pressure perturbations are plotted in Fig. \ref{2d-acoustic-1}. In order to check the result quantitatively, the density and pressure perturbations along $y=0$ are compared with the reference solution from the compact eighth-order GKS on a finer uniform mesh \cite{CGKSAIA}. The current 4th-order compact GKS gets exact result. This test case demonstrates that current GKS on triangular mesh can capture smooth acoustic waves accurately.

\begin{figure}[!htb]
	\centering
	\includegraphics[width=0.425\textwidth]{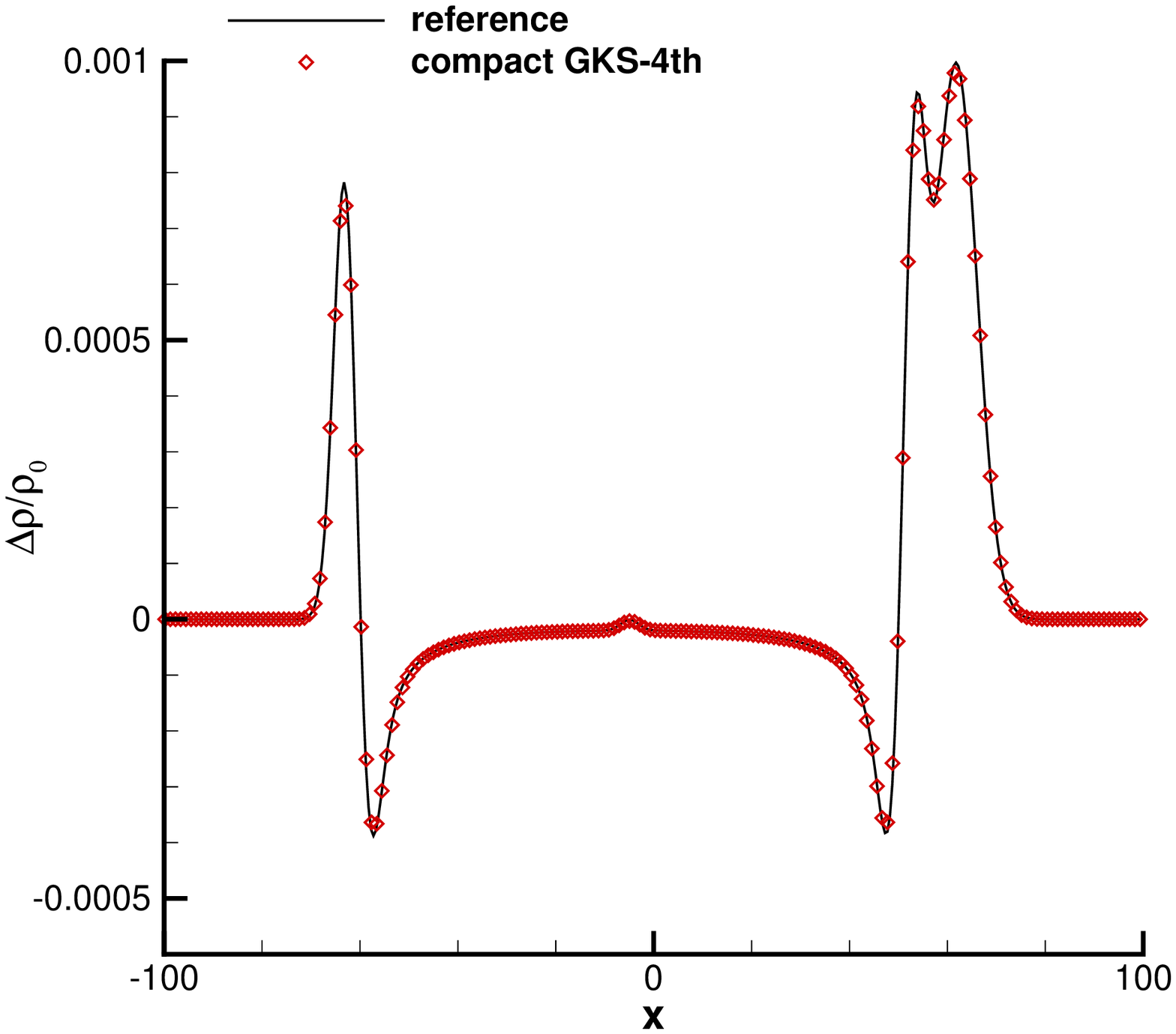}
	\includegraphics[width=0.425\textwidth]{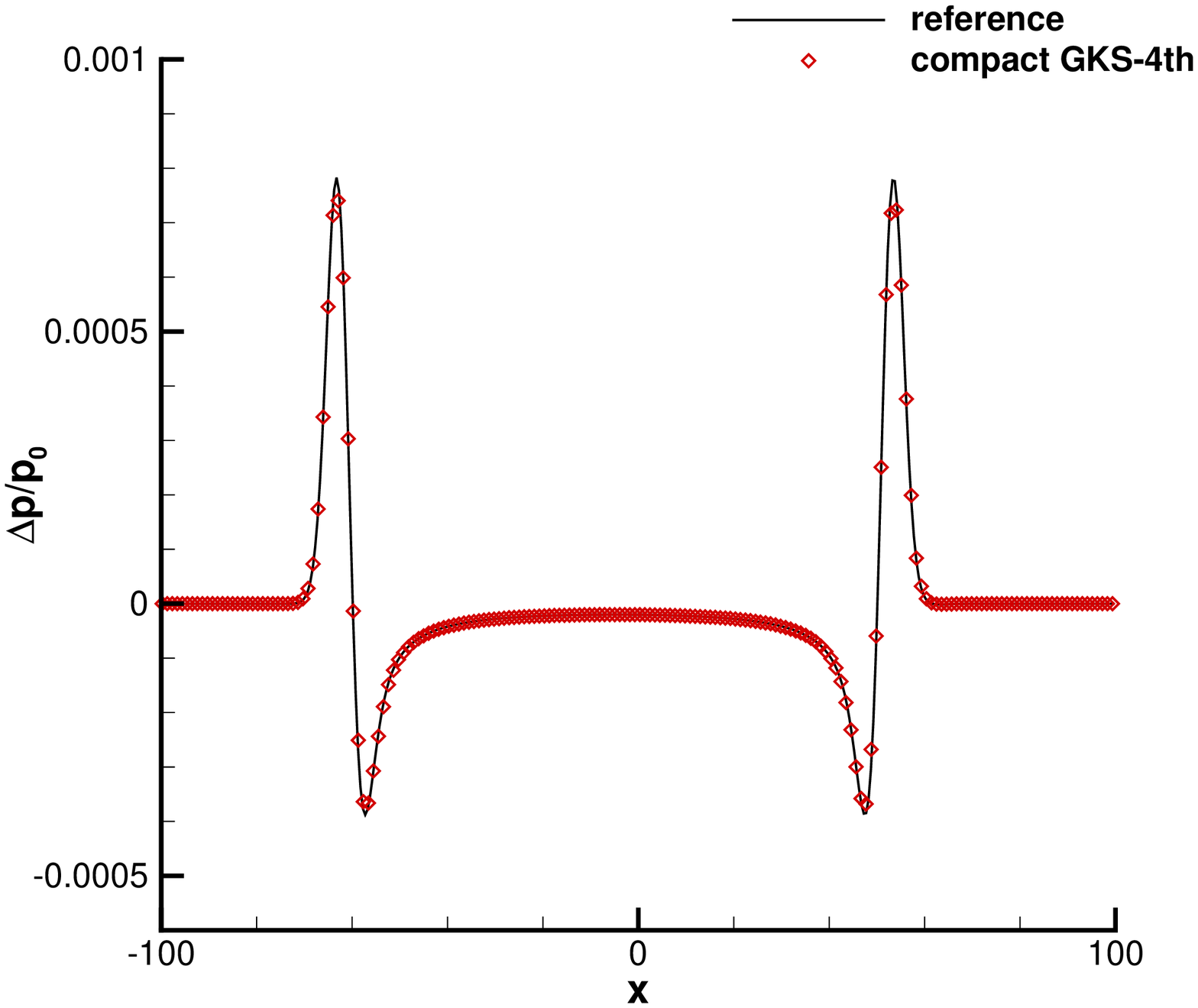}
	\caption{\label{2d-acoustic-2} Sound, entropy and vorticity waves test: density and pressure perturbations along $y=0$ at $t=56.9$ with cell size $h=1$ from compact fourth-order GKS.}
\end{figure}

\subsection{Sound generation by viscous shock-vortex interaction}
The case is the interaction of a shock wave with a single vortex in a viscous flow \cite{Compact-GKSacoustic,zhang-shockvor}.
The computational domain is $[-20,8]\times[-12,12]$.
The velocity of the initial counterclockwise vortex is
\begin{align*}
U_{\theta}(r)&=M_v r e^{(1-r^2)/2},~~U_r=0,
\end{align*}
where $U_{\theta}$ and $U_r$ are the tangential and radial velocity respectively.
The pressure and density distribution superposed by the isentropic vortex downstream of shock wave are
\begin{align*}
p(r)&=\frac{1}{\gamma}[1-\frac{\gamma-1}{2}M_v^2e^{(1-r^2)}]^{\gamma/(\gamma-1)},\\
\rho(r)&=[1-\frac{\gamma-1}{2}M_v^2e^{(1-r^2)}]^{1/(\gamma-1)}.
\end{align*}
The case of the vortex Mach number $M_v=1.0$ is computed. The Mach number of shock wave is $M_s=1.2$.
The Reynolds number is $Re=800$ defined by $Re=\rho_{\infty} L_{\infty} a_\infty /\mu_{\infty}$, where the
subscript $\infty$ denotes the quantity upstream of the shock wave,  $L_{\infty}=1$ is the critical radius of the vortex,
 and $a_\infty =1 $ is the incoming flow sound speed.  The initial location of vortex
is $(x_v,y_v)=(2,0)$, and the stationary shock is at $x=0$.
In the computation, the supersonic inflow boundary condition at $x=8$ as well as the reflecting boundary
condition at $y=\pm 12$ are imposed. The non-reflective boundary conditions are adopted at $x=-20$.
The triangular mesh with a total number of $180570$ cells is used. The mesh is refined near the shock wave with the mesh size $h=0.03$, and the mesh size at boundary is $h=1$.

\begin{figure}[!htb]
\centering
\includegraphics[width=0.495\textwidth]{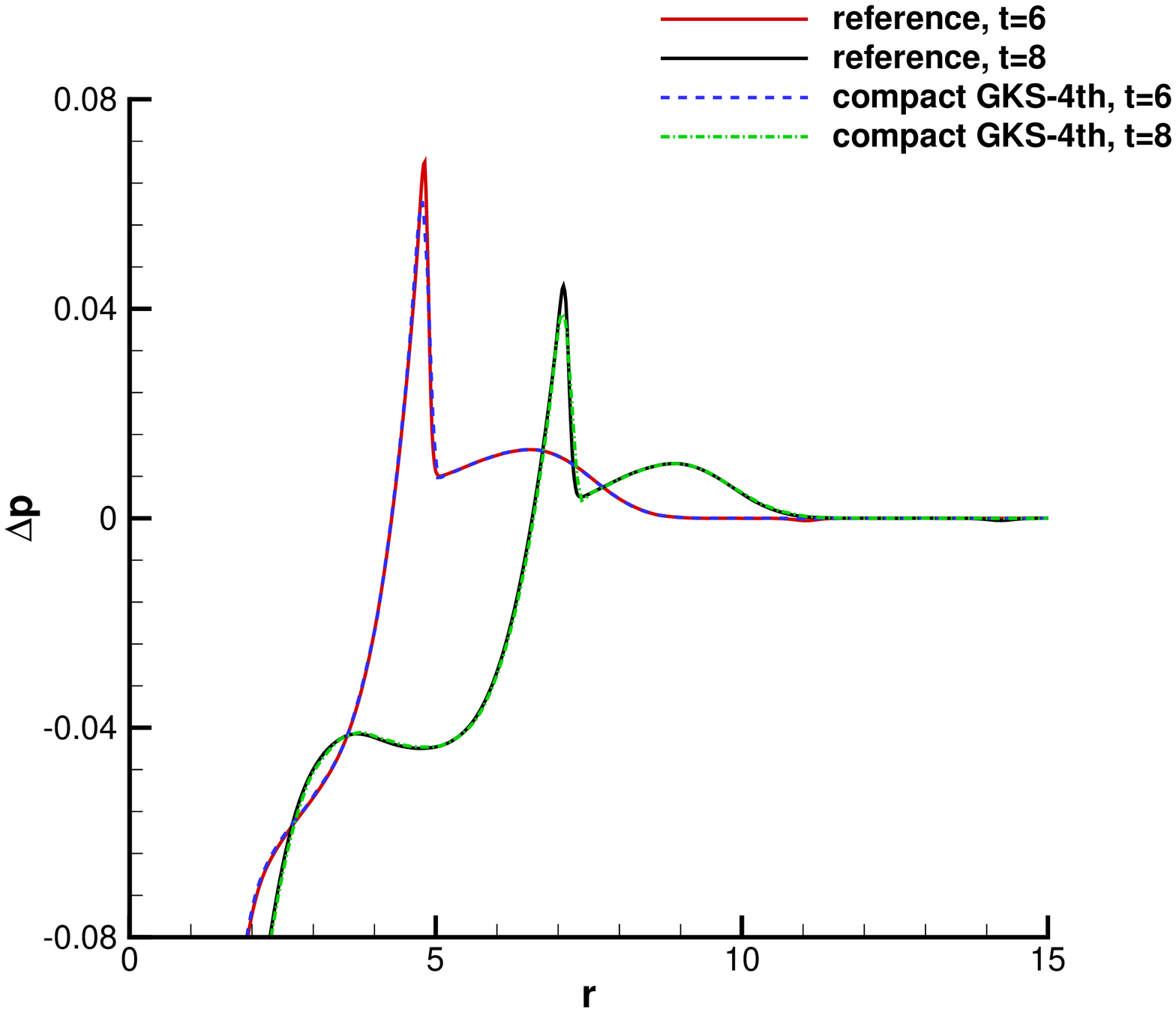}
\includegraphics[width=0.495\textwidth]{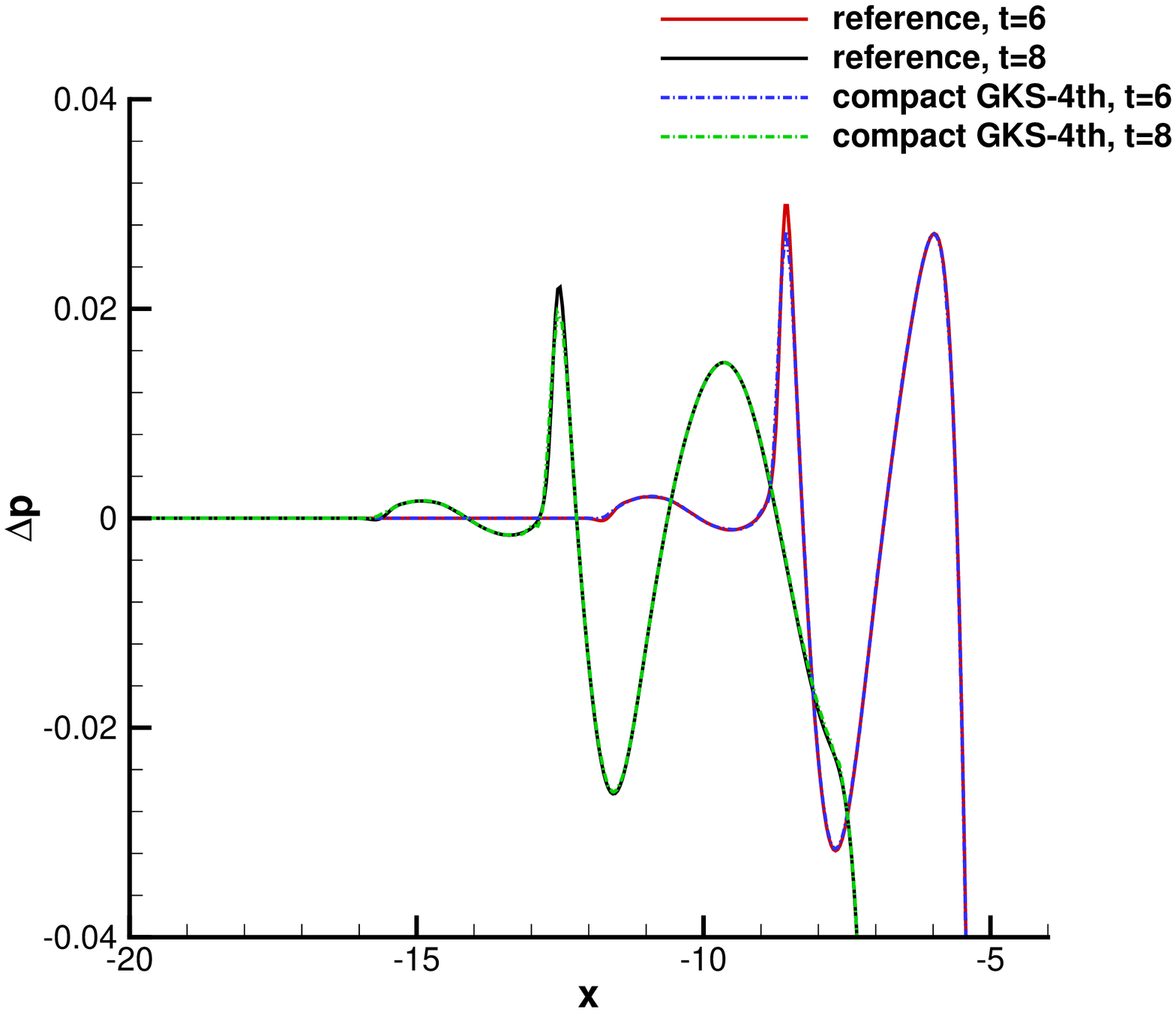}
\caption{\label{vis-shockvor-2} Sound generation by viscous shock-vortex interaction:  pressure perturbation $\triangle p$ along radial distribution direction (left) and the distribution $\triangle p$ along $y=0$ (right). The results of compact 8th-order GKS on the structured mesh are used as  reference solutions \cite{Compact-GKSacoustic}. The results of compact 4th-order GKS are from the triangular mesh. }
\end{figure}

The sound pressure contours at $t=6$ and $t=8$ are given in Fig. \ref{vis-shockvor-1}. The sound pressure is defined as $\triangle p=(p-p_{0})/p_{0}$, where $p_{0}$ is the pressure downstream of the shock wave.
The multiple sound waves with quadruple structure are generated. The sound pressure is smoothly distributed without apparent spurious oscillation.
The reflected shock wave extends and interacts with the vortex at core region.
As a result, complicated flow patterns are formed around the vortex.
Fig. \ref{vis-shockvor-2} is the distribution along radial direction with $135^0$ relative to the positive x-axis, and the distribution along $y=0$ of the sound pressure at different times.
The reflected shock wave produces a pressure jump between the precursor and the second sound wave.
The solutions obtained by the current scheme on triangular mesh are compared with reference solutions \cite{Compact-GKSacoustic}.
In comparison with the computation by high-order WNEO scheme on non-uniform structured mesh in \cite{zhang-shockvor}, the maximum mesh size used in the compact GKS is about $20$ times of that used in the WNEO scheme, and the minimum mesh size is about $10$ times of that in the WNEO scheme.

\begin{figure}[!htb]
\centering
\includegraphics[width=0.495\textwidth]{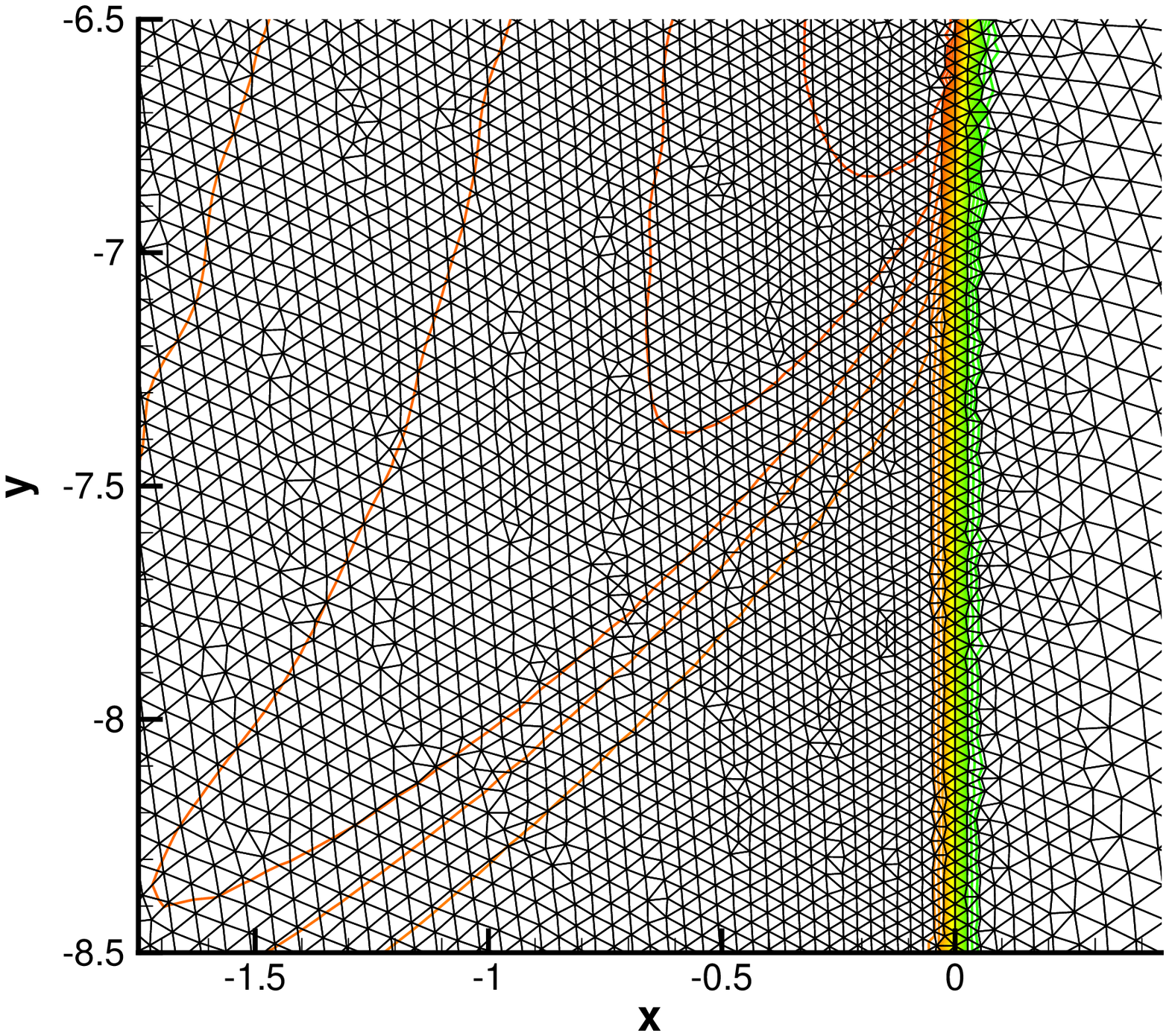}
\includegraphics[width=0.495\textwidth]{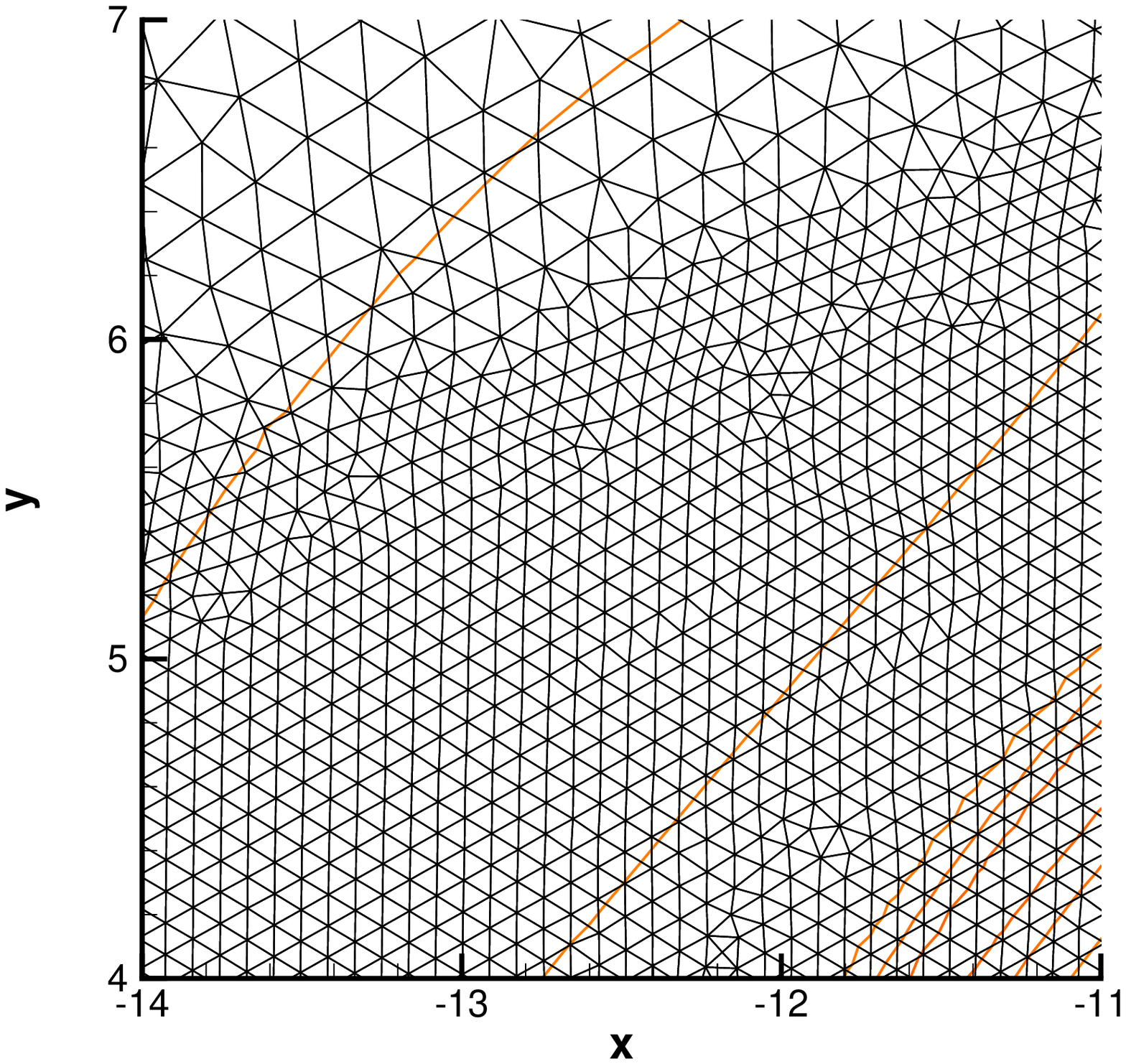}
\caption{\label{vis-shockvor-1} Sound generation by viscous shock-vortex interaction: mesh distribution near the shock wave and the center of vortex at $t=8$. The mesh is refined near the shock wave with the mesh size $h=0.03$, and the mesh size near the precursor wave at $t=8$ is about $h=0.1$. }
\end{figure}

\begin{figure}[!htb]
\centering
\includegraphics[width=0.495\textwidth]{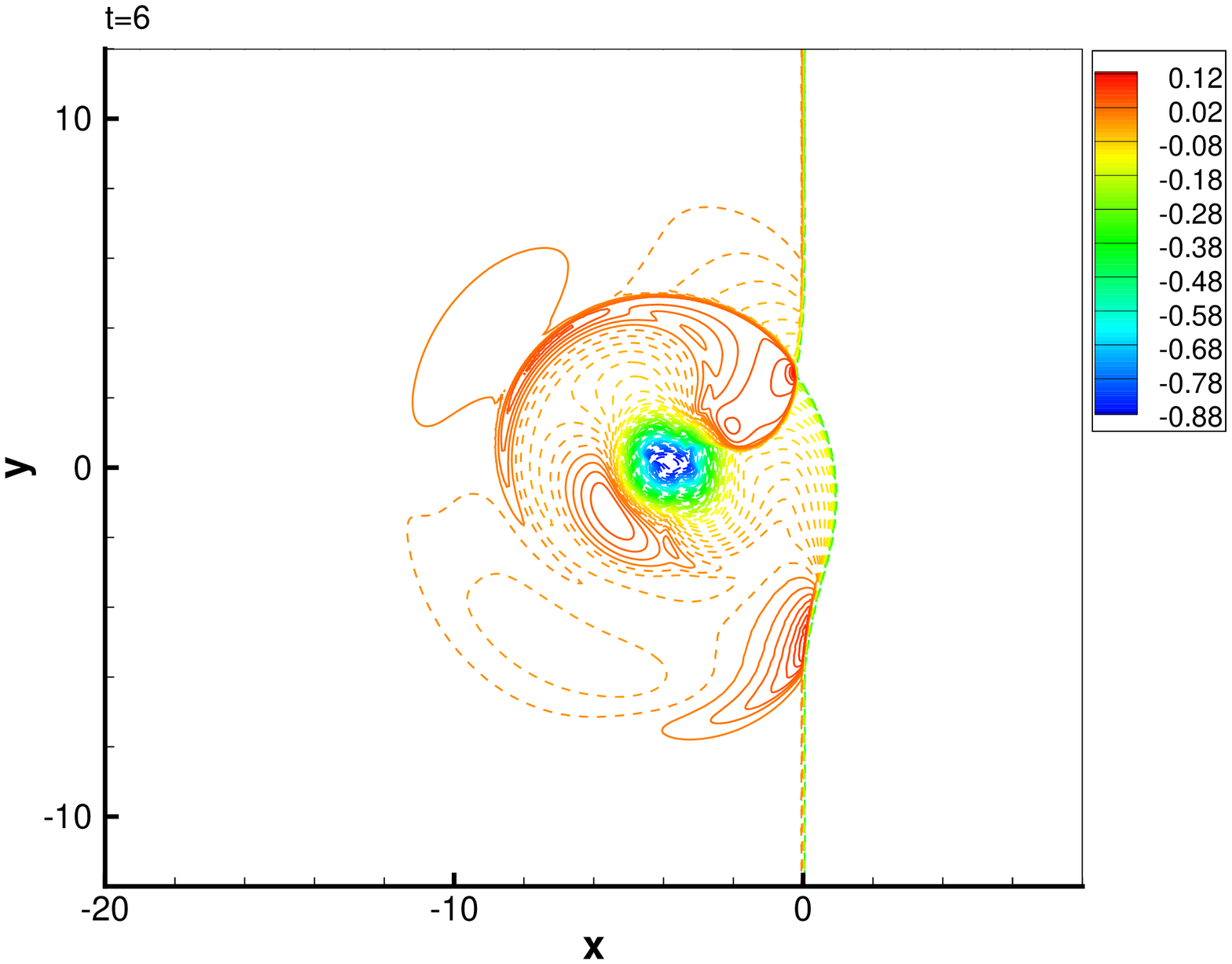}
\includegraphics[width=0.495\textwidth]{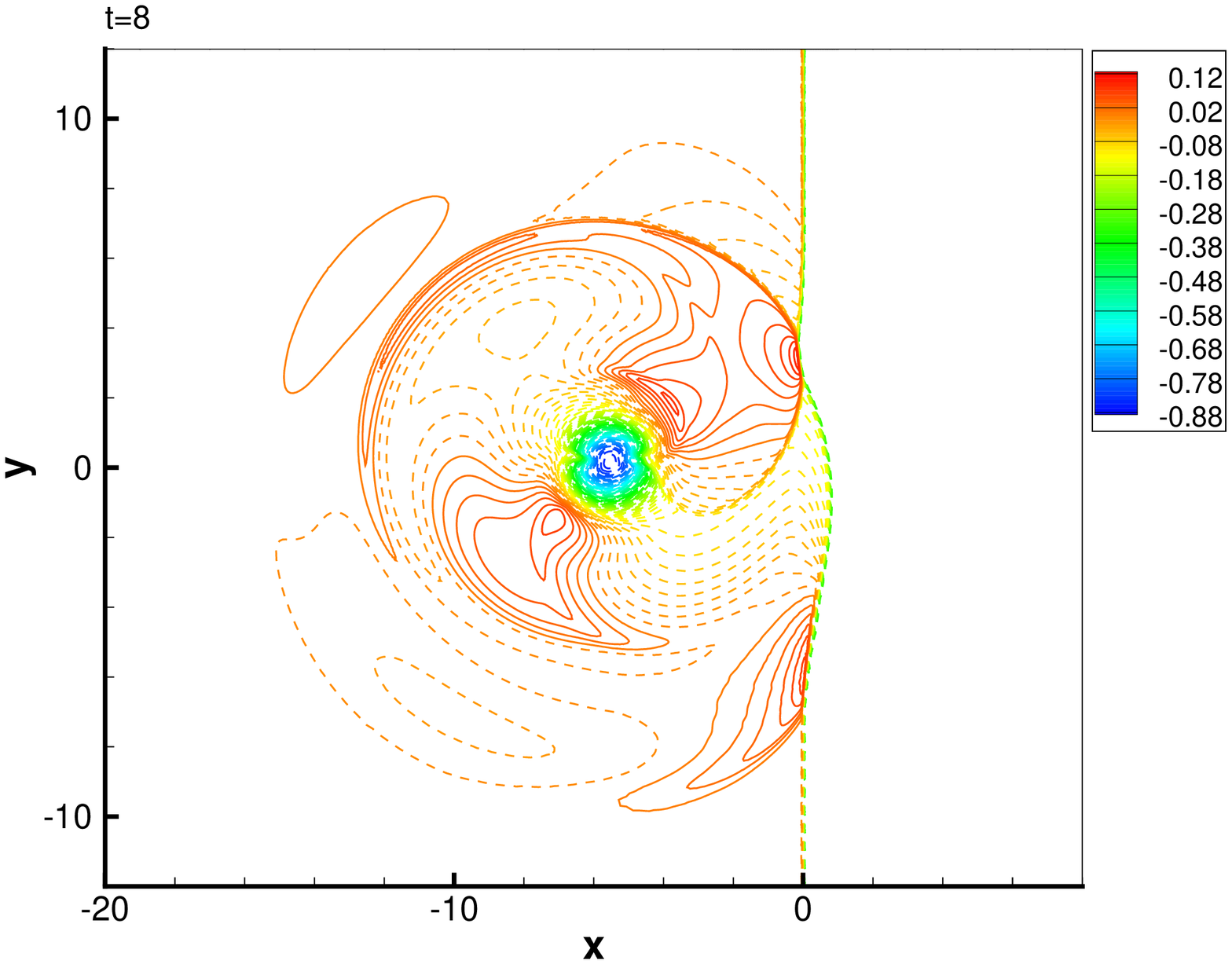}
\caption{\label{vis-shockvor-1} Sound generation by viscous shock-vortex interaction: $M_v=1.0, M_s=1.2$, $80$ equal-spaced sound pressure contours from $\Delta p_{min}=-0.88$ to $\Delta p_{max}=0.12$ on the triangular mesh. The $\triangle p$ is defined as $\triangle p=(p-p_{0})/p_{0}$, and $p_{0}$ is the downstream pressure of the shock wave. The dash line represents rarefaction region, and the solid line represents the compression region. }
\end{figure}

\section{Conclusion}
In this paper, a fourth- to sixth-order compact GKS has been constructed on unstructured mesh for the compressible flow simulation.
For the same order of accuracy, the stencils used in the current scheme is much less than the non-compact counterpart.
The success of the scheme is mainly coming from the high-order gas evolution model at a cell interface, where both time-accurate fluxes and flow variables are provided for the updates of cell-averaged conservative flow variables and their gradients.
Based on the cell averages and gradients, high-order linear and nonlinear reconstructions can be obtained.
The gas-kinetic scheme uses the nonlinear discontinuous reconstruction for the determination of non-equilibrium state
and colliding smooth reconstruction for the equilibrium state, and unifies their contributions
through a relaxation process from non-equilibrium to equilibrium evolution. In the evolution process, the inviscid and viscous terms are merged in
the single time-accurate gas distribution function. This property is preferred on the unstructured mesh for the NS solutions.
Due to the dynamic control of the contributions from the non-equilibrium and equilibrium state from the relaxation factor $\exp(-\Delta t/\tau)$,
both shock wave and acoustic wave solution can be properly captured. In the current scheme, there is no trouble cell indicator or additional limiter used.
Also, benefiting from the time-accurate flux solver,
a fourth-order time accuracy can be obtained through two stages,
instead of four stages in the traditional Runge-Kutta method. Therefore, the current scheme saves two additional reconstructions and its efficiency is no worse than any  traditional high-order DG and WENO scheme.
Many test cases have been used to validate the scheme. Excellent performance has been demonstrated in all cases from the high Mach number shock interaction to the acoustic wave propagation. The current scheme is being extended to 3D calculation on unstructured mesh.

\section*{Acknowledgements}
The current research is supported by National Numerical Windtunnel project and  National Science Foundation of China 11772281, 91852114.

\section*{References}

\bibliographystyle{ieeetr}
\bibliography{AIAbib}

\end{document}